\documentclass[11pt,reqno]{amsproc}
\usepackage[T1]{fontenc}
\usepackage{color}
\usepackage[semicolon,square,authoryear]{natbib}
\numberwithin{equation}{section}
\usepackage{MnSymbol}
\usepackage{stmaryrd}
\usepackage{fullpage}
\usepackage{graphicx}
\usepackage{multirow}
\usepackage{subfigure}
\usepackage{psfrag}
\usepackage{float}
\usepackage{cite}
\usepackage{enumerate}
\newtheorem{theorem}{Theorem}
\newtheorem{property}[theorem]{Property}

\usepackage{morefloats}

\newlength{\drop}
\definecolor{amethyst}{rgb}{0.6, 0.4, 0.8}
\definecolor{burgundy}{rgb}{0.5, 0.0, 0.13}

\title{Do current lattice Boltzmann methods for 
diffusion and advection-diffusion equations respect 
maximum principles and the non-negative constraint?}

\author{\textbf{S.~Karimi} and \textbf{K.~B.~Nakshatrala}\\
{\small Department of Civil and Environmental Engineering, 
University of Houston.}}

\date{\today}

\begin{document}


\begin{titlepage}
    \drop=0.1\textheight
    \centering
    \vspace*{\baselineskip}
    \rule{\textwidth}{1.6pt}\vspace*{-\baselineskip}\vspace*{2pt}
    \rule{\textwidth}{0.4pt}\\[\baselineskip]
    {\LARGE \textbf{\color{burgundy} Do current 
    lattice Boltzmann methods for diffusion and \\[\baselineskip] 
    advection-diffusion equations respect maximum 
    \\[\baselineskip] 
     principles and the non-negative constraint?}}\\[0.3\baselineskip]
    \rule{\textwidth}{0.4pt}\vspace*{-\baselineskip}\vspace{3.2pt}
    \rule{\textwidth}{1.6pt}\\[\baselineskip]
    \scshape
    An e-print of the paper is available on arXiv:~1503.08360 \par
    
    \vspace*{1\baselineskip}
    Authored by \\[\baselineskip]
    
    {\Large S.~Karimi\par}
    {\itshape Graduate Student, University of Houston.}\\[\baselineskip]
    
    {\Large K.~B.~Nakshatrala\par}
    {\itshape Department of Civil \& Environmental Engineering \\
    University of Houston, Houston, Texas 77204--4003. \\ 
    \textbf{phone:} +1-713-743-4418, \textbf{e-mail:} knakshatrala@uh.edu \\
    \textbf{website:} http://www.cive.uh.edu/faculty/nakshatrala\par}
    \vspace*{\baselineskip}
    \vspace{-0.25in}
    \begin{figure}[h]
    \centering
    \includegraphics[scale=0.3,clip]{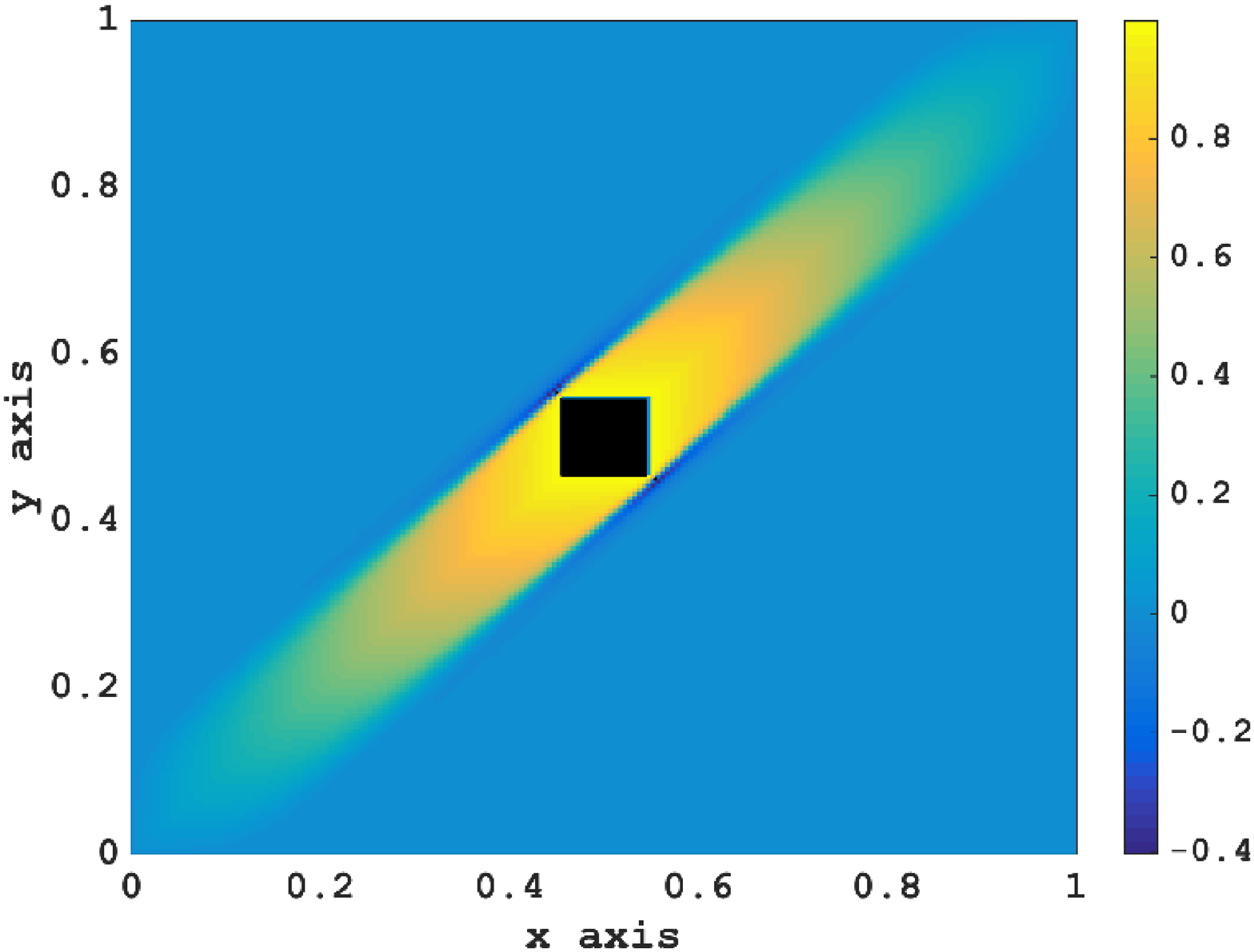}
\end{figure}
\emph{This figure shows that a popular multiple-relaxation-time 
  lattice Boltzmann method violates the non-negative constraint 
  for a transient diffusion equation.}
    \vfill
    {\scshape 2015} \\
    {\small Computational \& Applied Mechanics Laboratory} \par
  \end{titlepage}
	
\begin{abstract}
The lattice Boltzmann method (LBM) has established itself
as a valid numerical method in computational fluid dynamics. 
Recently, multiple-relaxation-time LBM has been proposed 
to simulate anisotropic advection-diffusion processes. The 
governing differential equations of advective-diffusive 
systems are known to satisfy maximum principles, comparison 
principles, the non-negative constraint, and the decay 
property. In this paper, it will be shown that current 
single- and multiple-relaxation-time lattice Boltzmann 
methods \emph{fail} to preserve these mathematical 
properties for transient diffusion and advection-diffusion 
equations. 
It will also be shown that the discretization of Dirichlet 
boundary conditions will affect the performance 
of lattice Boltzmann methods in meeting these 
mathematical principles. A new way of discretizing 
the Dirichlet boundary conditions is also proposed. 
Several benchmark problems have been solved to illustrate 
the performance of lattice Boltzmann methods and the effect 
of discretization of boundary conditions with respect 
to the aforementioned mathematical properties.  
\end{abstract}

\maketitle

\section{INTRODUCTION AND MOTIVATION}
The lattice Boltzmann method (LBM) has gained remarkable 
popularity as a versatile numerical method for fluid 
dynamics simulations \citep{Chen_Doolen_ARFM_1998}. LBM 
has its roots in the kinetic theory as opposed to 
the continuum theory. It needs to be emphasized that 
LBM solves the Boltzmann equation instead of solving 
the continuum field equations. On the other hand, the 
finite element method (FEM) and the finite volume method 
(FVM) solve the continuum field equations directly. 
Some of the attractive features of LBM are: It can 
easily handle irregular domains (e.g., unstructured 
pores and fractures in porous media applications), 
easy to implement even for complicated flow models, 
and natural to parallelize \citep{Succi_LBM}. Great 
advances have been made in extending LBM to simulate 
multi-phase flows \citep{Falcucci_CCP_2006}, reactive 
flows \citep{Rienzo_EPL_2012}, non-linear chemical 
reactions \citep{Ayodele_CCP_2006}, just to name a 
few. In this paper, we limit our scope to LBM-based 
formulations for advection-diffusion phenomena. 

In the recent years, several key advancements have been 
made to extend the LBM to simulate transport phenomena. 
To name a few: \citep{Steibler_Tolke_Krafczyk_CMA_2008}, 
\citep{Shi_Guo_2009_PRE}, \citep{Chai_Zhao_PRE_2013}, 
\citep{Yoshida_Nagaoka_JCP_2010} and \citep{Huang_Wu_JCP_2014}. 
Of these works, Yoshida and Nagaoka \citep{Yoshida_Nagaoka_JCP_2010}, 
and Huang and Wu \citep{Huang_Wu_JCP_2014} have proposed 
multiple-relaxation-time lattice Boltzmann methods to handle 
advection-diffusion equations with \emph{anisotropic} 
diffusivity tensors.

The governing equations for a transient advection-diffusion 
system are parabolic partial differential equations, which 
possess several important mathematical properties. These 
properties include the maximum principle and the comparison 
principle \citep{Protter_Weinberger,Pao}, which have 
crucial implications in modeling physical phenomena. 
For example, a key consequence of the maximum principle 
in modeling advection-diffusion systems is the non-negative 
constraint of the attendant chemical species.  
Several factors such as the physical properties 
of the medium, topology of the domain, and the 
spatial and temporal discretization determine 
the performance of a numerical solution in 
preserving the \emph{discrete} versions of the 
mentioned mathematical properties. A discussion 
on the influence of these factors in the context 
of the finite element method can be found in 
\citep{Nakshatrala_Valocchi_JCP_2009_v228_p6726}. 
Violations of these mathematical properties can 
make a numerical solution inappropriate for
scientific and engineering applications.
It has been shown that many popular finite element 
and finite volume formulations for diffusion-type 
equations violate the maximum principle and the 
non-negative constraint 
\citep{Liska_Shashkov_CiCP_2008_v3_p852,
Nakshatrala_Valocchi_JCP_2009_v228_p6726,
Nakshatrala_Nagarajan_Shabouei_Arxiv_2013}. 
Recently, numerical methodologies have been 
proposed under the finite element method to 
satisfy the non-negative constraint and the 
maximum principle by utilizing the underlying 
variational structure. Since the lattice 
Boltzmann method does not enjoy such a 
variational basis, these methodologies 
developed for the FEM cannot be extended 
to the lattice Boltzmann method.	

\emph{To the best of our knowledge, the performance of 
LBM-based formulations for unsteady diffusion-type 
equations with respect to comparison principles, 
maximum principles, and the non-negative constraint 
has not been documented in the literature. But, 
such a study is of paramount importance, as LBM-based 
formulations are being employed in predictive numerical 
simulations.} We shall, therefore, put some popular 
LBM-based formulations for diffusion-type equations to 
test, and particularly show that these formulations 
violate all the 
aforementioned mathematical principles and physical 
constraints. The LBM-based formulations of interest 
in this paper are the non-thermal single-relaxation-time 
LBM for advection with isotropic diffusion, and 
the multiple-relaxation-time methods for anisotropic 
diffusion proposed in \citep{Yoshida_Nagaoka_JCP_2010} 
and \citep{Huang_Wu_JCP_2014}. The mentioned LBM-based 
formulations have been chosen merely due to their 
popularity, and this paper does not pretend to be 
exhaustive.  

The remainder of this paper is organized as follows. 
Section \ref{Sec:LBM_S2_GE} presents the governing 
equations for a transient advective-diffusive system 
along with a brief discussion on the associated 
mathematical properties. Section \ref{Sec:LBM_S3_LBM} 
reviews the non-thermal single- and multiple-relaxation-time 
lattice Boltzmann methods for transient diffusion-type 
equations. 
In Section \ref{Sec:LBM_S4_NR}, several representative 
test problems are presented with a thorough discussion 
on the numerical results from LBM-based formulations. 
Section \ref{Sec:LBM_S5_Theoretical} performs a 
theoretical analysis on the lattice Boltzmann method 
to obtain a simple criterion on the time-step and 
lattice cell size to satisfy the non-negative constraint. 
Finally, conclusions are drawn in Section \ref{Sec:LBM_S6_CR}.
On the notational front, a quantity in the continuous 
setting will be denoted by upright font symbols (e.g., 
$\mathrm{u}$), and a quantity in the discrete setting 
will be denoted by italic font (e.g., $u$).

\section{UNSTEADY ANISOTROPIC DIFFUSION-TYPE EQUATIONS}
\label{Sec:LBM_S2_GE}
Consider a bounded open domain $\Omega$. We shall 
denote the boundary of the domain by $\Gamma = 
\overline{\Omega} - \Omega$, where $\overline{\Omega}$ 
is the set closure of $\Omega$. We assume that 
the boundary $\Gamma$ comprises of two parts 
$\Gamma^{\mathrm{N}}$and $\Gamma^{\mathrm{D}}$ such 
that $\Gamma^{\mathrm{N}} \cap \Gamma^{\mathrm{D}} 
= \emptyset$ and $\Gamma = \Gamma^{\mathrm{N}}\cup 
\Gamma^{\mathrm{D}}$. We denote the part of 
the boundary on which Dirichlet boundary condition 
is prescribed by $\Gamma^{\mathrm{D}}$. Neumann boundary 
condition is prescribed on $\Gamma^{\mathrm{N}}$. 
A spatial point will be denoted by $\mathbf{x}$. The 
unit outward normal to the boundary is denoted by 
$\widehat{\mathbf{n}}(\mathbf{x})$.
The time interval of interest will be denoted by 
$\left[0,\mathcal{T} \right]$, and the time is 
denoted by $\mathrm{t}$. For convenience, we 
introduce the following differential operator:   
\begin{align}
  \label{Eqn:LBM_differential_operator}
  \mathcal{L} \left[\mathrm{u}(\mathbf{x},\mathrm{t})\right] 
  := \frac{\partial \mathrm{u}\left( \mathbf{x}, \mathrm{t}\right)}{\partial 
    \mathrm{t}} + \mathrm{div} \left[ \mathbf{v}\left( \mathbf{x}, \mathrm{t}\right) 
    \mathrm{u}\left( \mathbf{x}, \mathrm{t}\right) - \mathbf{D}\left( \mathbf{x}\right) 
  \mathrm{grad} \left[ \mathrm{u}\left(\mathbf{x}, \mathrm{t}\right)\right] \right]
\end{align}
where $\mathrm{u}\left(\mathbf{x}, \mathrm{t}\right)$ 
is the concentration field. The diffusivity tensor is denoted 
by $\mathbf{D}\left(\mathbf{x}\right)$, which is assumed 
to be symmetric, positive definite, and bounded above. The 
advection velocity is denoted by $\mathbf{v}(\mathbf{x},
\mathrm{t})$, which is assumed to be solenoidal. That is, 
$\mathrm{div}[\mathbf{v}] = 0$. The divergence and gradient 
operators with respect to $\mathbf{x}$ are, 
respectively, denoted by $\mathrm{div}[\cdot]$ 
and $\mathrm{grad}[\cdot]$. The initial boundary value problem 
for a transient advection-diffusion system can be written as 
follows: 
\begin{subequations}
  \label{Eqn:Diffusion_Equation}
  \begin{align}
    \label{Eqn:Diffusion_Equation_Eq}
    &\mathcal{L}\left[ \mathrm{u}\left(\mathbf{x}, \mathrm{t}\right)\right] 
    = \mathrm{g}\left( \mathbf{x}, \mathrm{t}\right) \quad &&\left(\mathbf{x}, 
    \mathrm{t} \right)\in \Omega \times \left( 0 , \mathcal{T} \right] \\
      &\mathrm{u}\left( \mathbf{x},\mathrm{t}\right) = 
      \mathrm{u}^{\mathrm{p}}\left( \mathbf{x},\mathrm{t}\right) 
      \quad &&\left(\mathbf{x}, \mathrm{t} \right) \in 
      \Gamma^{\mathrm{D}}\times \left[ 0, \mathcal{T}\right]\\
      &\left(\mathbf{v}\left( \mathbf{x}, \mathrm{t}\right) 
      \mathrm{u}\left( \mathbf{x}, \mathrm{t}\right) - \mathbf{D} 
      \left(\mathbf{x}\right) \mathrm{grad} \left[ \mathrm{u}\left(\mathbf{x}, \mathrm{t}\right) \right]
      \right) \cdot \widehat{\mathbf{n}}\left(\mathbf{x}\right) 
      = \mathrm{q}^{\mathrm{p}}\left( \mathbf{x},\mathrm{t}\right) 
      \quad &&\left(\mathbf{x}, \mathrm{t} \right) \in \Gamma^{\mathrm{N}} 
      \times \left[ 0, \mathcal{T}\right] \\
      \label{Eqn:Diffusion_Equation_IC}
      &\mathrm{u}\left(\mathbf{x}, \mathrm{t} = 0\right) 
      = \mathrm{u}_{0} \left(\mathbf{x}\right) \quad 
      &&\mathbf{x} \in \Omega
  \end{align}
\end{subequations}
where the source/sink is denoted by $\mathrm{g}(\mathbf{x},
\mathrm{t})$, $\mathrm{u}^{\mathrm{p}}(\mathbf{x},\mathrm{t})$ 
is the prescribed concentration, $\mathrm{q}^{\mathrm{p}}
(\mathbf{x},\mathrm{t})$ is the prescribed flux, and 
the initial concentration is denoted by $\mathrm{u}_{0} 
(\mathbf{x})$. 
It is easy to check that equation \eqref{Eqn:Diffusion_Equation_Eq} 
is a linear parabolic partial differential equation. The initial 
boundary value problem given by equations 
\eqref{Eqn:Diffusion_Equation_Eq}--\eqref{Eqn:Diffusion_Equation_IC} 
satisfies several important mathematical properties, which will be 
discussed next. 

\subsection{Mathematical properties}
We introduce the following function space:
\begin{align}
  \mathsf{C}^{2}_{1} \left( \Omega \times 
  \left( 0, \mathcal{T}\right] \right) := \left\{ 
  \mathrm{u} : \Omega \times \left[ 0, \mathcal{T}\right] 
  \rightarrow \mathbb{R} \; \big| \; 
  \mathrm{u}, 
  \frac{\partial \mathrm{u}}{\partial \mathrm{t}}, 
  \frac{\partial \mathrm{u}}{\partial \mathrm{x}_i}, 
  \frac{\partial^2 \mathrm{u}}{\partial \mathrm{x}_i 
    \partial \mathrm{x}_j}
  \in \mathsf{C} \left( \Omega \times 
  \left( 0, \mathcal{T}\right]\right)\right\}
\end{align}
where $\mathsf{C}(\Omega \times \left(0,\mathcal{T}
\right])$ is the set of all continuous functions 
defined on $\Omega \times \left(0,\mathcal{T}\right]$. 
Similarly, one can define $\mathsf{C}(\overline{\Omega} 
\times [0,\mathcal{T}])$.

  \begin{property}[The maximum principle]
    \label{Thm:Max_Principle}
    Let $\mathrm{u}(\mathbf{x},\mathrm{t})\in 
    \mathsf{C}^{2}_{1}\left(\Omega \times 
    \left( 0, \mathcal{T}\right] \right) \cap 
      \mathsf{C}\left( \overline{\Omega} \times 
      \left[ 0, \mathcal{T}\right]\right)$ be 
      a solution of the initial boundary value 
      problem \eqref{Eqn:Diffusion_Equation} 
      with $\partial \Omega = \Gamma^{\mathrm{D}}$.
      If $\mathrm{g}(\mathbf{x},\mathrm{t}) \geq 0$
      then
      \begin{align}
        \min_{\left( \mathbf{x}, \mathrm{t}\right) \in \overline{\Omega} 
          \times \left[ 0, \mathcal{T}\right]} \mathrm{u} \left( \mathbf{x}, 
        \mathrm{t}\right) = \min \left[
          \min_{\left( \mathbf{x}, \mathrm{t}\right) \in \Gamma \times 
            \left[ 0 , \mathcal{T}\right]} \mathrm{u} (\mathbf{x}, 
          \mathrm{t}), \; \min_{\mathbf{x} \in \Omega} \; 
          \mathrm{u}_0(\mathbf{x}) \right]
      \end{align}
\end{property}

\begin{property}[The comparison principle]
\label{Thm:Comp_Principle}
	Let $\mathrm{u}_1 \left( \mathbf{x}, \mathrm{t}\right)$ and
	$\mathrm{u}_2 \left( \mathbf{x}, \mathrm{t}\right)$ belong 
	to $\mathsf{C}^{2}_{1} 	\left( \Omega \times \left(0,\mathcal{T}\right]\right) \cap 
	\mathsf{C} \left( \overline{\Omega} \times \left[ 0, \mathcal{T}\right]\right)$. If $\mathcal{L} 
	\left[ \mathrm{u}_1 \left( \mathbf{x}, \mathrm{t}\right)\right] 
	\geq \mathcal{L} \left[ \mathrm{u}_2 \left( \mathbf{x}, 
	\mathrm{t}\right)\right]$ on $\Omega \times \left( 0 , 
	\mathcal{T}\right]$, and $\mathrm{u}_1^{\mathrm{p}} \left( \mathbf{x}, 
	\mathrm{t}\right) \geq \mathrm{u}_2^{\mathrm{p}} \left( \mathbf{x}, 
	\mathrm{t}\right)$ on $\Gamma \times \left[ 0, \mathcal{T}
	\right]$ then $\mathrm{u}_1 \left( \mathbf{x}, \mathrm{t}\right) 
	\geq \mathrm{u}_2 \left( \mathbf{x}, \mathrm{t}\right)$ on 
	$\overline{\Omega}\times \left[ 0 , \mathcal{T}\right]$. 
\end{property}
Mathematical proofs to the  maximum principle and the 
comparison principle can found in \citep{Evans_PDE}. 
We now show that the non-negative constraint for the 
concentration can be obtained as a consequence of 
the maximum principle under certain assumptions on 
the input data. One could alternatively obtain the 
non-negative constraint from the comparison principle,
which we would not present here.

\begin{property}[The non-negative constraint]
  If $\mathrm{g}(\mathbf{x},\mathrm{t}) \geq 0$ 
  in $\Omega$, $\mathrm{u}^{\mathrm{p}}(\mathbf{x},t) 
  \geq 0$ on $\Gamma$, and $\mathrm{u}_0(\mathbf{x}) 
  \geq 0$ in $\Omega$ then 
  \begin{align}
    \mathrm{u}(\mathbf{x},t) \geq 0 \quad 
    \mathrm{in} \; \overline{\Omega} \; 
    \mathrm{and} \; \forall t
  \end{align}
\end{property}
\begin{proof}
  Based on the hypothesis, we have the 
  following two results:
  \begin{align}
    &\min_{(\mathbf{x},\mathrm{t}) \in \Gamma 
      \times [0,\mathcal{T}]} \mathrm{u} 
    (\mathbf{x},\mathrm{t}) = 
    \min_{(\mathbf{x},\mathrm{t}) \in \Gamma 
      \times [0,\mathcal{T}]} \mathrm{u}^{\mathrm{p}} 
    (\mathbf{x},\mathrm{t}) \geq 0 \\
    &\min_{\mathbf{x} \in \Omega} \; \mathrm{u}_0 
    (\mathbf{x}) \geq 0 
  \end{align}
  These further imply that 
  \begin{align}
    \min \left[
      \min_{\left( \mathbf{x}, \mathrm{t}\right) \in \Gamma \times 
        \left[ 0 , \mathcal{T}\right]} \mathrm{u} \left( \mathbf{x}, 
      \mathrm{t}\right), \; \min_{\mathrm{x} \in \Omega} \; 
      \mathrm{u}_0(\mathbf{x}) \right]
    \geq 0 
  \end{align}
  From the maximum principle, we can conclude that 
  \begin{align}
    \min_{(\mathbf{x}, \mathrm{t}) \in \overline{\Omega} 
      \times [0,\mathcal{T}]} \mathrm{u} 
    (\mathbf{x},\mathrm{t}) \geq 0 
  \end{align}
  which gives the desired result that 
  $\mathrm{u}(\mathbf{x},\mathrm{t}) 
  \geq 0$ in $\overline{\Omega}$ and 
  $\forall t$. 
\end{proof}
The following integrals will be used 
in the remainder of this paper: 
\begin{align}
  \label{Eqn:Integralu}
  \mathcal{J}_{1}(\mathrm{u};\Omega;\mathrm{t}) 
  := \int_{\Omega} \mathrm{u}(\mathbf{x},\mathrm{t}) 
  \; \mathrm{d}\Omega, \; \;
  \mathcal{J}_{1}^{+}(\mathrm{u};\Omega;\mathrm{t}) 
  &:= \int_{\Omega} \max \left[\mathrm{u}(\mathbf{x},
    \mathrm{t}), 0 \right] \; \mathrm{d}\Omega, \nonumber \\
  \mathcal{J}_{2}(\mathrm{u};\Omega;\mathrm{t}) 
  := \int_{\Omega} \mathrm{u}^2(\mathbf{x},\mathrm{t}) 
  \; \mathrm{d}\Omega, \; \;
  \mathcal{J}_{2}^{+}(\mathrm{u};\Omega;\mathrm{t}) 
  &:= \int_{\Omega} \left( \max[\mathrm{u}(\mathbf{x},\mathrm{t}),0] 
  \right) ^2 \; \mathrm{d}\Omega
\end{align}

\begin{property}[The decay property]
  If $\mathbf{v}(\mathbf{x},t) = \mathbf{0}$, 
  $\mathrm{u}^{\mathrm{p}}(\mathbf{x},\mathrm{t}) 
  = 0$ on the entire $\Gamma$, and $\mathrm{g}
  (\mathbf{x},\mathrm{t}) = 0$ in $\Omega$ then 
  \begin{align}
    \frac{\mathrm{d}}{\mathrm{d}\mathrm{t}}\mathcal{J}_{2}
    (\mathrm{u};\Omega;\mathrm{t}) \leq 0 \quad \forall 
    \mathrm{t}
  \end{align}
\end{property}
\begin{proof} 
  Noting that $\mathrm{g}(\mathbf{x},\mathrm{t}) = 0$ 
  and $\mathbf{v}(\mathbf{x},\mathrm{t}) = \mathbf{0}$, 
  equation \eqref{Eqn:Diffusion_Equation_Eq} implies 
  the following:
  \begin{align}
    \int_{\Omega} \mathrm{u} \frac{\partial \mathrm{u}}{\partial \mathrm{t}} 
    \; \mathrm{d} \Omega 
    - \int_{\Omega} \mathrm{u} \; \mathrm{div}[\mathbf{D}(\mathbf{x})
      \mathrm{grad}[\mathrm{u}]] \; \mathrm{d} \Omega 
    = 0
  \end{align}
  Using Green's identity and noting that 
  $\mathrm{u}(\mathbf{x},\mathrm{t}) = 
  \mathrm{u}^{\mathrm{p}}(\mathbf{x},\mathrm{t}) 
  = 0$ on $\Gamma$, we obtain the following:
  \begin{align}
    \int_{\Omega} \mathrm{u} \frac{\partial \mathrm{u}}{\partial \mathrm{t}} 
    \; \mathrm{d} \Omega 
    + \int_{\Omega} \mathrm{grad}[\mathrm{u}] \cdot  \mathbf{D}(\mathbf{x})
    \mathrm{grad}[\mathrm{u}] \; \mathrm{d} \Omega 
    = 0
  \end{align}
  Noting that $\mathbf{D}(\mathbf{x})$ is a positive 
  definite second-order tensor, we have:
  \begin{align}
    \int_{\Omega} \mathrm{u} \frac{\partial \mathrm{u}}{\partial \mathrm{t}} 
    \; \mathrm{d} \Omega \leq 0 
  \end{align}
  This further implies that 
  \begin{align}
    \frac{\mathrm{d}}{\mathrm{d}\mathrm{t}} 
    \int_{\Omega} \mathrm{u}^{2}(\mathbf{x},\mathrm{t}) 
    \; \mathrm{d} \Omega 
    \equiv \frac{\mathrm{d}}{\mathrm{d}\mathrm{t}} 
    \mathcal{J}_2(\mathrm{u};\Omega;\mathrm{t}) 
    \leq 0 
  \end{align}
  In order to obtain the above result, we 
  have assumed that $\Omega$ is independent 
  of $\mathrm{t}$, which is the case in this 
  paper.
\end{proof}
In the subsequent sections we will illustrate the 
performance of some popular LBM-based formulations 
with respect to the aforementioned mathematical 
properties in the discrete setting. We will also 
compare the performance of the lattice Boltzmann 
method with the finite element method in this 
regard.

\section{THE LATTICE BOLTZMANN METHOD}
\label{Sec:LBM_S3_LBM}
The lattice Boltzmann method is a numerical method for 
solving the Boltzmann equation. A numerical 
simulation based on LBM will provide an 
approximate distribution of particles in the discrete 
configuration-momentum space. Macroscopic quantities 
such as concentration can then 
be calculated using the obtained discrete distributions. 
In this paper, we will consider the Bhatnagar-Gross-Krook 
(BGK) model for the collision term:
\begin{align}
\label{Eqn:BoltzmannEq}
	\frac{\partial \mathrm{f}}{\partial \mathrm{t}} + 
	\frac{1}{\mathrm{m}} \mathbf{p} \cdot \mathrm{grad}
	\left[ \mathrm{f} \right] + \mathbf{F} \cdot 
	\frac{\partial \mathrm{f}}{\partial \mathbf{p}} = 
	\frac{1}{\lambda} \left(  
	\mathrm{f}^{\mathrm{eq}} - \mathrm{f}\right)
\end{align}
where the continuous distribution function is denoted 
by $\mathrm{f} = \mathrm{f}(\mathbf{p},\mathbf{x},\mathrm{t})$, 
the equilibrium distribution is denoted by 
$\mathrm{f}^{\mathrm{eq}}$, 
$\lambda$ is the relaxation time, the macroscopic momentum 
is denoted by $\mathbf{p}$, and $\mathbf{F}$ is the external 
force. We shall use the symbol $\tau$ for the non-dimensional 
relaxation time.
We will now specialize on the LBM-based formulations for 
isotropic and anisotropic advection-diffusion equations. 
However, for an in-depth discussion on the lattice Boltzmann 
method, one can consult the references \citep{He_Luo_PRE_1997,Succi_LBM}. 
We first present the single-relaxation-time lattice Boltzmann 
method, and then followed by multiple-relaxation-time methods. 

\subsection{\textbf{The single-relaxation-time lattice Boltzmann method}}
\begin{figure}
  \centering
  \psfrag{0}{0}
  \psfrag{1}{1}
  \psfrag{2}{2}
  \psfrag{3}{3}
  \psfrag{4}{4}
  \psfrag{5}{5}
  \psfrag{6}{6}
  \psfrag{7}{7}
  \psfrag{8}{8}
  \psfrag{ln}{lattice node}
  \psfrag{1d}{$D1Q3$}
  \psfrag{2d}{$D2Q9$}
  \psfrag{2d2}{$D2Q5$}
  \includegraphics[scale=0.7,clip]{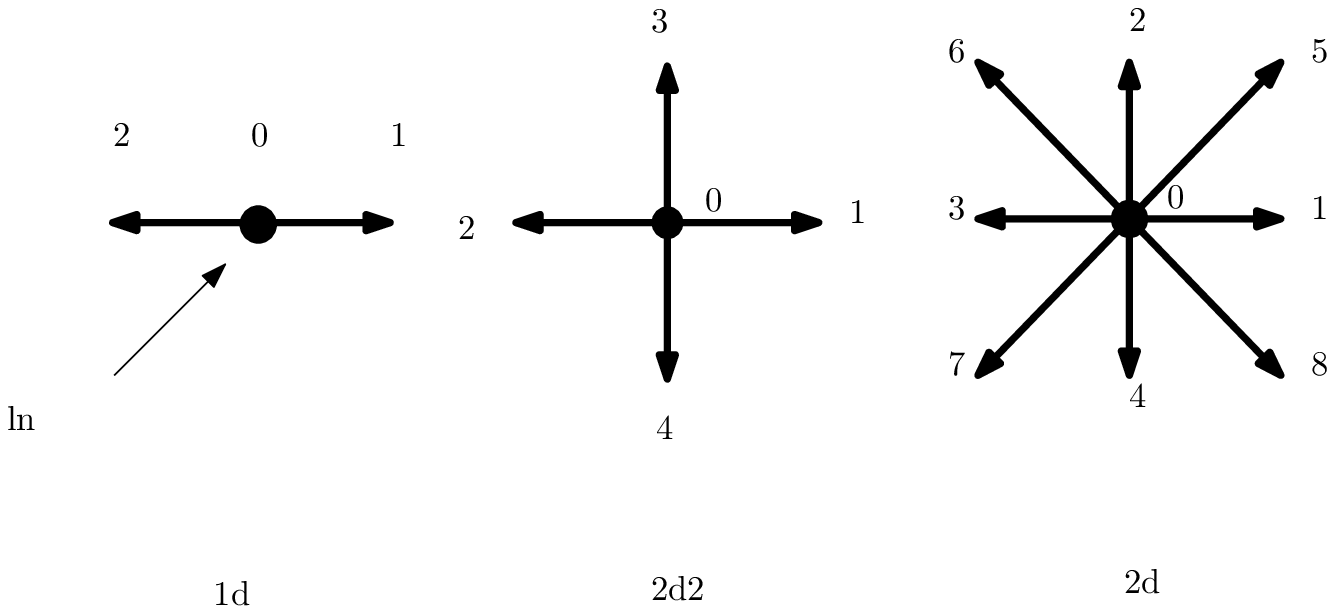}
  \caption{This figure shows the one- and two-dimensional 
    lattice models that have been employed in this paper. 
    \label{Fig:LatticeModels}}
\end{figure}

We will use a uniform and structured discretization of the spatial
domain $\Omega$. The lattice cell size is denoted by $\Delta x$, 
and $\Delta t$ is the time-step. We will consider a $DnQm$ lattice 
model, where $n$ is the number of spatial dimensions and $m$ 
is the number of discretized directions for momentum. 
The lattice models that are employed in this paper 
are shown in Figure \ref{Fig:LatticeModels}.
The distribution of particles at a lattice node 
$\boldsymbol{x}$, at time $t$ and along the 
velocity direction $\boldsymbol{e}_i$ will be 
denoted by $f_i(\boldsymbol{x},t)$. 
The evolution of the discretized distributions 
is governed by the following formula: 
\begin{align}
\label{Eqn:Translation}
	f_i \left( \boldsymbol{x} + \boldsymbol{e}_i \Delta t, t + 
	\Delta t\right) = \widehat{f}_i \left( \boldsymbol{x}, t\right)
	\qquad \text{[translation step]}
	\end{align}
where the post-collision distribution is denoted 
by $\widehat{f}_i$ and is defined as follows:
\begin{align}
\label{Eqn:Collision}
	\widehat{f}_i\left( \boldsymbol{x}, t\right) =  
	f_i \left( \boldsymbol{x}, t\right) - \frac{1}{\tau} \left( 
	f_i \left( \boldsymbol{x}, t\right) - f_i^{\mathrm{eq}} 
	\left( \boldsymbol{x}, t\right)\right) + 
	w_i \Delta t g\left( \boldsymbol{x}, t\right) \qquad 
	\text{[collision step]}
\end{align} 
where $g\left( \boldsymbol{x}, t\right)$ is the source/sink 
at a lattice node. The weight associated with the $i$-th 
velocity direction is denoted by $w_i$. The equilibrium 
distribution for the $i$-th velocity 
direction will be denoted by $f_i^{\mathrm{eq}}$. A popular choice 
for the equilibrium distribution of an advective-diffusive system 
takes the following form:
\begin{align}
	f_i^{\mathrm{eq}} \left( \boldsymbol{x}, t\right) = w_i 
	u \left( \boldsymbol{x}, t\right) \left( 1 + \frac{\mathbf{v} 
	\left( \boldsymbol{x},	t\right) \cdot \boldsymbol{e}_i}
	{c_s^2}\right)
\end{align}
where advection velocity is denoted by $\mathbf{v}$, and 
$c_s$ is the lattice (sound) velocity. The relaxation time 
is given by the following relation:
\begin{align}
\label{Eqn:relaxation_SRT}
	\tau = \frac{\mathrm{D}}{\Delta t c_s^2} + \frac{1}{2}
\end{align}
where $\mathrm{D}$ is the (isotropic) diffusion coefficient. 
The concentration at a lattice node can then be calculated 
as:
\begin{align}
\label{Eqn:Distribution_Concnetration}
	u \left( \boldsymbol{x}, t\right) = \sum_{i = 1}^{m} 
	f_i \left( \boldsymbol{x}, t\right)
\end{align} 
We now describe the discretization of
boundary conditions. 

\begin{enumerate}[(i)]
\item \emph{Standard method for Dirichlet boundary conditions}: 
Let $\{j\}$ be the set of unknown distributions at a point 
$\boldsymbol{x} \in \Gamma^{\mathrm{D}}$, and $\{i\}$ the set of 
known distributions at that point. We can then write
\begin{align}
\label{Eqn:Dirichlet_BC_Std}
	f_{\alpha} \left( \boldsymbol{x}, t \right) = 
	\frac{w_{\alpha}}{\sum_{p \in \{ j\}} w_p} \left( \mathrm{u}^{\mathrm{p}} 
	\left( \boldsymbol{x},t \right) - \sum_{q \in \{ i\}} f_q 
	\left( \boldsymbol{x}, t  \right)\right)
	\quad \forall \alpha \in \{ j \}
\end{align}
It will be shown that this standard way of discretizing 
Dirichlet boundary conditions can contribute to the 
violation of the non-negative constraint. We, therefore, 
propose a new way of discretizing the Dirichlet boundary 
conditions.
\item \emph{Weighted splitting method for Dirichlet boundary conditions}:
Let $\boldsymbol{x}$ be the spatial
coordinates of a lattice node on $\Gamma^{\mathrm{D}}$. 
The discretized distributions at each time-level will be 
assigned as follows:
\begin{align}
\label{Eqn:Dirichlet_BC_ad_hoc}
	f_{\alpha} \left( \boldsymbol{x}, t\right) = w_{\alpha} 
	\mathrm{u}^{\mathrm{p}}\left( \boldsymbol{x}, t\right) 
	\quad \alpha = 0, \cdots, m-1 
\end{align}
where $\mathrm{u}^{\mathrm{p}}$ is the 
prescribed concentration.
\item \emph{Neumann boundary conditions}: Let $\{j\}$ be 
the set of unknown distributions at a point $\boldsymbol{x} 
\in \Gamma^{\mathrm{N}}$. We can then discrete the 
Neumann boundary conditions as follows:
\begin{align}
\label{Eqn:Neumann_BC}
	f_{\alpha} \left( \boldsymbol{x}, t + \Delta t\right) = 
		\widehat{f}_{\beta} \left( \boldsymbol{x}, t\right) + 
		\frac{1}{c_s} \mathrm{q}^{\mathrm{p}} 
		\left( \boldsymbol{x}, t + \Delta t\right)
		\frac{\boldsymbol{e}_{\alpha} \cdot \widehat{\mathbf{n}}\left( \boldsymbol{x}\right)}{\sum_{k \in {j}} \boldsymbol{e}_{k} \cdot \widehat{\mathbf{n}}\left( \boldsymbol{x}\right)}
		 \quad
		\forall \alpha \in \{ j \}
\end{align}
where $\beta$ is the index of a discrete velocity
direction such that $\boldsymbol{e}_{\alpha} = - 
\boldsymbol{e}_{\beta}$, and $\mathrm{q}^{\mathrm{p}}$ 
is the prescribed flux. Note that $\boldsymbol{e}_{\alpha} 
\cdot \widehat{\mathbf{n}} \left( \boldsymbol{x}\right) 
< 0 \; \text{for} \;\alpha \in \left\{ j\right\}$.
\end{enumerate}

\subsection{\textbf{The multiple-relaxation-time 
    lattice Boltzmann methods}}
The single-relaxation-time lattice Boltzmann methods are 
limited to isotropic diffusion. Multiple-relaxation-time 
lattice Boltzmann methods offer a more suitable framework 
for simulating anisotropic diffusion, and are more easily 
presentable in the momentum space 
\citep{Lallemand_Lou_PRE_v61_2000,Li_He_Tang_Tao_PRE_v81_2010}. 
The distributions can be transformed to moments 
using the following linear operation:
\begin{align}
  | \varrho_i \rangle \left( \boldsymbol{x}, t \right) = 
  \boldsymbol{M} | f_i\rangle \left( \boldsymbol{x}, t\right)
\end{align}
where $\boldsymbol{M}$ is an $m \times m$ orthogonal transformation
matrix, and $| f_i \rangle \left( \boldsymbol{x}, t\right)$ is
the column vector of size $m \times 1$ of the distributions at
lattice node $\boldsymbol{x}$ and at time $t$. 
The moment corresponding to $f_i$ is denoted by 
$\varrho_i$. The BGK lattice Boltzmann equation 
can then be written as follows:
\begin{align}
  | \varrho_i\rangle \left(\boldsymbol{x} + \boldsymbol{e}_i \Delta t, 
  t + \Delta t\right)
  = | \varrho_i\rangle \left( \boldsymbol{x}, t \right) -
  \boldsymbol{S} \left( | \varrho_i\rangle \left( \boldsymbol{x}, t \right)
  - | \varrho_i^{\mathrm{eq}}\rangle \left( \boldsymbol{x}, t \right)\right)
  + \Delta t\boldsymbol{M} | w_i \rangle g \left( \boldsymbol{x}, t\right) 
\end{align}
where $\boldsymbol{S}$ is the $m \times m$ matrix of 
relaxation times. In the case of the single-relaxation-time 
lattice Boltzmann method, the matrix 
$\boldsymbol{S}$ will be $\boldsymbol{S} = \mathbf{I}_{m}/\tau$, 
where $\mathbf{I}_{m}$ is the $m \times m$ identity matrix, 
and $\tau$ is defined in equation \eqref{Eqn:relaxation_SRT}. 
The matrix $\boldsymbol{S}$ need not remain 
diagonal in the case of anisotropic diffusion. 
In this paper, we shall use the multiple-relaxation-time
methods proposed in \citep{Yoshida_Nagaoka_JCP_2010,
Huang_Wu_JCP_2014}. The treatment of boundary 
conditions will remain unchanged to what was presented earlier 
(i.e., equations \eqref{Eqn:Dirichlet_BC_Std}--\eqref{Eqn:Neumann_BC}).

\section{REPRESENTATIVE NUMERICAL RESULTS}
\label{Sec:LBM_S4_NR}
In this section, we employ the single- and 
multiple-relaxation-time lattice Boltzmann 
methods to solve representative diffusion 
and advection-diffusion problems. 
%
For brevity, we shall refer to the multiple-relaxation-time 
method proposed in \citep{Yoshida_Nagaoka_JCP_2010} as the 
Y-N method, and to the multiple-relaxation-time method 
proposed in \citep{Huang_Wu_JCP_2014} as the H-W method.
In all the numerical simulations that employ the 
H-W method, we have taken $c_1 = 1$, $c_2 = -2$, 
$\alpha_1 = 8$, $\alpha_2 = -8$ and $s_0 = s_1 
= \cdots = s_m = 1$ (cf., equations (15)--(17) 
in \citep{Huang_Wu_JCP_2014}).
In the rest of this section, we shall use 
the following notation:
\begin{align}
  u_{\min}(t) &= \min_{\boldsymbol{x} \in \Omega} 
  \; u(\boldsymbol{x},t) \\
  u_{\max}(t) &= \max_{\boldsymbol{x} \in \Omega} 
  \; u(\boldsymbol{x},t)
\end{align}
In all the problems that follow, the distributions 
$f_i$ are initialized as follows:
\begin{align}
  f_{i} \left( \boldsymbol{x}, t = 0\right) = 
  w_i \mathrm{u}_0 \left( \boldsymbol{x}\right)
\end{align}

\subsection{One-dimensional problems}
Consider the computational domain to be $\Omega = (0,1)$. 
The diffusion coefficient is taken as $\mathrm{D} = 1/3$, 
and the advection is neglected. The governing equations 
take the following form: 
\begin{subequations}
  \label{Eqn:1D_Problem}
  \begin{align}
    &\frac{\partial \mathrm{u}(\mathrm{x},\mathrm{t})}
    { \partial \mathrm{t}} -\mathrm{div} \left[ \mathrm{D} 
      \mathrm{grad} \left[ \mathrm{u}(\mathrm{x},\mathrm{t}) 
        \right] \right] = \mathrm{g}
    \left(\mathrm{x}, \mathrm{t}\right) \quad &&\left( \mathrm{x}, \mathrm{t}\right) 
    \in \Omega \times \left( 0 , \mathcal{T}\right] \\
    &\mathrm{u} \left( \mathrm{x}, \mathrm{t}\right) = \mathrm{u}^{\mathrm{p}} 
    \left( \mathrm{x}, \mathrm{t}\right) \quad &&\left( \mathrm{x}, \mathrm{t}\right) 
    \in \Gamma^{\mathrm{D}} \times \left[ 0 , \mathcal{T}\right] \\ 
    -&\widehat{\mathbf{n}}\cdot\mathrm{D} \mathrm{grad} \left[ \mathrm{u} \right] = 
    \mathrm{q}^{\mathrm{p}} \left( \mathrm{x}, \mathrm{t}\right)\quad 
    &&\left( \mathrm{x}, \mathrm{t}\right) \in \Gamma^{\mathrm{N}} \times 
    \left[ 0 , \mathcal{T}\right] \\ 
    &\mathrm{u} \left( \mathrm{x}, \mathrm{t} = 0\right) = \mathrm{u}_0 
    \left( \mathrm{x}\right) \quad &&\mathrm{x} \in \overline{\Omega}
  \end{align}
\end{subequations}
We employ the $D1Q3$ lattice model with 
weights $w_1 = 1/2$, and $w_2 = w_3 = 1/4$. 
Several different cases are considered below. 

\subsubsection{Uniform initial conditions.} 
\begin{figure}
  \centering
  \psfrag{time}{time}
  \psfrag{max c}{$u_{\max}\left( t \right)$}
  \psfrag{c1}{$\Delta x = 10^{-1}$}
  \psfrag{c2}{$\Delta x = 10^{-2}$}
  \psfrag{c3}{$\Delta x = 10^{-3}$}
  \includegraphics[clip,scale=0.4]{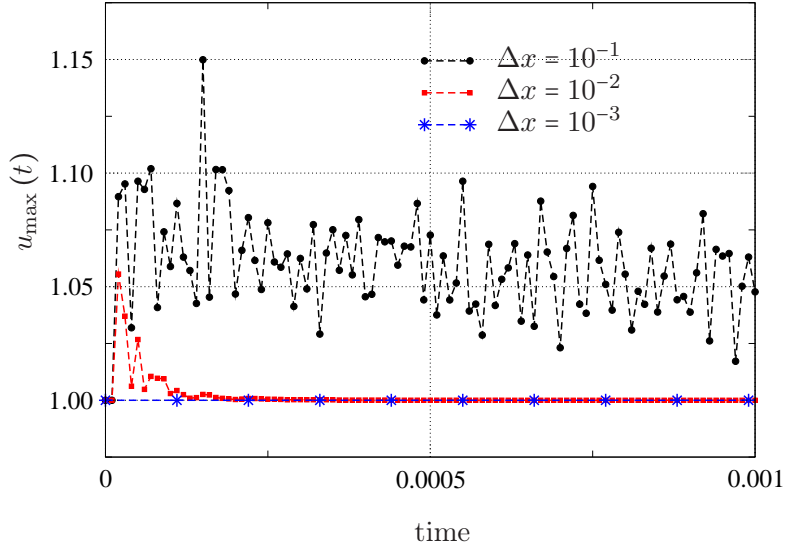}
  \caption{\textsf{One-dimensional problem with 
      uniform initial condition:}~The maximum 
    nodal concentration is plotted against time 
    for various lattice cell sizes. The time-step 
    is arbitrarily chosen as $\Delta t = 10^{-5}$ 
    for all cases. \emph{It can be observed that, 
      for a chosen time-step, reducing the lattice 
      cell size can eliminate the violation of the 
      maximum principle for one-dimensional problems.} 
    \label{Fig:1D_1_Max_c_dx}}
\end{figure}

Consider the following initial and boundary conditions:
\begin{subequations}
  \begin{align}
    &\mathrm{u}_0 \left( \mathrm{x}, \mathrm{t} = 0\right) = 
    1 \quad \forall \mathrm{x} \in \Omega \\
    &\mathrm{q}^{\mathrm{p}}(\mathrm{x} = 0,\mathrm{t}) = 0 \\
    &\mathrm{u}^{\mathrm{p}}(\mathrm{x} = 1,\mathrm{t}) = 0
  \end{align}
\end{subequations}
where $\mathrm{t} \in \left[0,10^{-2}\right]$. According 
to the maximum principle, the concentration should remain 
in $[0, 1]$ at all times. The following conclusions 
can be drawn from Figures \ref{Fig:1D_1_Max_c_dt} and 
\ref{Fig:1D_1_Max_c_dx}: 
\begin{enumerate}[(a)]
\item  The single-relaxation-time lattice Boltzmann 
  method for diffusion equations violates the maximum 
  principle. Although the minimum observed nodal 
  concentration is zero, some of the nodal 
  concentrations exceeded unity. 
\item For a given time-step, the violation 
  of the maximum principle can be reduced 
  by reducing the lattice cell size $\Delta 
  x$. As we will see later, this need not be 
  the case when the diffusion is anisotropic. 
\end{enumerate}

\subsubsection{Non-uniform initial conditions.}
\begin{figure}
  \centering
  \psfrag{time}{time}
  \psfrag{max c}{$u_{\max}\left( t \right)$}
  \psfrag{min c}{$u_{\min}\left( t \right)$}
  \psfrag{c1}{$\Delta t = 10^{-3}$}
  \psfrag{c2}{$\Delta t = 10^{-4}$}
  \psfrag{c3}{$\Delta t = 10^{-5}$}
  \includegraphics[scale=0.4,clip]{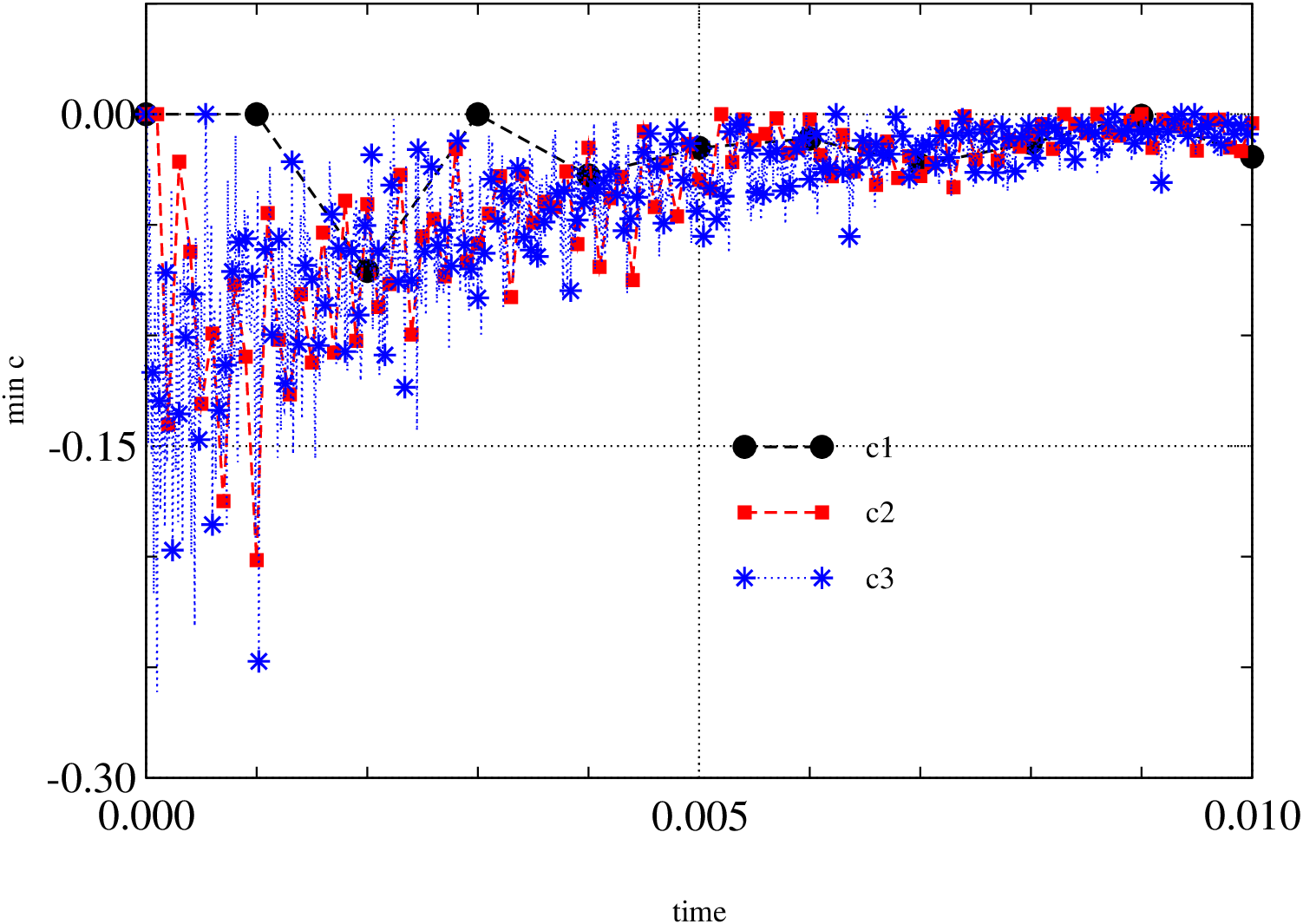}
  \caption{\textsf{One-dimensional problem with non-uniform 
      initial condition:}~The minimum nodal concentration is 
    plotted against time for various time-steps. The lattice 
    cell size is taken as $\Delta x = 0.1$. For the given 
    input data, the maximum principle asserts that the 
    concentration should be between $0$ and $1$ for all 
    times. \emph{The numerical solutions under the 
    LBM have violated the non-negative constraint. 
    It can be observed that reducing the time-step 
    may not eliminate the violation of the non-negative 
    constraint.} \label{Fig:1D_2_Min_Max_dt}}
\end{figure}

\begin{figure}
  \centering
  \psfrag{time}{time}
  \psfrag{max c}[l]{$u_{\max}\left( t \right)$}
  \psfrag{min c}[l]{$u_{\min}\left( t \right)$}
  \psfrag{c1}{$\Delta x = 10^{-1}$}
  \psfrag{c2}{$\Delta x = 10^{-2}$}
  \psfrag{c3}{$\Delta x = 10^{-3}$}
  \subfigure[]{
    \includegraphics[scale = 0.4, clip]{Figures/1D_2/max_c_dx.eps}}
  \subfigure[]{
    \includegraphics[scale = 0.4, clip]{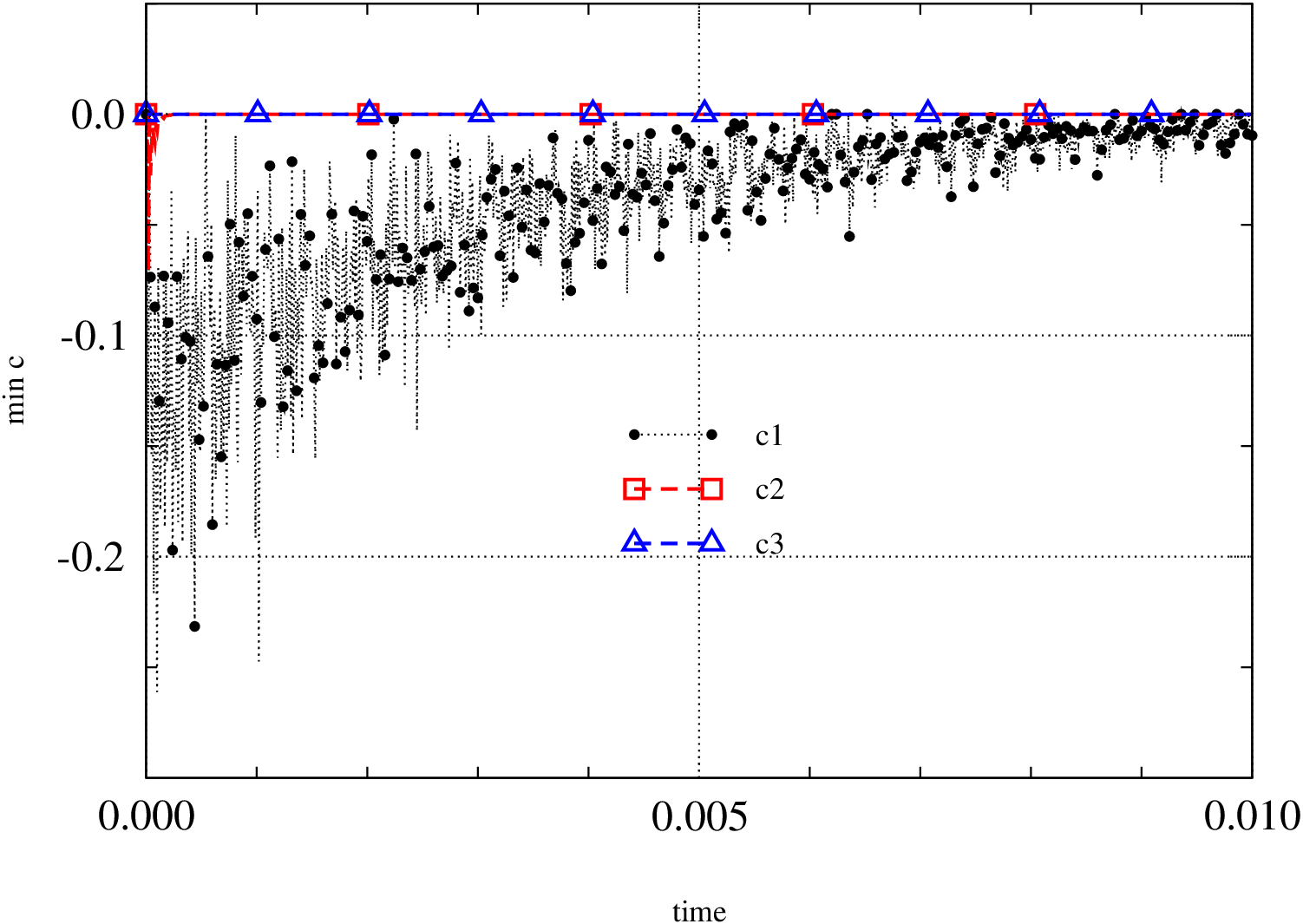}}
  \caption{\textsf{One-dimensional problem with non-uniform 
      initial condition:}~Maximum and minimum nodal 
    concentrations are plotted against the time for various 
    lattice cell sizes. The lattice time-step is taken 
    as $\Delta t = 10^{-5}$. The concentration should be 
    between $0$ and $1$. The numerical solutions have 
    clearly violated the maximum (see near $t = 0$) and 
    minimum bounds. \emph{However, it is observed that, 
      for a fixed time-step, reducing the lattice cell 
      size $\Delta x$ can decrease the violation.} 
    \label{Fig:1D_2_Min_Max_dx}}
\end{figure}

~Consider 
the following initial condition:
\begin{align}
  \mathrm{u}_0 \left( \mathrm{x}, \mathrm{t} = 0\right) = \left\{
  \begin{array}{l l}
    1  & 	\mathrm{x} \in \left[  0.4,  0.6\right] \\
    0  & 	\mathrm{otherwise}
  \end{array} \right.
\end{align}
The boundary conditions are:
\begin{subequations}
  \begin{align}
    &\mathrm{u}^{\mathrm{p}}\left( \mathrm{x} = 1,\mathrm{t} \right) = 0 \\
    &\mathrm{q}^{\mathrm{p}}\left( \mathrm{x} = 0,\mathrm{t} \right) = 0 
  \end{align}
\end{subequations}
Even for this problem, the concentration should remain 
in $[0, 1]$. Maximum and minimum nodal concentrations 
($u_{\min}(t)$ and $u_{\max}(t)$) are shown in Figures 
\ref{Fig:1D_2_Min_Max_dt} and \ref{Fig:1D_2_Min_Max_dx}. 
The following conclusions can be made from these figures: 
\begin{enumerate}[(a)]
\item The maximum and minimum bounds can be violated 
simultaneously. 
\item Decreasing the time-step while keeping the 
lattice cell size constant will not alleviate the 
violations of the non-negative constraint and the 
maximum principle. 
\item Just like the previous problem, the violations 
vanish with the refinement of the lattice cell size 
for a fixed time-step. 
\end{enumerate}

\subsubsection{Uniform initial condition with constant source.} 
We have taken $\mathrm{u}_{0}(\mathrm{x}) = 0$ on the entire 
domain $\Omega = (0,1)$, and $\Gamma^{\mathrm{D}} = \Gamma$ 
with $\mathrm{u}^{\mathrm{p}}(\mathrm{x} \in \Gamma,\mathrm{t}) 
= 0$. The time interval of interest is taken as $\mathcal{T} = 10^{-2}$. 
The source is taken as $\mathrm{g}(\mathrm{x},\mathrm{t}) = 1$ 
for $\mathrm{x} \in \Omega$ and $\mathrm{t} \in (0,\mathcal{T}]$. 
The error will be calculated as follows:
\begin{align}
  \mathrm{Error}(t) = \frac{1}{N} \sqrt{\sum_{i = 1}^{N} 
    \left(u(x_i,t) - \mathrm{u}_{\mathrm{exact}}(x_i,t)\right)^2} 
\end{align}
where $N$ is the total number of lattice nodes, 
and $x_i$ is the spatial coordinate of the $i$-th 
lattice node. The exact solution to this problem 
is denoted by $\mathrm{u}_{\mathrm{exact}}(x,t)$.
\emph{
The obtained numerical results are summarized in Table 
\ref{Tbl:1D_3}, and one can conclude that the weighted 
splitting method for discretizing the Dirichlet boundary 
conditions produces non-negative nodal concentrations 
for one-dimensional problems. However, the method comes 
at an expense of marginal decrease in accuracy.}

\begin{table}
  \caption{\textsf{One-dimensional problem with 
      uniform initial condition and constant source:}~ 
    This table compares the numerical solutions under 
    the standard and weighted splitting methods for 
    discretizing Dirichlet boundary conditions for 
    several choices of $\Delta x$ and $\Delta t$. In 
    this numerical experiment, one observes that the 
    weighted splitting method produces non-negative 
    solutions, and Its accuracy is comparable to the 
    standard method of discretization. \label{Tbl:1D_3}}
  \begin{tabular}{| c | c | c | c | c | c | c |} \hline
    \multirow{2}{*}{$\Delta t$} & \multirow{2}{*}{$\Delta x$} & \multirow{2}{*}
    {$1 / \tau$} & \multicolumn{2}{|c|}{standard method} & \multicolumn{2}{|c|}
    {weighted splitting method} \\ 
    \cline{4-5} \cline{6-7} 
    & & & violation & Error$\left( \mathcal{T}\right)$ & violation & Error$\left( \mathcal{T}\right)$  \\ \hline
    $10^{-3}$  &  $10^{-3}$ & $9.995 \times 10^{-4}$ & Yes & $3.14\times10^{-4}$ & No &  $3.15\times10^{-4}$ \\ \hline
    $10^{-4}$  &  $10^{-3}$ & $1.000 \times 10^{-1}$ & Yes & $2.98\times10^{-4}$ & No & $3.10\times10^{-4}$ \\ \hline
    $10^{-5}$  &  $10^{-3}$ & $9.520 \times 10^{-2}$ & Yes & $2.90\times10^{-4}$ & No &  $2.94\times10^{-4}$  \\ \hline
    $10^{-6}$  &  $10^{-3}$ & $6.667 \times 10^{-1}$ & Yes & $2.90\times10^{-4}$ & No &  $2.90\times10^{-4}$ \\ \hline
    $10^{-7}$  &  $10^{-3}$ & $1.667 \times 10^{0}$ & Yes & $2.90\times10^{-4}$ & No &  
    $2.90\times10^{-4}$ \\ \hline 
  \end{tabular}
\end{table}

\subsubsection{On comparison principle.} Consider 
the following boundary and initial conditions:
\begin{subequations}
  \begin{alignat}{2}
    &\mathrm{u}_0 (\mathrm{x},\mathrm{t} = 0) = 0 
    && \quad \mathrm{x} \in \Omega \\
    &\mathrm{u}^{\mathrm{p}}(\mathrm{x} = 0,\mathrm{t}) 
    = \mathrm{u}_{\mathrm{L}} 
    && \quad \mathrm{t} \in [0,\mathcal{T}] \\
    &\mathrm{u}^{\mathrm{p}} (\mathrm{x} = 1,\mathrm{t}) = 0 
    && \quad \mathrm{t} \in \left[ 0, \mathcal{T}\right]
  \end{alignat}
\end{subequations}
In this numerical experiment, we have taken 
$\mathcal{T} = 0.01$, and the source term is 
taken to be zero (i.e., $\mathrm{g}(\mathrm{x},
\mathrm{t}) = 0$). We shall solve the problem 
using several cases of $\Delta x$, $\Delta t$, 
and $\mathrm{u}_{\mathrm{L}}$. Both the standard 
and weighted splitting methods are employed in 
separate test runs.

\begin{figure}
  \centering
  \psfrag{x}{$x$-axis}
  \psfrag{con}{$u(x,\mathcal{T})$}
  \psfrag{c1}{$\mathrm{u}^{\mathrm{p}}(\mathrm{x}=0,\mathrm{t})=1$}
  \psfrag{c2}{$\mathrm{u}^{\mathrm{p}}(\mathrm{x}=0,\mathrm{t})=2$}
  \psfrag{c3}{$\mathrm{u}^{\mathrm{p}}(\mathrm{x}=0,\mathrm{t})=3$}
  \subfigure[Violation of the comparison principle near $x=0.3$.]{
    \includegraphics[clip,scale=0.4]{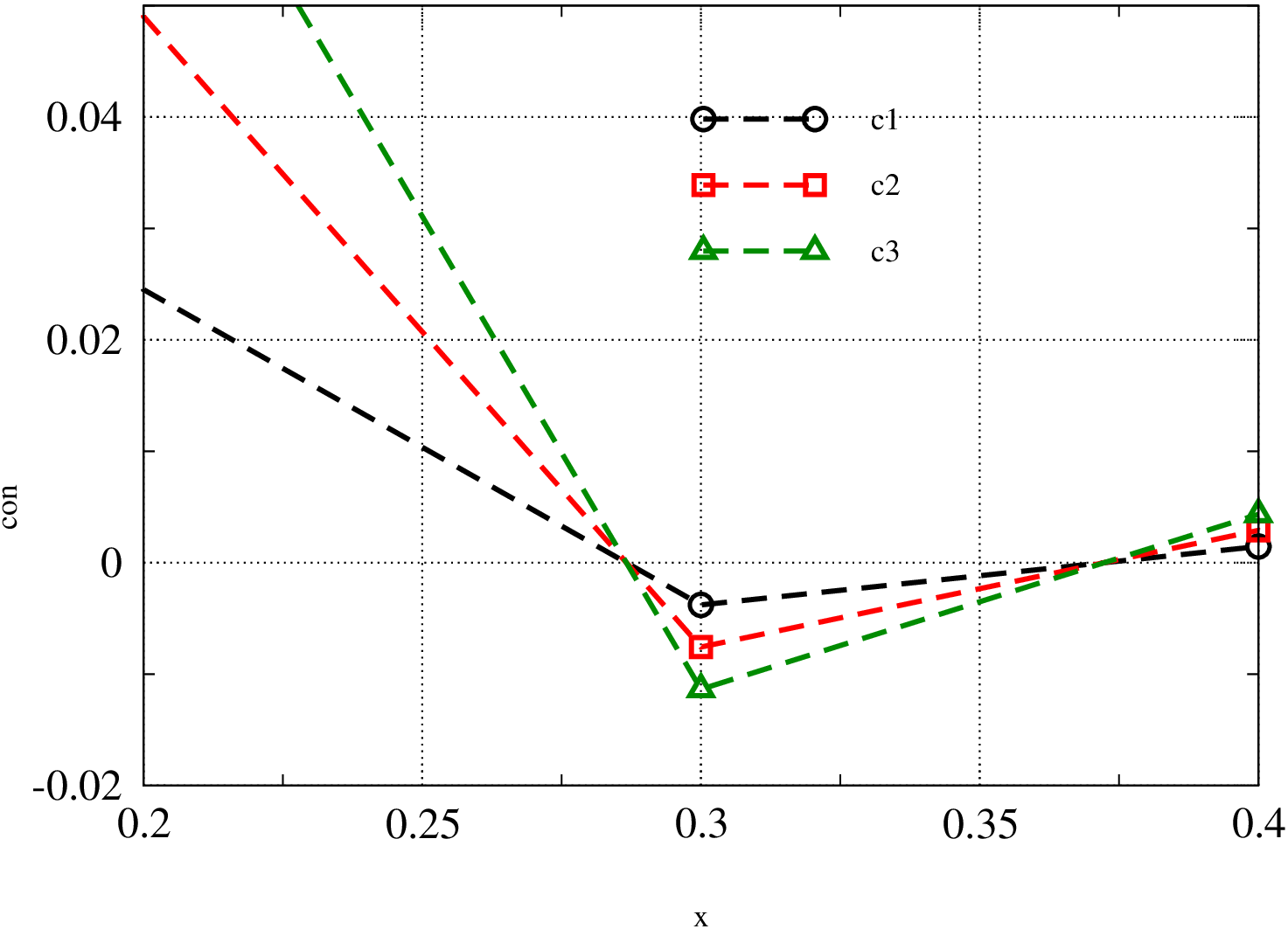}}
  \subfigure[Violation of the comparison principle near $x = 0.5$ and 
    $x = 0.7$.]{\includegraphics[clip,scale=0.4]{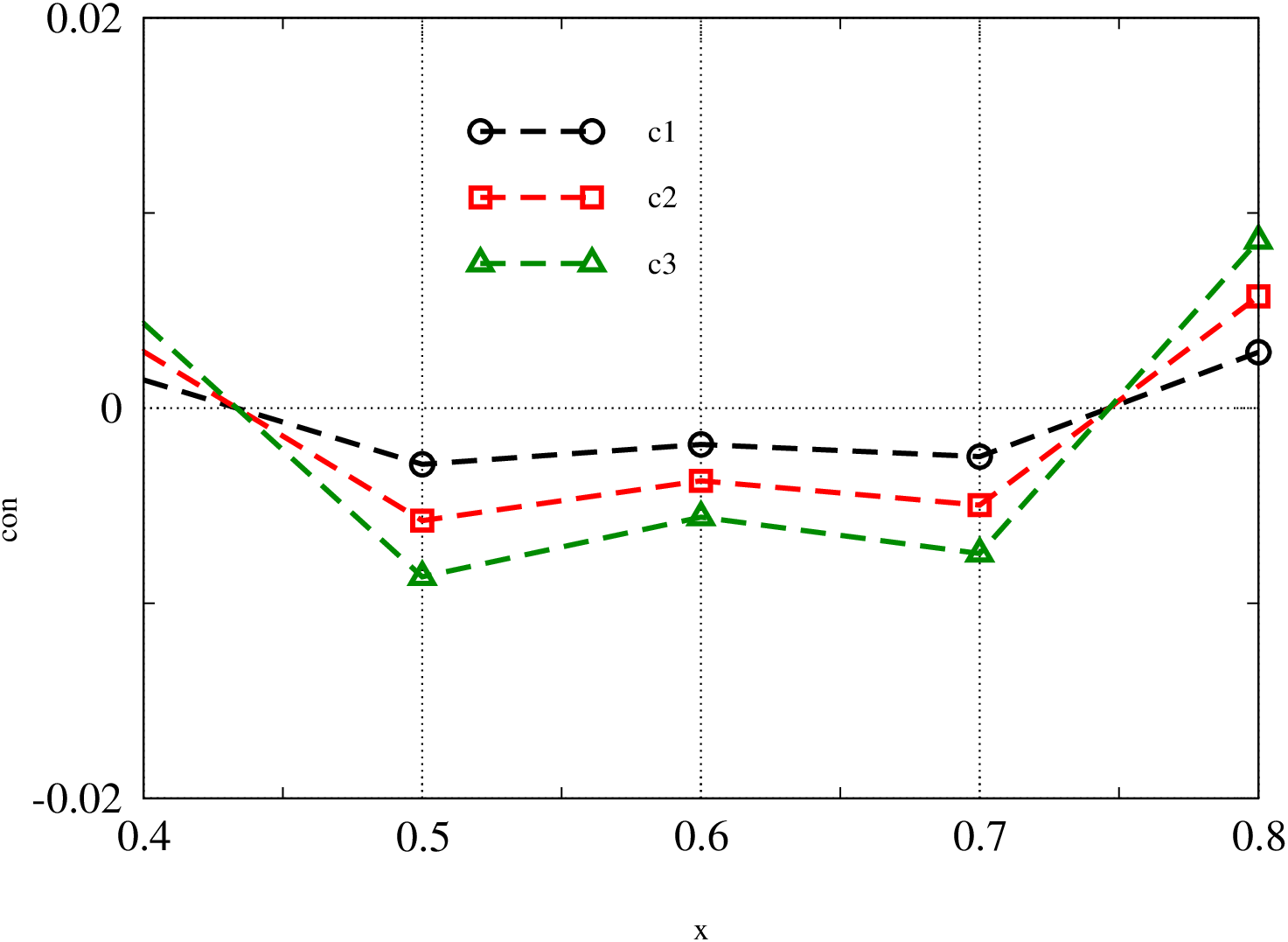}}
  \caption{\textsf{One-dimensional problem and the 
      comparison principle:}~The source is taken to 
    be zero (i.e., $\mathrm{g}(\mathrm{x},\mathrm{t}) 
    = 0$), and the initial concentration is zero in 
    the entire domain. The lattice cell size is $\Delta 
    x = 0.1$ and the time-step is $\Delta t = 10^{-5}$. 
    The prescribed concentration at $\mathrm{x} = 1$ 
    is zero, and three different cases are considered 
    for the prescribed concentration at $\mathrm{x}=1$, 
    as indicated in the figure. 
    \emph{It can be observed that the LBM based on 
      the $D1Q3$ lattice model violates the comparison 
      principle.} \label{Fig:1D_Comp}}
\end{figure}

For this problem, if $\mathrm{u}_{\mathrm{L}}^{1} \leq 
\mathrm{u}_{\mathrm{L}}^{2}$ then the comparison principle 
implies that $\mathrm{u}^{1}(\mathrm{x},\mathrm{t}) 
\leq \mathrm{u}^{2}(\mathrm{x},\mathrm{t})$ for 
$\forall (\mathrm{x}, \mathrm{t}) \in \overline{\Omega} \times [0,\mathcal{T}]$. 
A sample result is presented in Figure \ref{Fig:1D_Comp}.
Several other numerical experiments are performed using 
different choices of $\Delta x$ and $\Delta t$. Figure 
\ref{Fig:1D_Comp} and these numerical experiments reveal 
the following conclusions:
\begin{enumerate}[(a)]
\item LBM-based formulations do not, in general, 
  respect the comparison principle. 
\item For one-dimensional problems, refining the lattice 
  cell size $\Delta x$ for a given time-step $\Delta t$ 
  can remove the violations of the comparison principle.
\item The weighted splitting method does not guarantee 
  the satisfaction of the comparison principle even in 
  1D. 
\end{enumerate}

\subsection{Two-dimensional problem with anisotropic 
  diffusion on a non-convex domain}
\begin{figure}
  \centering
  \psfrag{ly}{$L$}
  \psfrag{lx}{$L$}
  \psfrag{nf}{zero flux}
  \psfrag{d}{$L/10$}
  \psfrag{gin}{$\Gamma_{\mathrm{interior}}$}
  \psfrag{gout}{$\Gamma_{\mathrm{exterior}}$}
  \includegraphics[scale=0.9,clip]{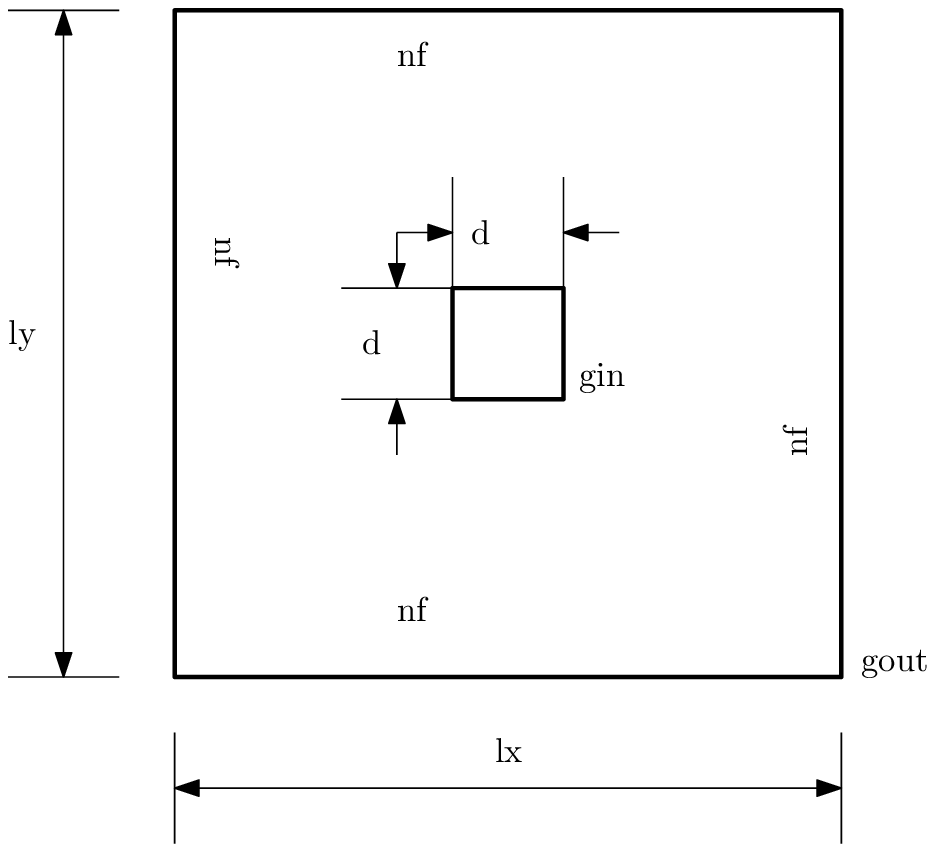}
  \caption{\textsf{Two-dimensional problem with 
      anisotropic diffusion in a non-convex 
      domain:}~This figure provides a pictorial 
    description of the test problem. A concentration 
    of $\mathrm{u}^{\mathrm{p}} = 1$ is prescribed 
    on the inner boundary. \label{Fig:2DYN_des}}
\end{figure}

\begin{figure}
  \centering
  \includegraphics[scale=0.4,clip]{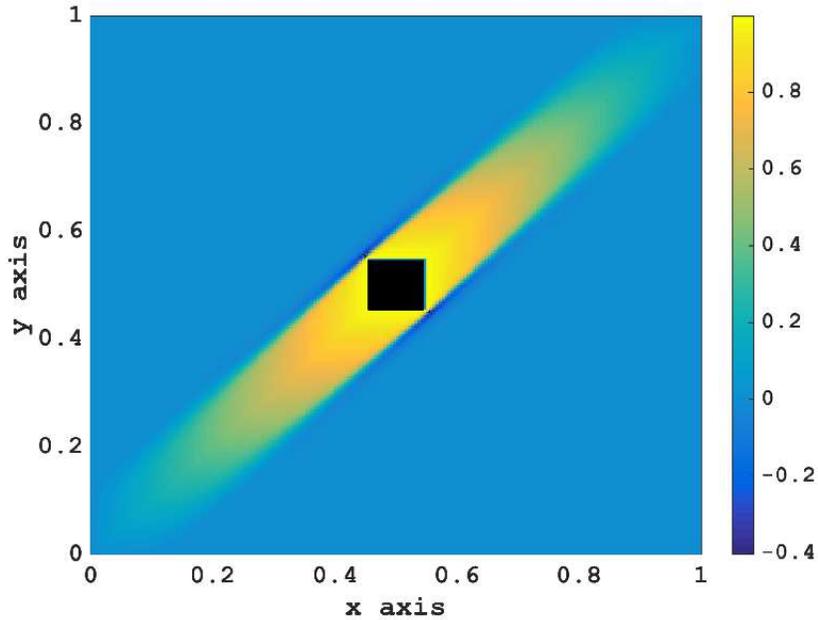}
  \caption{\textsf{Two-dimensional problem with 
      anisotropic diffusion in a non-convex 
      domain:}~This figure shows the concentration 
    profile at time $t = 0.01$ under the Y-N 
    method for the data given by Case 3 in Table 
    \ref{Tbl:2DYN_params}. The diffusivity tensor 
    is anisotropic but spatially homogeneous. 
    \emph{Significant negative values for the 
      concentration are observed in the 
      computational domain.} \label{Fig:2DYN_c}}
\end{figure}

\begin{figure}
  \centering
  \subfigure[Case 1]{
    \includegraphics[scale=0.3,clip]{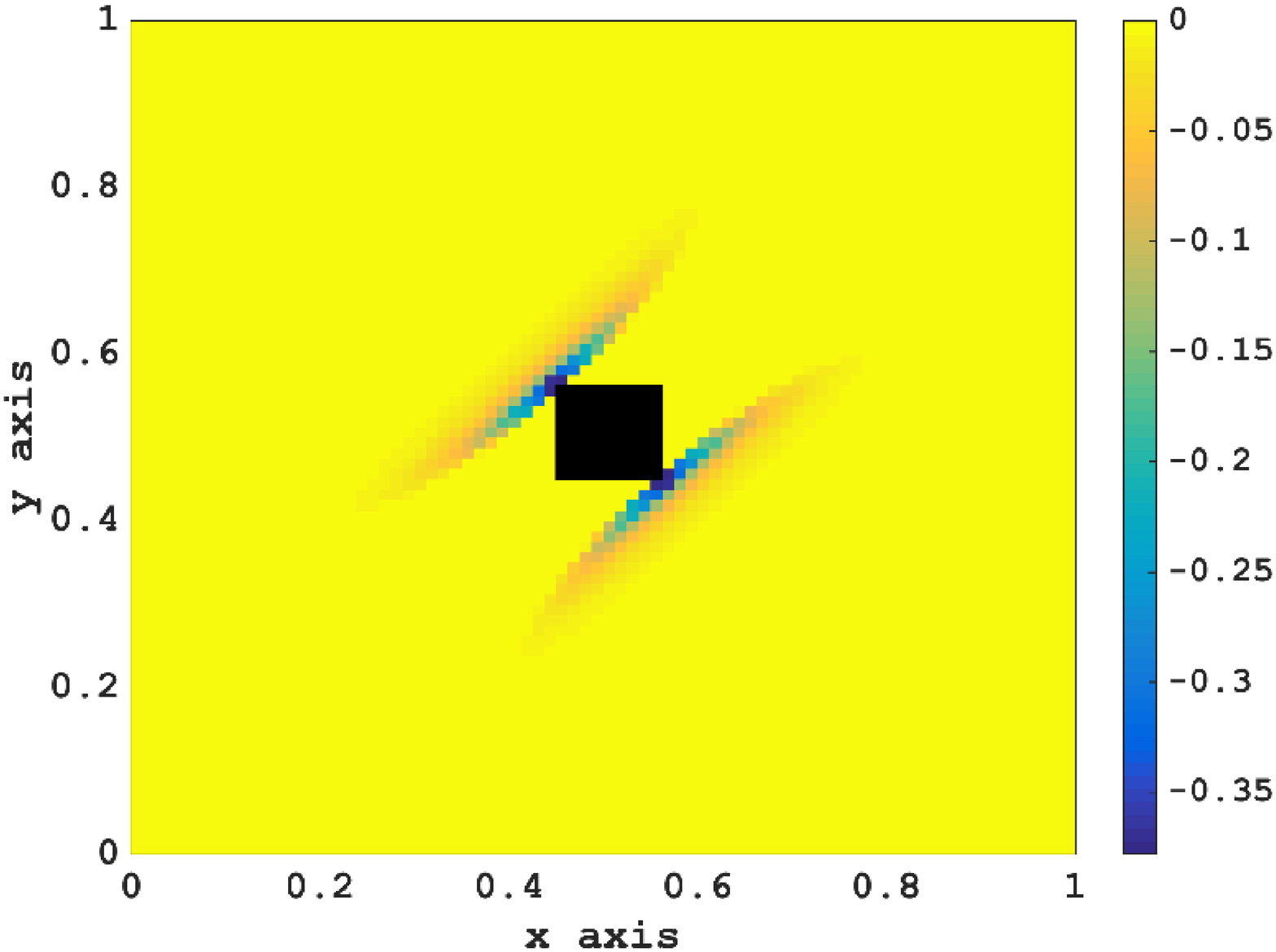}}
  \subfigure[Case 2]{
    \includegraphics[scale=0.3,clip]{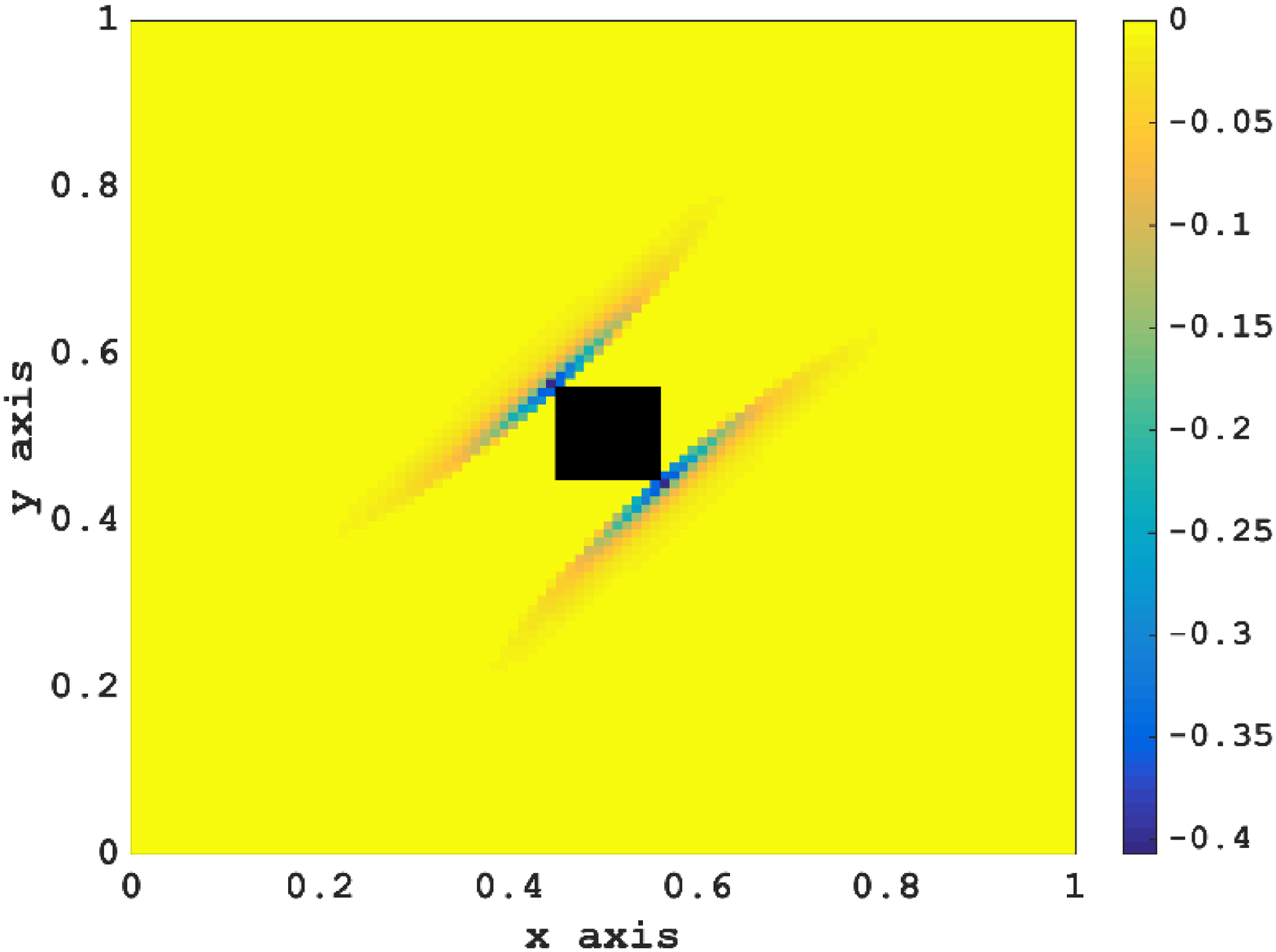}}
  \subfigure[Case 3]{
    \includegraphics[scale=0.3,clip]{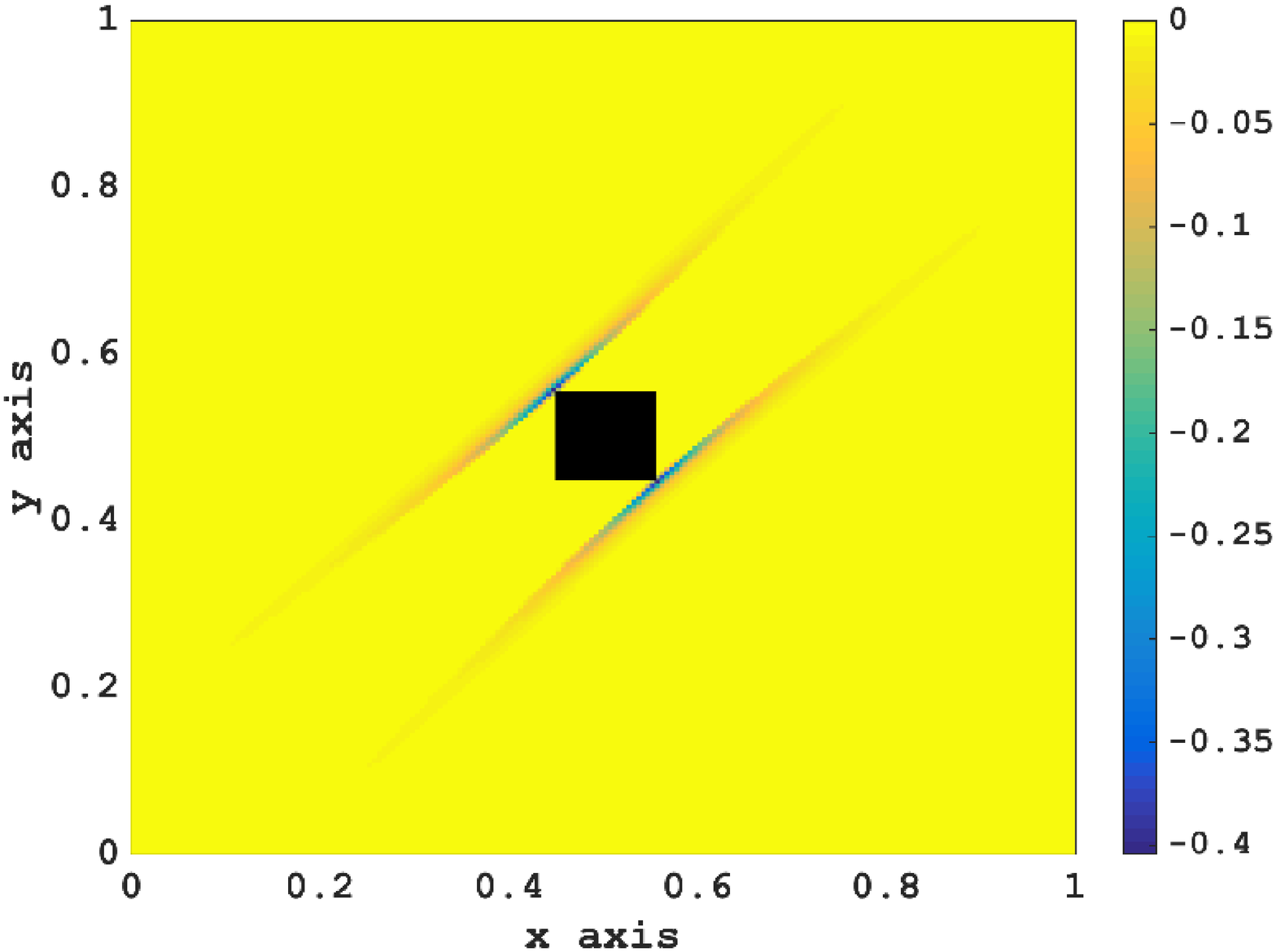}}
  \caption{\textsf{Two-dimensional problem with 
      anisotropic diffusion in a non-convex 
      domain:}~The figure shows the regions 
    where the non-negative constraint is 
    violated under the Y-N method. 
    \label{Fig:2DYN_neg}}
\end{figure}

\begin{figure}
  \centering
  \psfrag{c1}{Case 1}
  \psfrag{c2}{Case 2}
  \psfrag{c3}{Case 3}
  \psfrag{time}{time}
  \psfrag{minc}{$u_{\min}\left( t \right)$}
  \includegraphics[scale=0.4, clip]{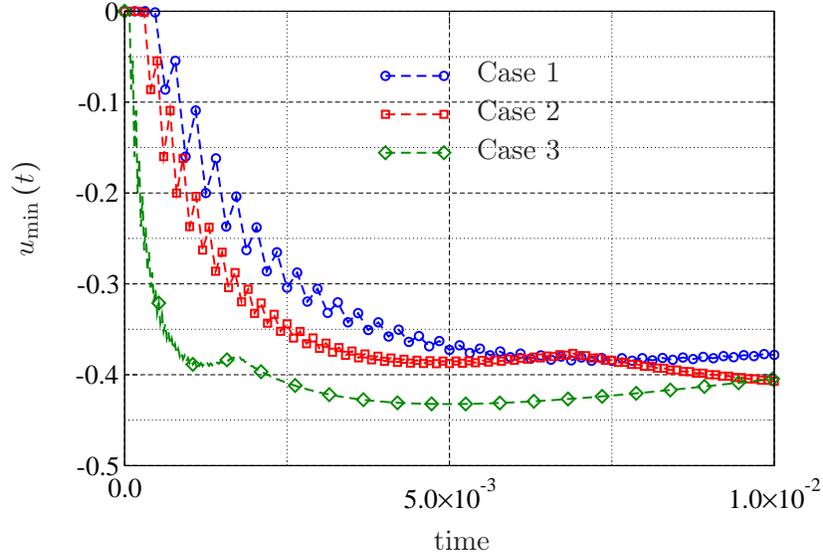}
  \caption{\textsf{Two-dimensional problem with anisotropic 
      diffusion in a non-convex domain:}~The figure shows 
    the variation of minimum value of the concentration 
    with respect to time under the Y-N multiple-relaxation-time 
    lattice Boltzmann method. \emph{It is evident from 
      this numerical experiment that refining the discretization 
      parameters, $\Delta t$ and $\Delta x$, does not alleviate 
      the violation of the non-negative constraint in the 
      case of anisotropic diffusion.} \label{Fig:2DYN_minc}}
\end{figure}

We now examine the Y-N multiple-relaxation-time method 
for anisotropic diffusion tensor. The computational 
domain is shown in Figure \ref{Fig:2DYN_des}. We 
have taken $L = 1$ and $\Gamma^{\mathrm{D}} = 
\Gamma_{\mathrm{outer}} \cup \Gamma_{\mathrm{inner}}$.
On the inner boundary the prescribed, concentration is prescribed to be unity 
(i.e., $\mathrm{u}^{\mathrm{p}}\left( \mathbf{x}, \mathrm{t}\right) = 1\;
\text{for} \; \mathbf{x}\in \Gamma_{\mathrm{inner}}$). The flux is
prescribed to be zero on the outer boundary (i.e., 
$\mathrm{q}^{\mathrm{p}}\left( \mathbf{x}, \mathrm{t}\right) = 0\;
\text{for}\;\mathbf{x}\in \Gamma_{\mathrm{outer}}$).
The anisotropic diffusion tensor is taken as:
\begin{align}
  \mathbf{D}(\mathbf{x}) = \mathbf{R}_{\theta}^{\mathrm{T}} 
  \mathbf{D}_0 \mathbf{R}_{\theta}
\end{align}
where 
\begin{align}
  \mathbf{D}_0 = \left[ \begin{array}{c c}
      10 	& 	0 \\
      0 	& 	10^{-3}
    \end{array} \right]
\end{align}
The orthogonal rotation matrix is denoted by 
$\mathbf{R}_{\theta}$, where $\theta$ denotes 
the angle of rotation. Herein, we have taken 
$\theta = \pi/4$. We employed the $D2Q5$ lattice 
model. The discrete velocity directions are given 
by
\begin{align}
  \boldsymbol{e}_i^{\mathrm{T}} = \left\{ \begin{array}{ll}
    \left[ 0 , 0\right] \quad &i = 0 \\
    c\left[ \mathrm{cos}\left( (i - 1)\pi\right), 
      \mathrm{sin}\left( (i - 1)\pi\right)\right] \quad &i=1,2 \\
    c\left[ \mathrm{cos}\left( (2i - 5)\pi/2\right), 
      \mathrm{sin}\left( (2i - 5)\pi/2\right)\right] \quad &i=3,4
  \end{array} \right.
\end{align}
where $c = \Delta x/\Delta t$. The respective 
weights are taken as:
\begin{align}
  w_i = \left\{ \begin{array}{ll}
    1/3 	\quad &i = 0 \\
    1/6 	\quad &i = 1,2,3,4
  \end{array} \right.
\end{align}
The time interval of interest is taken as $\mathcal{T} = 10^{-2}$. 
We employed the standard method of enforcing Dirichlet boundary 
conditions (see equation \eqref{Eqn:Dirichlet_BC_Std}).
Table \ref{Tbl:2DYN_params} provides the discretization 
parameters employed in this paper. 
Figures \ref{Fig:2DYN_c}--\ref{Fig:2DYN_minc} show that 
the Y-N method violates the 
non-negative constraint. In fact, the obtained minimum 
concentration is about $-0.4$, which is a significant 
violation given the fact that the concentration should 
be between $0$ and $1$. Another noticeable feature in 
all the cases considered, the minimum concentration 
converged to a negative value as the 
time progressed.

\begin{table}
\caption{\textsf{Two-dimensional problem with 
    anisotropic diffusion tensor on a non-convex 
    domain:}~This table provides the minimum 
  concentrations for various discretization 
  parameters (i.e., $\Delta x$ and $\Delta t$). 
  We have taken $\Delta x^2 = \Delta t$.}
\label{Tbl:2DYN_params}
\begin{tabular}{| c | c | c | c |}
\hline
Case 	&	$\Delta x$	& 	$\Delta t$  & $u_{\min}\left( \mathcal{T}\right)$\\ \hline
1		& 	$1.25 \times 10^{-2}$	&	$1.5625 \times 10^{-4}$ &  -0.3781 \\ \hline
2 		&	$1.00 \times 10^{-2}$	&	$1.0000 \times 10^{-4}$ &  -0.4072 \\ \hline
3		&	$5.00 \times 10^{-3}$	&	$2.5000 \times 10^{-5}$ &  -0.4044 \\ \hline
\end{tabular}
\end{table}

\subsection{Two-dimensional problem with anisotropic and 
heterogeneous diffusion tensor}
\begin{figure}
  \centering
  \includegraphics[clip,scale=0.4]{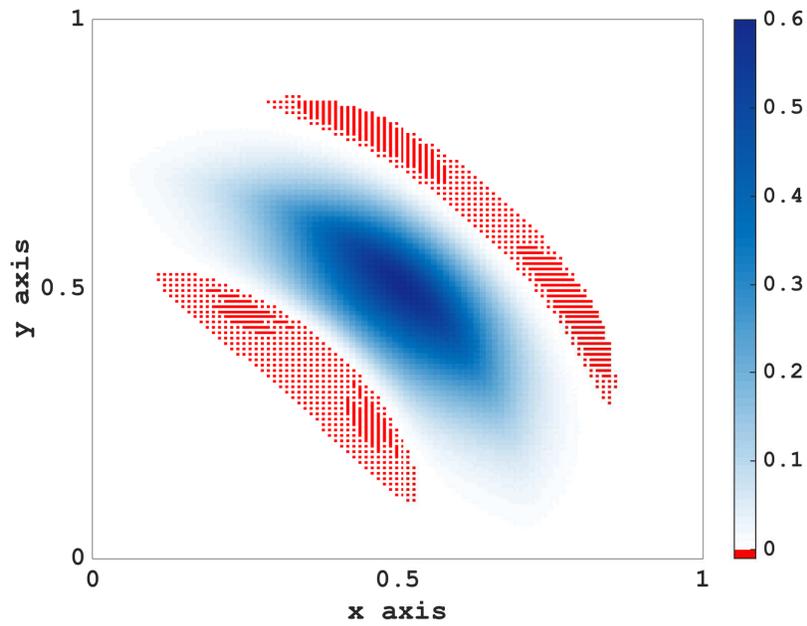}
  \caption{\textsf{Two-dimensional problem with 
      anisotropic and heterogeneous diffusion 
      tensor:}~This figure shows the concentration 
    profile at $t=0.025$ under the H-W method. The 
    parameters for this numerical experiment are 
    given by Case 5 of Table \ref{Tbl:H_W_2D_Anisotropic}. 
    The minimum value of $u$ is negative with $u_{\min}
    (\mathcal{T}) = -6.8\times10^{-3}$. (See the red 
    regions on the online color version of this paper.)
    \label{Fig:H_W_2D_Anisotropic}}
\end{figure}

\begin{figure}
  \centering
  \subfigure[Case 1]{
    \includegraphics[clip,scale=0.2]{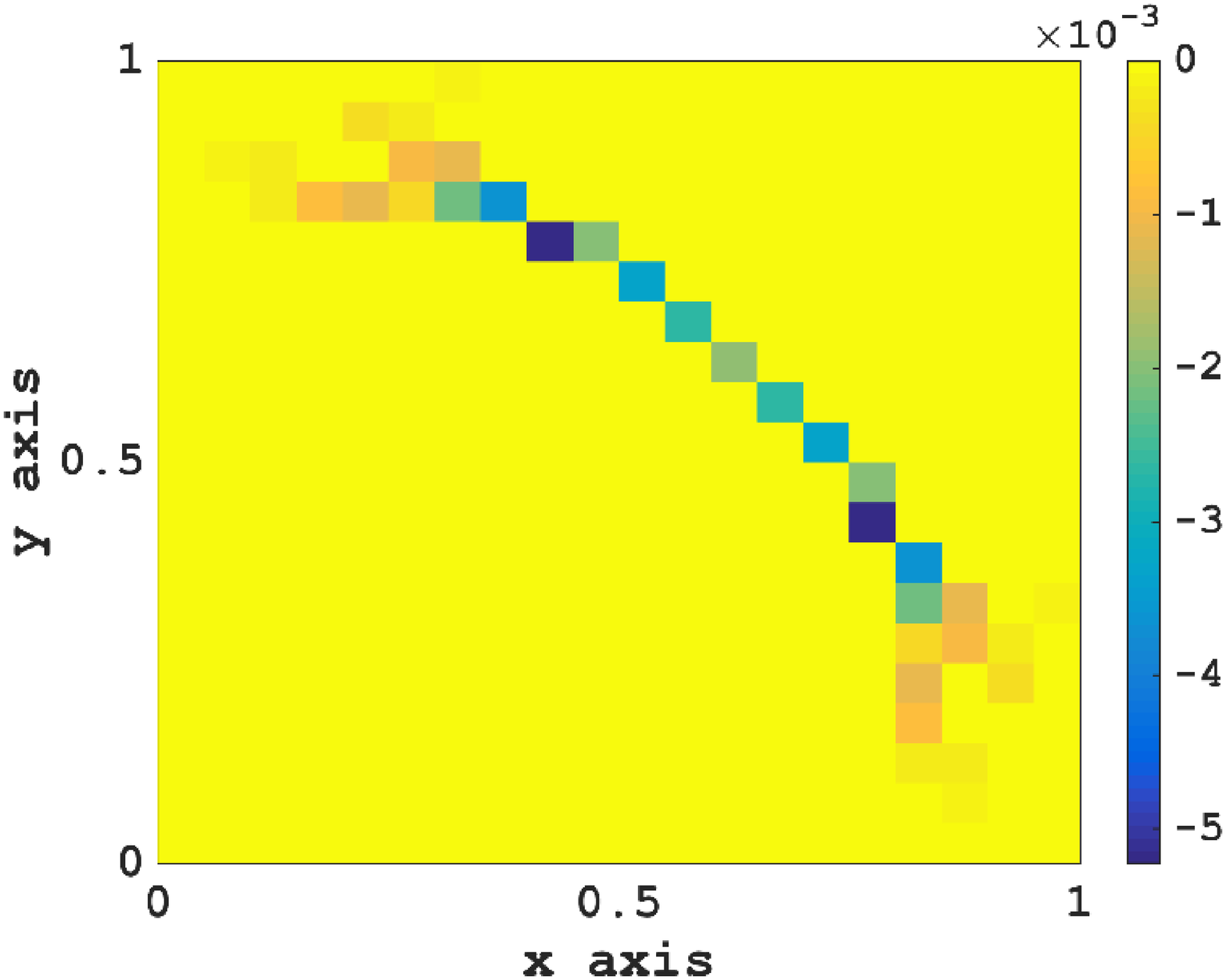}}
  \subfigure[Case 2]{
    \includegraphics[clip,scale=0.2]{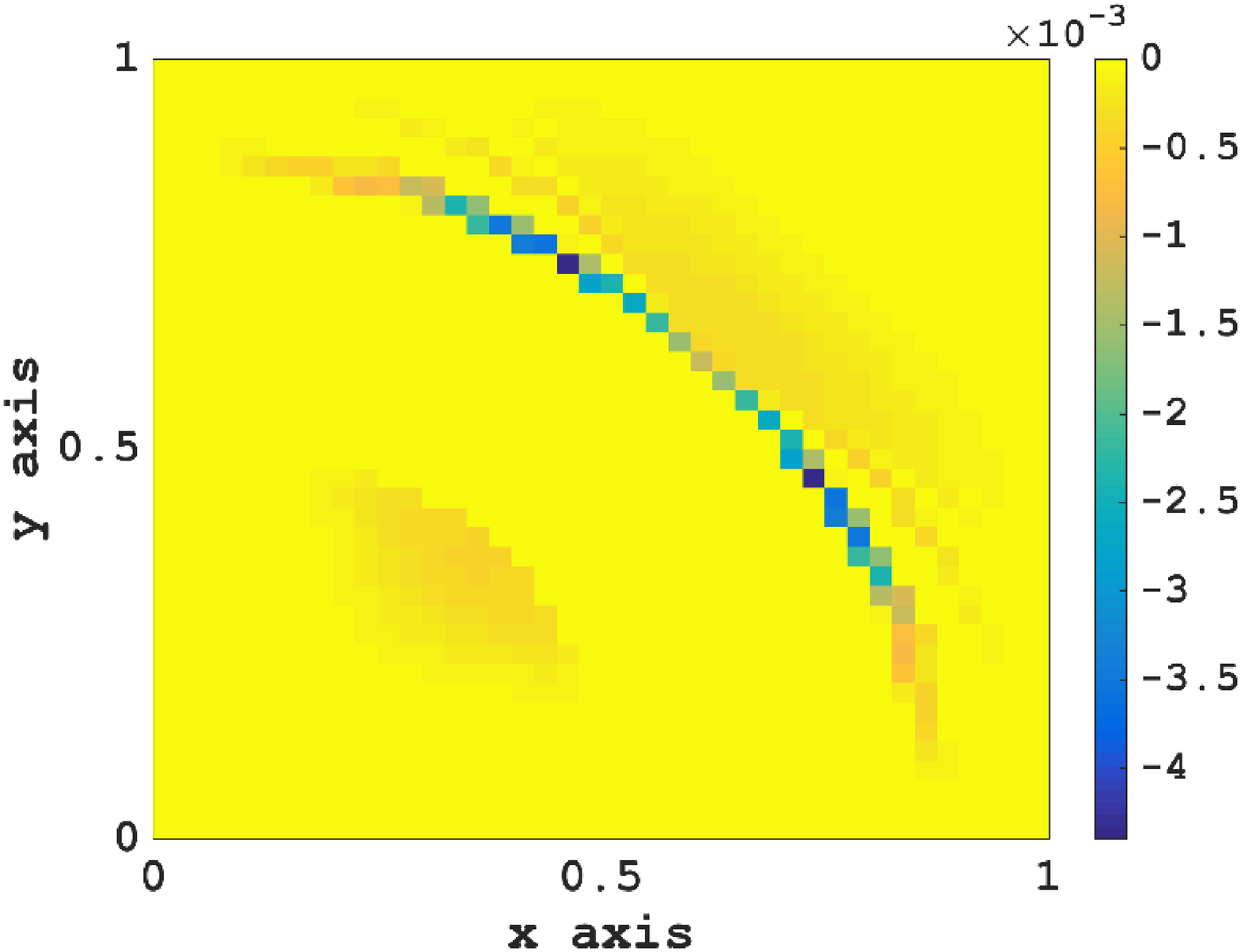}}
  \subfigure[Case 3]{
    \includegraphics[clip,scale=0.2]{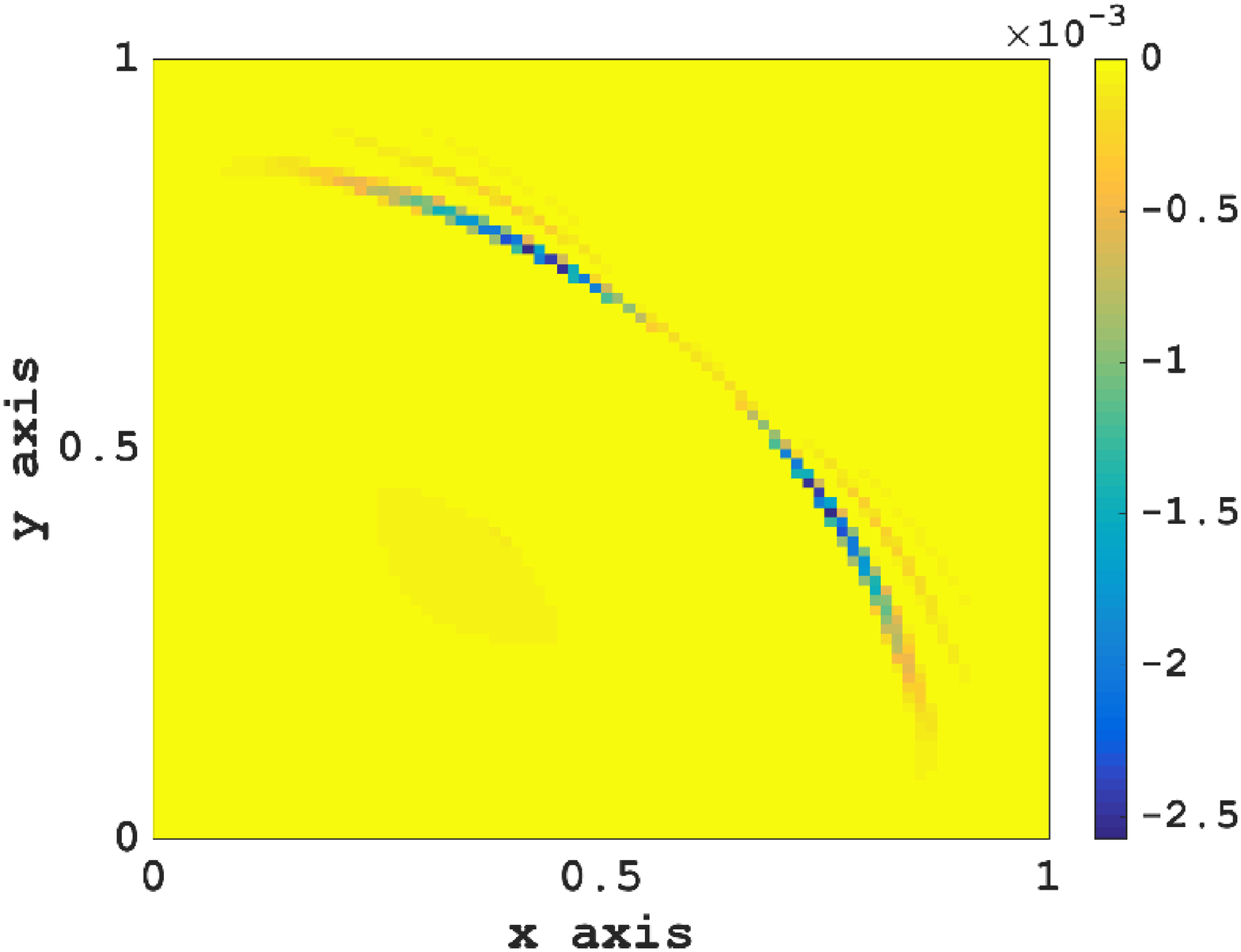}}
  \subfigure[Case 4]{
    \includegraphics[clip,scale=0.2]{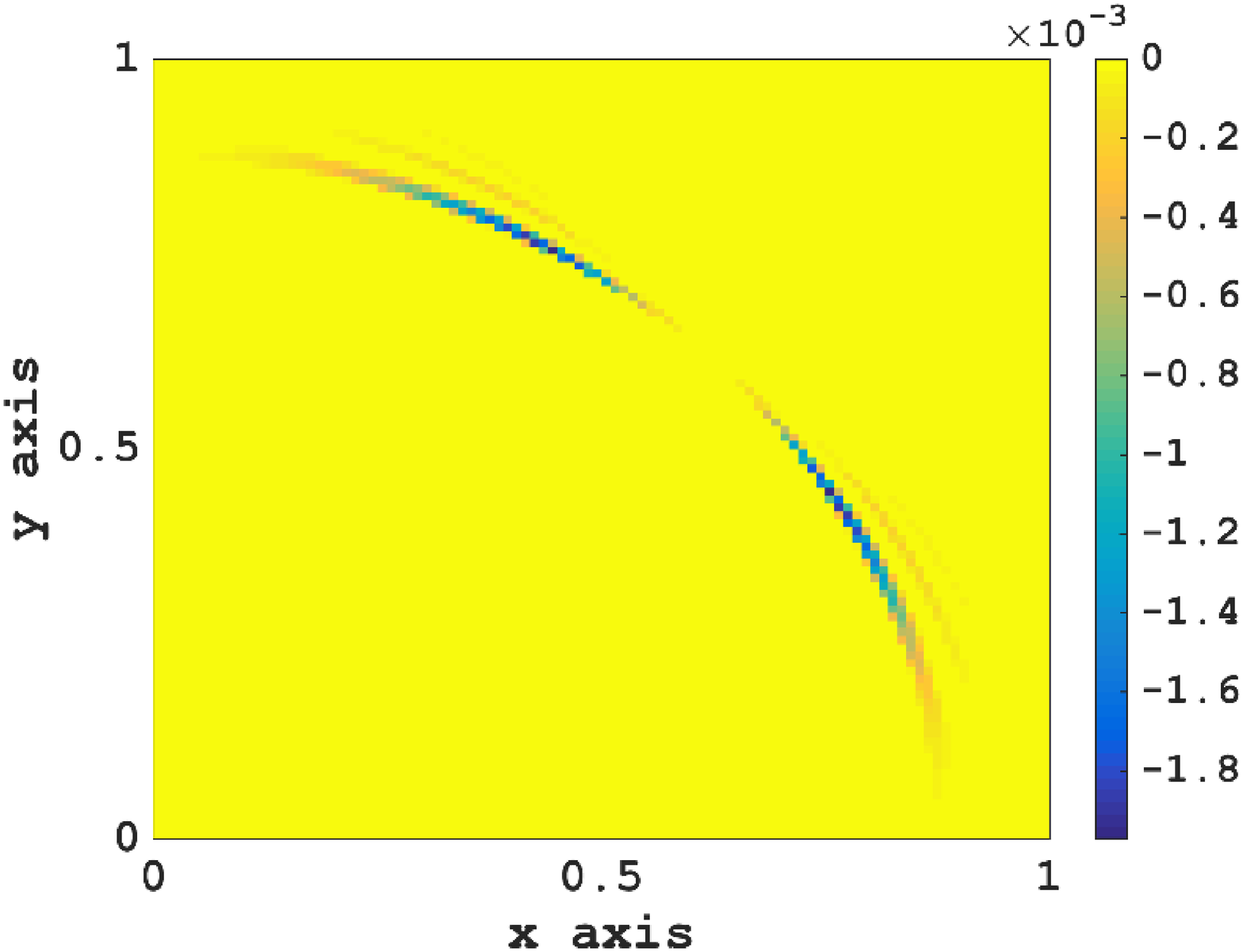}}
  \subfigure[Case 5]{
    \includegraphics[clip,scale=0.2]{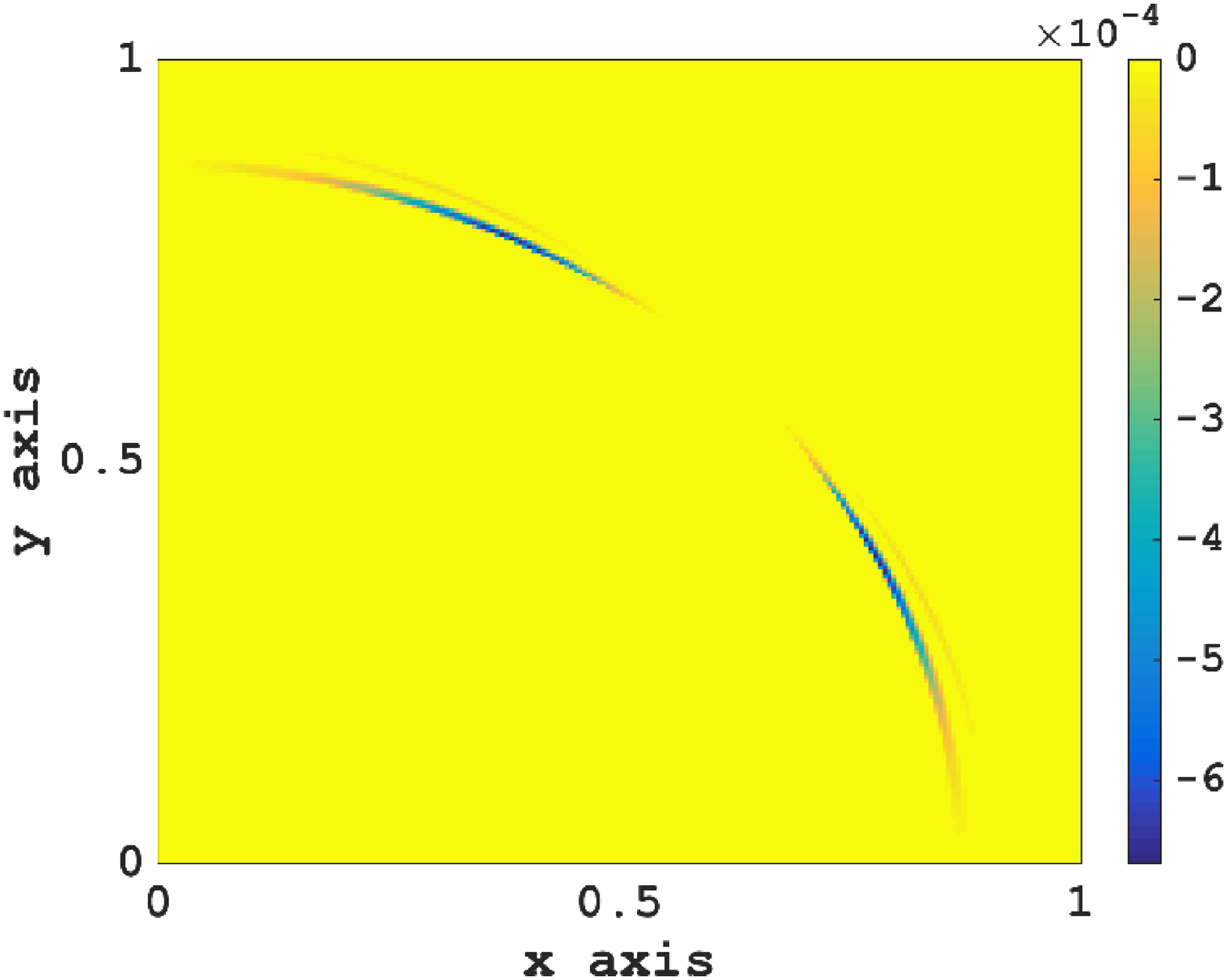}}
  \caption{\textsf{Two-dimensional problem with 
      anisotropic and heterogeneous diffusion 
      tensor:}~This figure shows the regions 
    that have negative values for the concentration 
    under the H-W method. \label{Fig:H_W_2D_Neg}}
\end{figure}

\begin{figure}
  \centering
  \psfrag{time}{time}
  \psfrag{minPhi}{$u_{\min}(t)$}
  \psfrag{intPhi}{$\mathcal{J}_{1}\left(u;\Omega;t\right)$}
  \psfrag{intPhipos}{$\mathcal{J}^{+}_{1}\left(u;\Omega;t\right)$}
  \psfrag{c1}{Case 1}
  \psfrag{c2}{Case 2}
  \psfrag{c3}{Case 3}
  \psfrag{c4}{Case 4}
  \psfrag{c5}{Case 5}
  \subfigure[Minimum nodal concentration.]{
    \includegraphics[clip, scale=0.32]{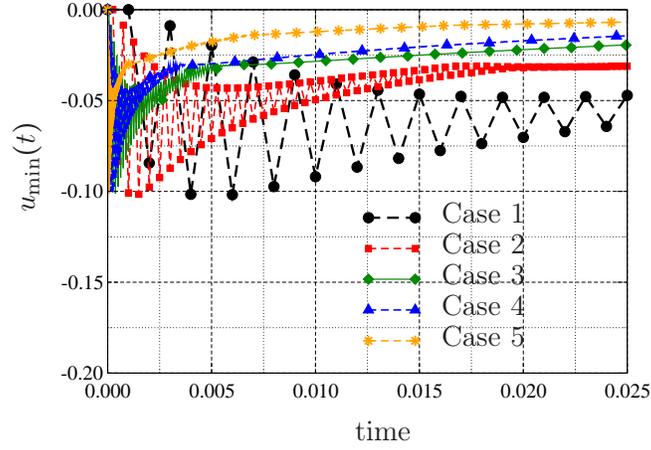}
  }
  \subfigure[Total amount of the dispersing chemical species with the 
    negative values.]{
    \includegraphics[clip, scale=0.32]{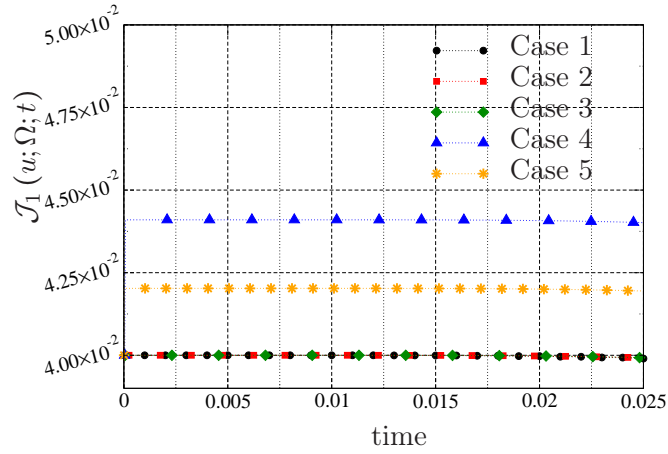}
  }
  \subfigure[Total amount of the dispersing chemical species \emph{without}
    the negative values.]{
    \includegraphics[clip, scale=0.32]{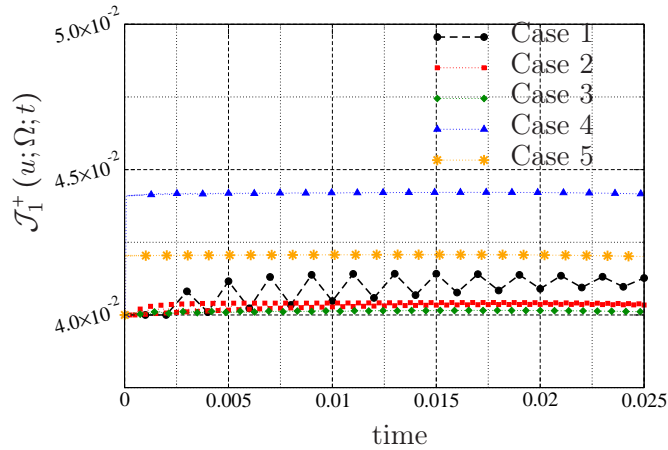}
  }
  \caption{\textsf{Two-dimensional problem with anisotropic 
  and heterogeneous diffusion tensor:}~The total 
    amounts of the diffusing species under the 
    single-relaxation-time LBM and the clipping 
    procedure are plotted against time. The 
    variation of the minimum concentration over 
    time is shown in the bottom figure. The 
    parameters in various cases are provided 
    in Table \ref{Tbl:H_W_2D_Anisotropic}.
    \label{Fig:H_W_2D_Anisotropic_plots_1}}
\end{figure}

\begin{figure}
  \centering
  \psfrag{time}{time}
  \psfrag{intPhi2}{$\mathcal{J}_{2}\left(u;\Omega;t\right)$}
  \psfrag{intPhi2pos}{$\mathcal{J}^{+}_{2}\left(u;\Omega;t\right)$}
  \psfrag{diff}{$\big| \mathcal{J}^{+}_{2}\left(u;\Omega;t\right) - 
  \mathcal{J}_{2}\left(u;\Omega;t\right)\big|$}
  \psfrag{c1}{Case 1}
  \psfrag{c2}{Case 2}
  \psfrag{c3}{Case 3}
  \psfrag{c4}{Case 4}
  \psfrag{c5}{Case 5}
  \subfigure[Negative values are \emph{not} clipped.]{
    \includegraphics[clip, scale=0.32]{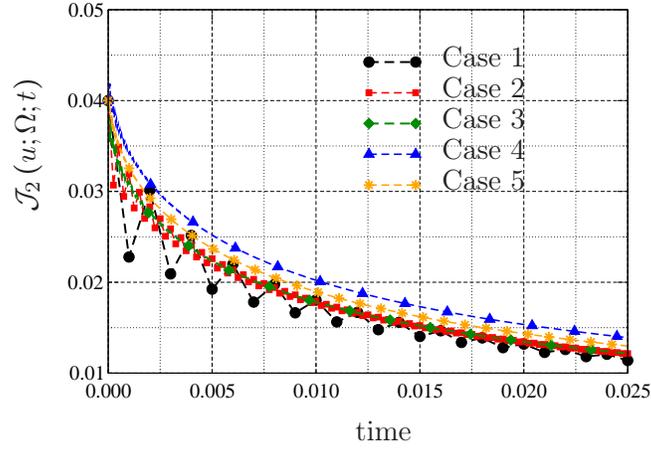}
  }
  \subfigure[Negative values are clipped.]{
    \includegraphics[clip, scale=0.32]{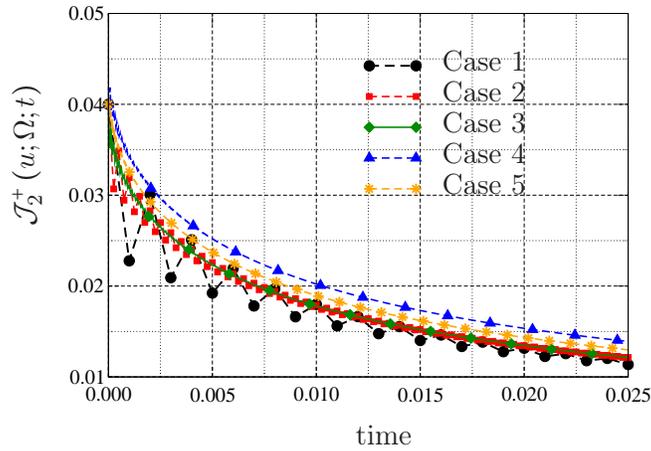}
  }
  \subfigure[Resulting difference from clipping the negative values.]{
    \includegraphics[clip, scale=0.32]{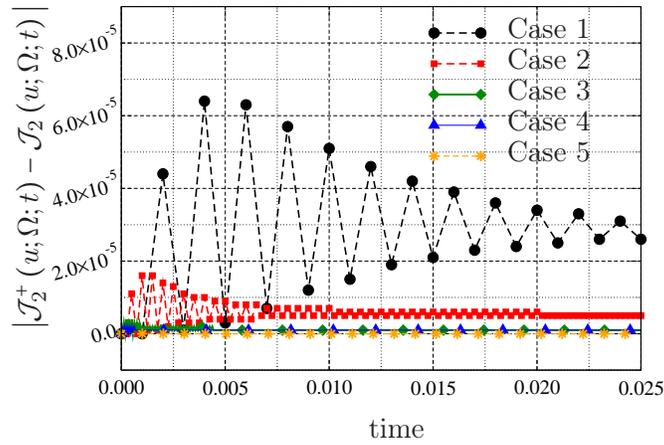}
  }
  \caption{\textsf{Two-dimensional problem with anisotropic 
      and heterogeneous diffusion tensor:}~This figure shows 
    the variation of various integrals defined in equation 
    \eqref{Eqn:Integralu} with respect to time under the 
    H-W method. The simulation parameters are provided in 
    Table \ref{Tbl:H_W_2D_Anisotropic}. \emph{It is evident 
      from the figure that the H-W method violated the decay 
      property.} \label{Fig:H_W_2D_Anisotropic_plots_2}}
\end{figure}

Consider the spatial domain to be $\Omega = (0,1) 
\times (0,1)$. We have taken the following anisotropic 
and heterogeneous diffusivity tensor:
\begin{align}
  \mathbf{D}\left( \mathrm{x}, \mathrm{y}\right) = \epsilon' 
  \mathbf{I}_{2} + \left[ \begin{array}{c c}
      \epsilon \mathrm{x}^2 + \mathrm{y}^2 & -\left( 1 - \epsilon\right)
      \mathrm{x}\mathrm{y} \\
      -\left( 1 - \epsilon\right)\mathrm{x}\mathrm{y} & \mathrm{x}^2 + 
      \epsilon \mathrm{y}^2
    \end{array} \right]
\end{align}
where $\epsilon \ll 1$ and $\epsilon' \ll 1$ are 
arbitrary constants, and $\mathbf{I}_2$ denotes the 
$2 \times 2$ identity matrix. For this numerical 
experiment, we have taken $\epsilon = 10^{-3}$ and 
$\epsilon' = 10^{-10}$. 
The prescribed concentration on the entire boundary 
is taken to be zero. The initial concentration is 
taken as: 
\begin{align}
  \mathrm{u}_0 \left( \mathrm{x}, \mathrm{y}\right) = 
  \left\{ \begin{array}{c c}
    1 	\quad 	&\left( \mathrm{x}, \mathrm{y}\right) \in 
    \left[0.4,0.6\right]\times\left[0.4,0.6\right] \\
    0   \quad 	&\mathrm{otherwise} \end{array} \right.
\end{align}
The time interval of interest is taken as $\mathcal{T} 
= 0.025$. The H-W method based on the $D2Q9$ lattice 
model is employed with the lattice velocity $c = \Delta 
x/\Delta t$. 
The discrete momenta are taken as:
\begin{align}
  \label{Eqn:D2Q9Lattice}
  \boldsymbol{e}_i^{\mathrm{T}} = \left\{ \begin{array}{ll}
    \left[ 0 , 0 \right] \quad &i = 0 \\
    c\left[ \mathrm{cos}\left( (i - 1)\pi/2\right), 
      \mathrm{sin}\left( (i - 1)\pi/2\right) \right]  \quad &i = 1,2,3,4 \\
    \sqrt{2} c \left[ \mathrm{cos}\left( (2i - 9)\pi/4\right), 
      \mathrm{sin}\left( (2i - 9)\pi/4\right)\right] \quad &i = 5,6,7,8
  \end{array}\right.
\end{align}
with the following weights: 
\begin{align}
  w_i = \left\{ \begin{array}{ll}
    4/9 \quad &i=0 \\
    1/9 \quad &i=1,2,3,4 \\
    1/36 \quad &i=5,6,7,8
  \end{array}\right.
\end{align}
The problem is solved using different choices of 
$\Delta x$ and $\Delta t$, which are provided in 
Table \ref{Tbl:H_W_2D_Anisotropic}. This table also 
provides insight on the performance of the H-W method. 
The spread of the lattice nodes that experience 
violation of non-negative constraint is shown in Figures 
\ref{Fig:H_W_2D_Anisotropic} and \ref{Fig:H_W_2D_Neg}. 
The following conclusions can be drawn from Figures 
\ref{Fig:H_W_2D_Anisotropic}--\ref{Fig:H_W_2D_Anisotropic_plots_2} 
and Table \ref{Tbl:H_W_2D_Anisotropic}: 
\begin{enumerate}[(a)]
\item The H-W multiple-relaxation-time method violates 
  the non-negative constraint when the diffusion is 
  anisotropic. 
\item As discussed earlier, the integral $\mathcal{J}_2$ 
  should decrease monotonically with time for pure 
  diffusion equations. However, the H-W method does 
  not respect the decay property. Clipping procedure 
  can eliminate the violation of the non-negative 
  constraint but does not eliminate the violation 
  of the decay property. However, it has been observed 
  that refining the discretization parameters (i.e., 
  $\Delta x$ and $\Delta t$) can improve the performance 
  of numerical solutions with respect to the decay property. 
\item In the case of the integral $\mathcal{J}_1$, 
  noticeable differences appear between the numerical 
  solution with negative values and the numerical 
  solution with the negative values clipped.  
\end{enumerate}

\begin{table}[h]
  \caption{\textsf{Two-dimensional problem with 
      anisotropic and heterogeneous diffusion 
      tensor:}~In this table, the number of 
      nodes that experience violations of the 
      non-negative constraint $N_{\mathrm{neg}}$,
      the minimum and maximum observed concentration 
      for various choices of $\Delta x$ and $\Delta t$
      are shown}. Note that in all the cases, 
    despite refining the lattice cell size, the 
    non-negative constraint is violated. The number 
    of nodes that have negative values for the 
    concentration is denoted by $N_{\mathrm{neg.}}$. 
    The total number of lattice nodes is denoted 
    by $N$. \label{Tbl:H_W_2D_Anisotropic}
  \begin{tabular}{| c | c | c | l | c | c | c | c |}
    \hline
    Case & $\Delta x$ & $\Delta t$ & $N_{\mathrm{neg}}/N \times 100$ 
    & $u_{\min}(\mathcal{T})$ & $u_{\max}(\mathcal{T})$ \\ \hline
    1 & $5.00\times10^{-2}$& $1.00\times10^{-3}$& $90/441      
    \times 100 = 20.41\%$ & -0.0472& 0.5379 \\ \hline
    2 & $2.50\times10^{-2}$& $2.50\times10^{-4}$& $762/1681    
    \times 100 = 45.33\%$ & -0.0311& 0.5576 \\ \hline
    3 & $1.25\times10^{-2}$& $6.25\times10^{-5}$& $2867/6561   
    \times100  = 43.70\%$ & -0.0193& 0.5768 \\ \hline
    4 & $1.00\times10^{-2}$& $4.00\times10^{-5}$& $4207/10201  
    \times100  = 41.24\%$ & -0.0144& 0.6188 \\ \hline
    5 & $5.00\times10^{-3}$& $1.00\times10^{-5}$& $11894/40401 
    \times100  = 29.44\%$ & -0.0068& 0.6021 \\ \hline
  \end{tabular}
\end{table}

\subsection{Fast bimolecular reaction}
Consider a simple chemical reaction of the form: 
\begin{align}
  n_A A + n_B B \rightarrow n_c C
\end{align}
where $A$, $B$ and $C$ are the participating chemical 
species; and $n_A$, $n_B$ and $n_C$ are their respective 
stoichiometry coefficients. We are interested in the 
fate of the product $C$ when the time-scale of the 
chemical reaction is much faster than that of the 
transport processes (i.e., diffusion and advection). 
A detailed description of this mathematical model can 
be found in \citep{Nakshatrala_Mudunuru_Valocchi_JCP_2013}, 
and will not be repeated here. However,  
the mentioned paper neglected advection in all their 
numerical examples, and the entire paper is devoted to 
the finite element method. 

\subsubsection{Fast bimolecular reaction in a porous medium} 
\begin{figure}
  \centering
  \psfrag{cpA}{$\mathrm{u}^{\mathrm{p}}_{A}$}
  \psfrag{cpB}{$\mathrm{u}^{\mathrm{p}}_{B}$}
  \psfrag{zt}{zero traction}
  \psfrag{vin}{$\mathbf{v}^{\mathrm{inlet}}$}
  \psfrag{nf}{zero flux}
  \psfrag{rb}{rigid solid}
  \psfrag{lx}{$L_{\mathrm{x}}$}
  \psfrag{ly}{$L_{\mathrm{y}}/2$}
  \includegraphics[clip,scale=0.75]{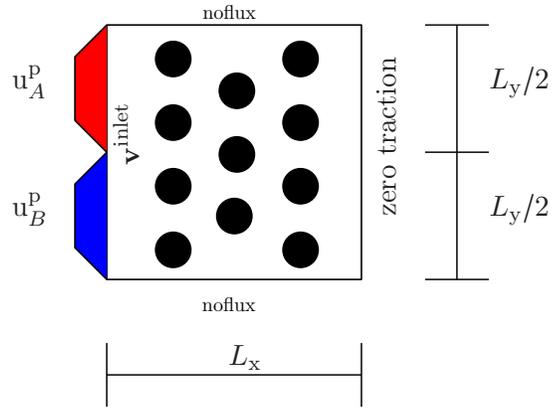}
  \caption{\textsf{Fast bimolecular reaction in a porous 
      medium:}~This figure provides a pictorial description 
    of the test problem. We have taken $L_{\mathrm{x}} = 1/2$ 
    and $L_{\mathrm{y}} = 2$ in the numerical experiment.
    \label{Fig:2DPorous_desfig}}
\end{figure}

We consider the combined effect of transport and 
chemical reactions in a porous medium. Problems 
of this type are frequently encountered in 
precipitation studies and ground-water hydrology 
\citep{Willigham_Werth_Valocchi_EST_2008_v42_p3185}. 
It needs to be mentioned that Willingham \emph{et al.} 
\citep{Willigham_Werth_Valocchi_EST_2008_v42_p3185} 
have performed numerical modeling but employed the 
finite volume method for the transport problem and 
LBM for the flow problem. We will investigate whether 
the non-negative constraint will be violated in the 
numerical simulations of such situations under LBM 
for both flow and transport problems. 
In the first numerical experiment, we shall assume 
the diffusivity tensor to be isotropic. The advection 
velocity will be determined by solving the incompressible 
Navier-Stokes equations in the pore structure (see Figure 
\ref{Fig:2DPorous_desfig}). That is, obtain $\mathbf{v}
(\mathbf{x},\mathrm{t})$ and $\mathrm{p}(\mathbf{x},
\mathrm{t})$ by solving: 
\begin{subequations}
  \begin{align}
    &\rho \left( \frac{\partial \mathbf{v}}{\partial \mathrm{t}} + 
    \mathbf{v} \cdot \mathrm{grad}\left[ \mathbf{v} \right]\right) = 
    -\mathrm{grad}\left[ \mathrm{p} \right] + 
    \mu \mathrm{div} \left[ \mathrm{grad}\left[ \mathbf{v}\right] + 
      \mathrm{grad}\left[ \mathbf{v}\right]^{\mathrm{T}}\right] \\
    &\mathrm{div}\left[ \mathbf{v} \right] = 0
  \end{align}
\end{subequations}
where $\mathbf{v} = \mathbf{v}(\mathbf{x},\mathrm{t})$ 
is the velocity, $\mathrm{p} = \mathrm{p}(\mathbf{x}, 
\mathrm{t})$ is the pressure, $\rho$ is the density of 
the fluid, and $\mu$ is the coefficient of viscosity. 
The density is taken as $\rho = 1$, and the viscosity of the 
fluid is taken to be $\mu = 10^{-2}$. The components of the 
inlet velocity are $\mathrm{v}_{\mathrm{x}}^{\mathrm{inlet}} = 1$ 
and $\mathrm{v}_{\mathrm{y}}^{\mathrm{inlet}} = 0$. Diffusion 
coefficient of all the chemical species is taken as 
$\mathrm{D} = 10^{-2}$. A pictorial description of the 
problem is given in Figure \ref{Fig:2DPorous_desfig}. 
The prescribed concentrations for the chemical species 
$A$ and $B$ are $\mathrm{u}_{A}^{\mathrm{p}} = 1$ and 
$\mathrm{u}_{B}^{\mathrm{p}} = 1$. The stoichiometry 
coefficients are taken as $n_A = 1$, $n_B = 2$, and 
$n_C = 1$. The dimensions of the computational domain 
are $L_{\mathrm{x}} = 1/2$ and $L_{\mathrm{y}} = 2$. The 
time interval of interest is taken as 
$\mathcal{T} = 0.5$.

\begin{figure}
  \centering
  \subfigure[$u_A$ at $t=0.125$]{
    \includegraphics[clip,scale=0.43]{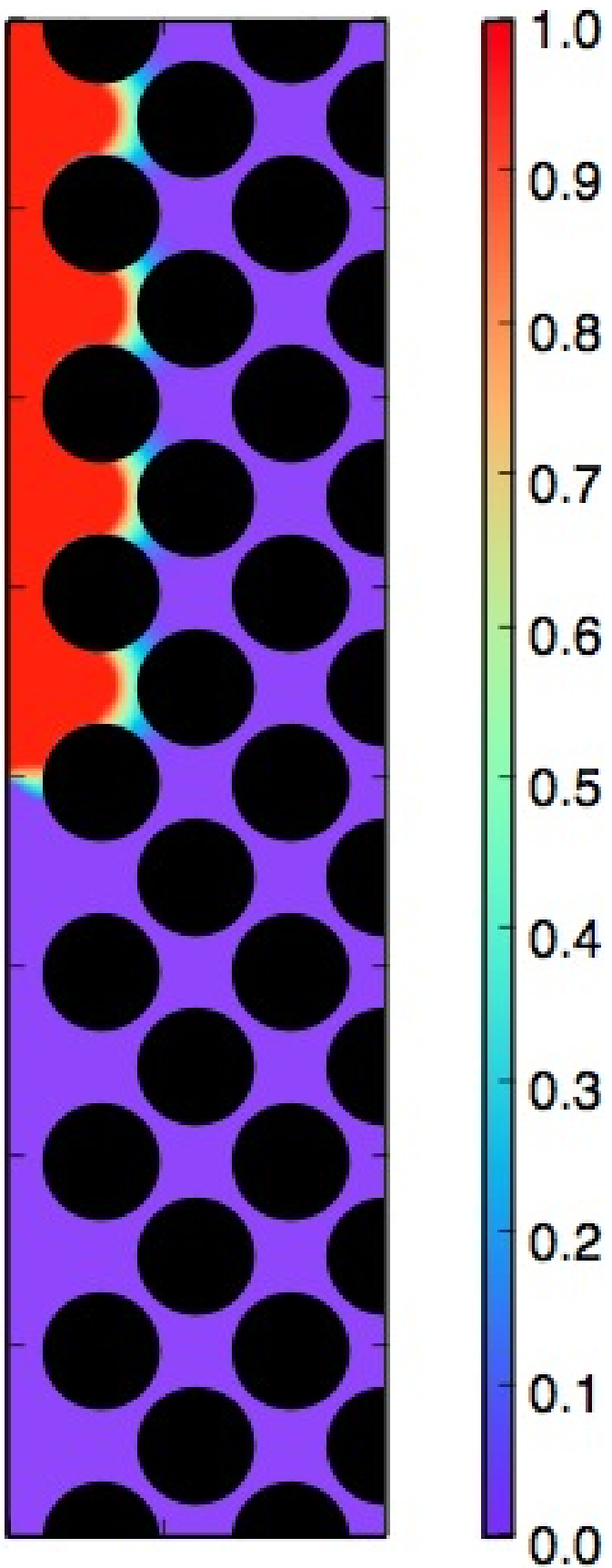}}
  \subfigure[$u_A$ at $t=0.250$]{
    \includegraphics[clip,scale=0.43]{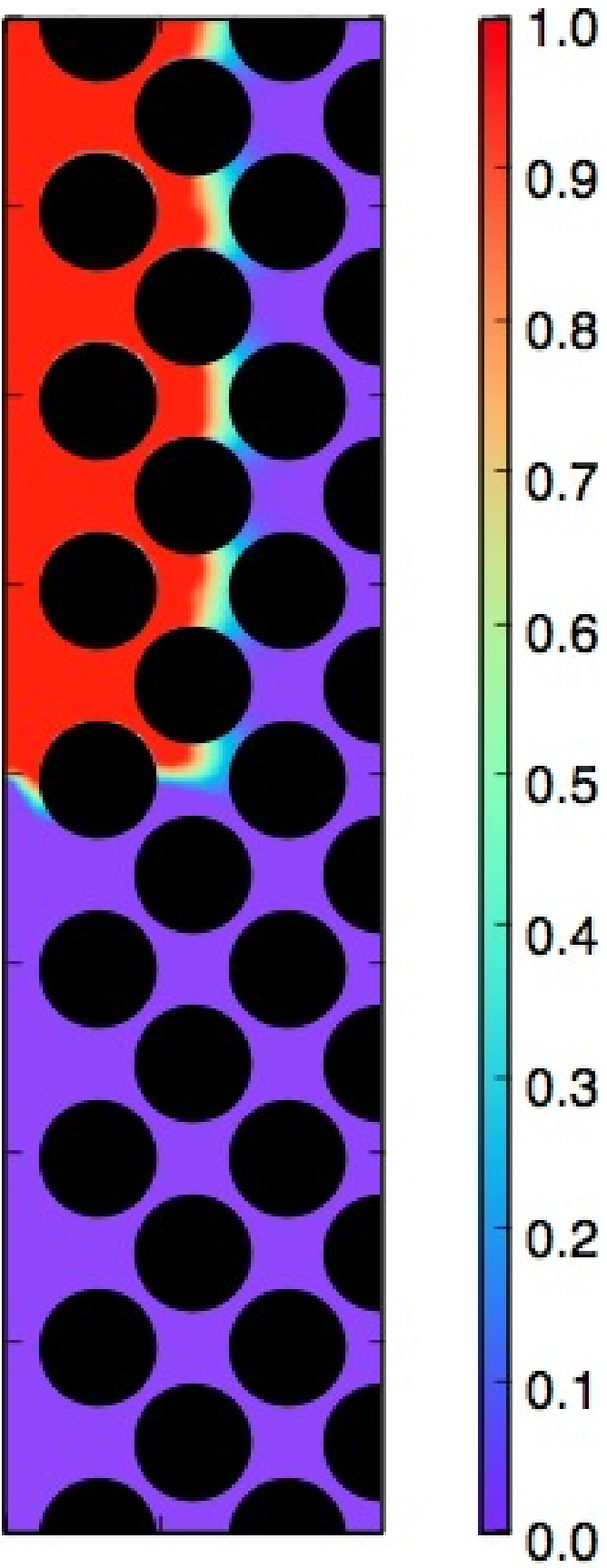}}
  \subfigure[$u_A$ at $t=0.375$]{
    \includegraphics[clip,scale=0.43]{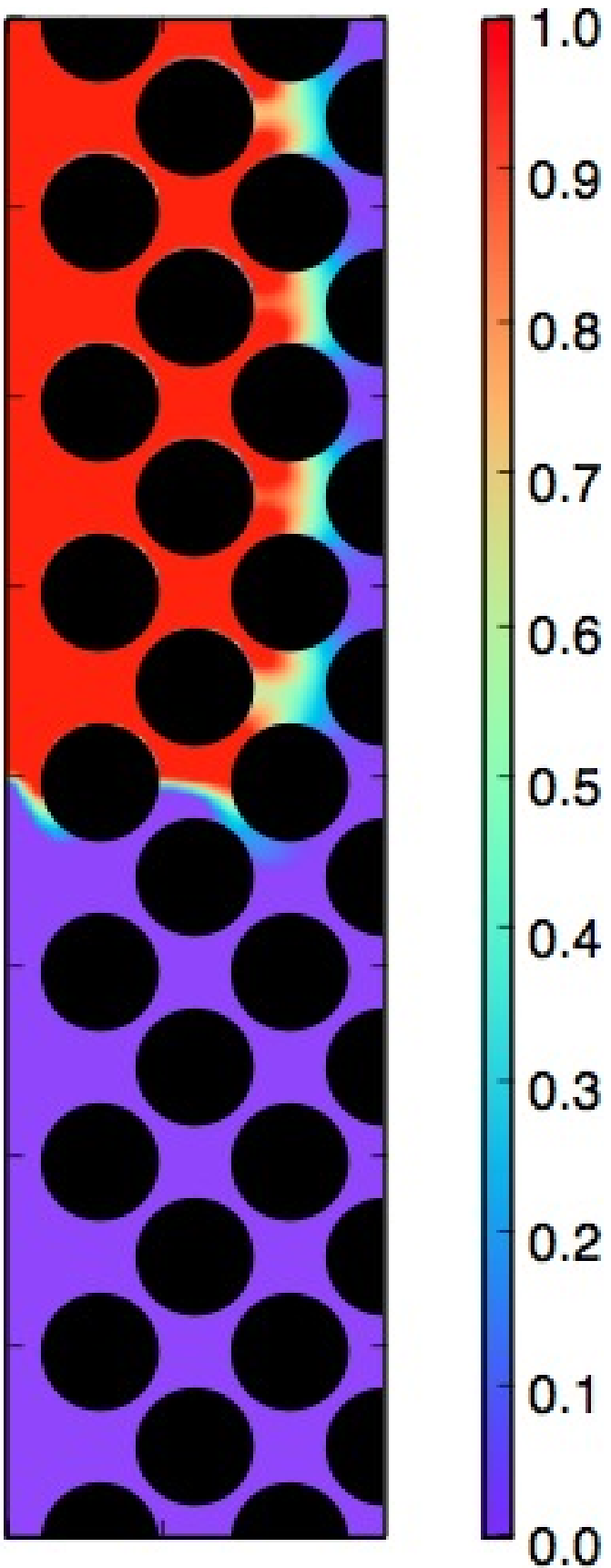}}
  \subfigure[$u_A$ at $t=0.500$]{
    \includegraphics[clip,scale=0.43]{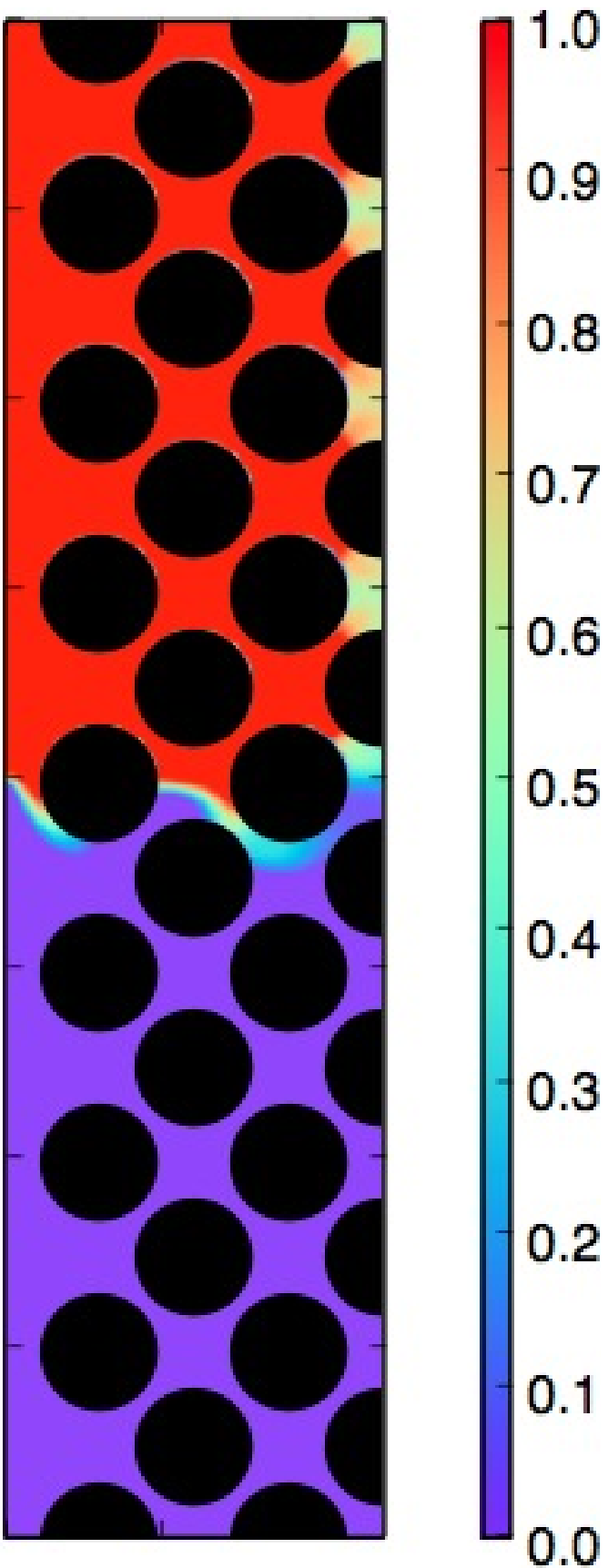}}
  \hspace{0.4in}
  \subfigure[$u_B$ at $t=0.125$]{
    \includegraphics[clip,scale=0.43]{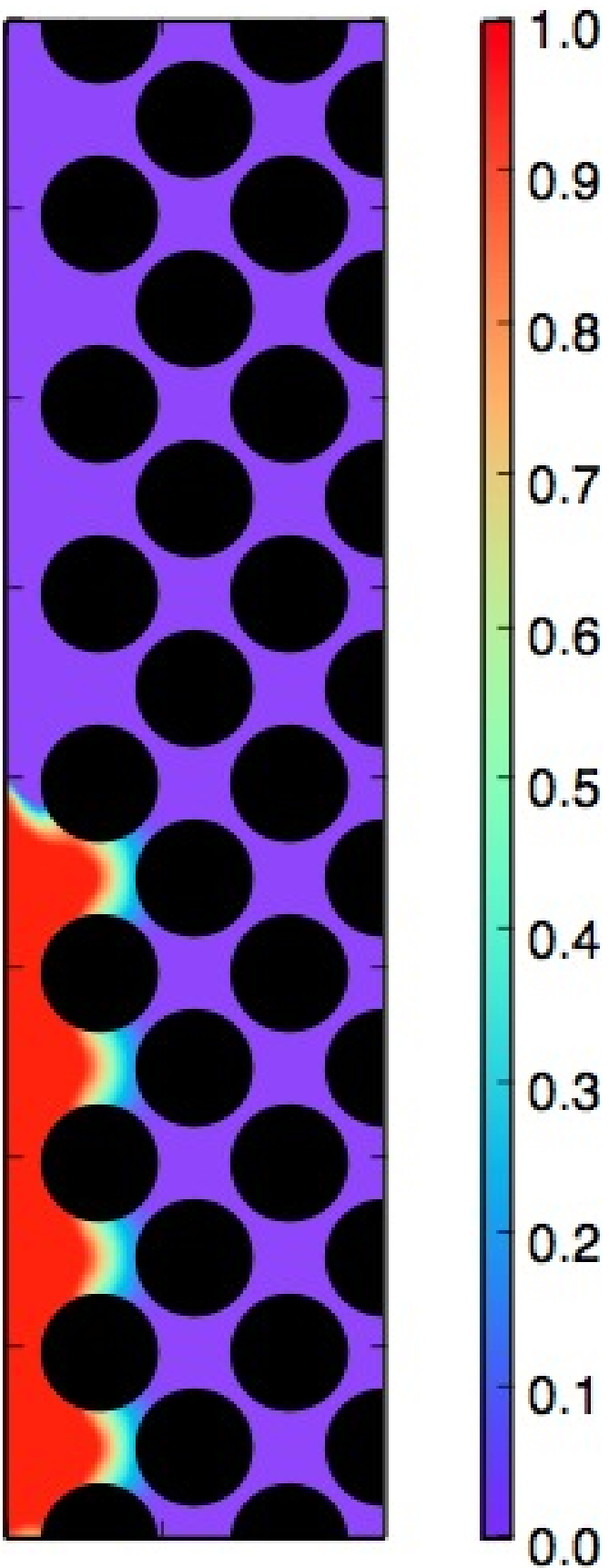}}
  \subfigure[$u_B$ at $t=0.250$]{
    \includegraphics[clip,scale=0.43]{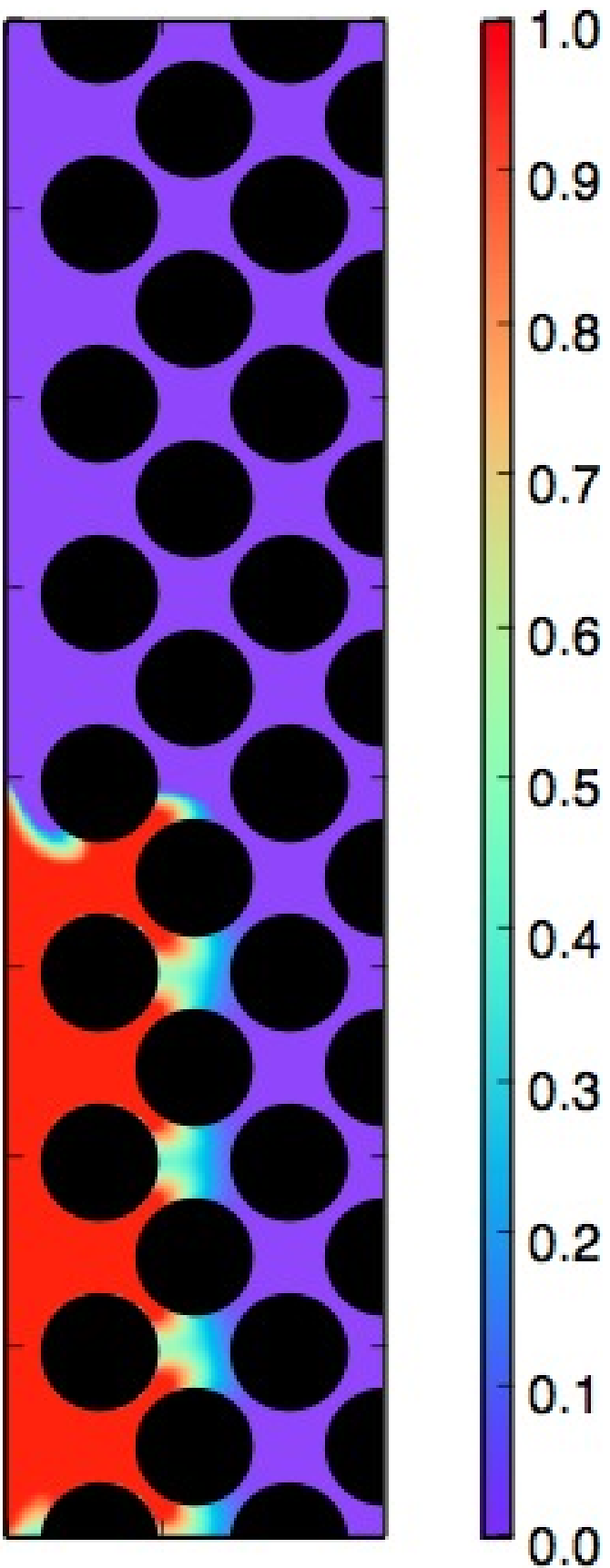}}
  \subfigure[$u_B$ at $t=0.375$]{
    \includegraphics[clip,scale=0.43]{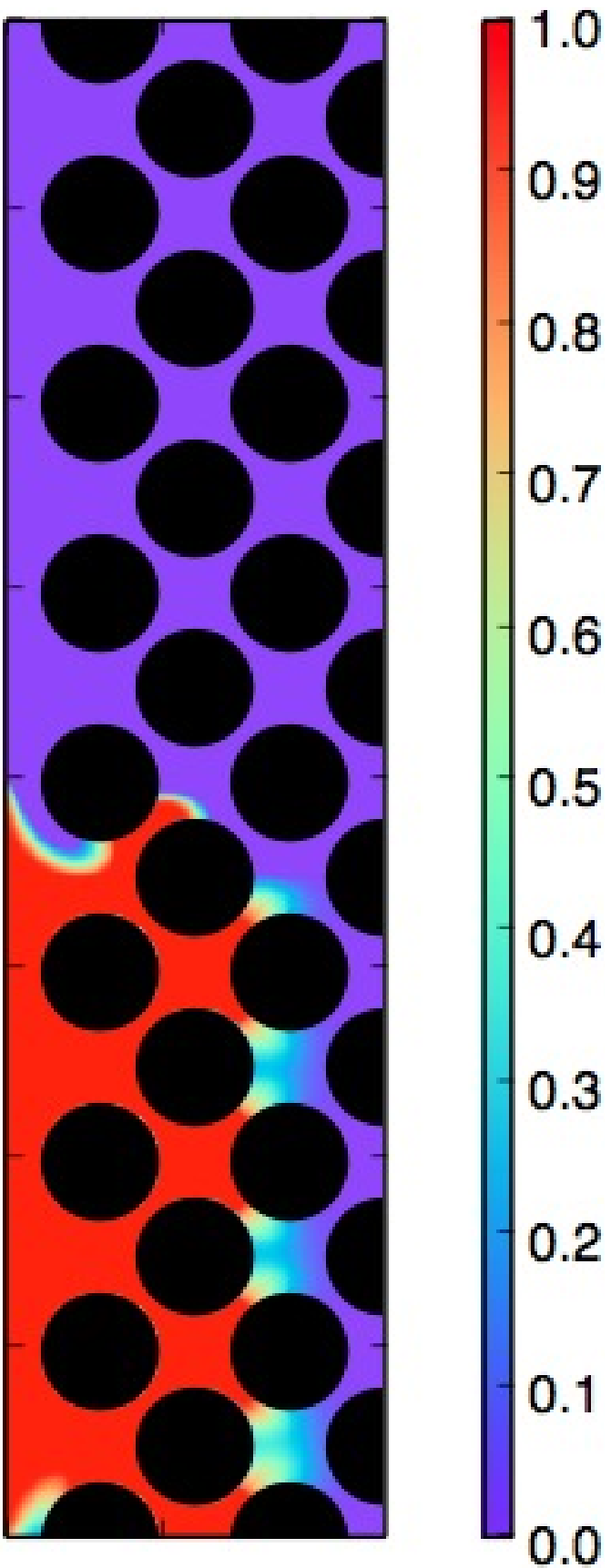}}
  \subfigure[$u_B$ at $t=0.500$]{
    \includegraphics[clip,scale=0.43]{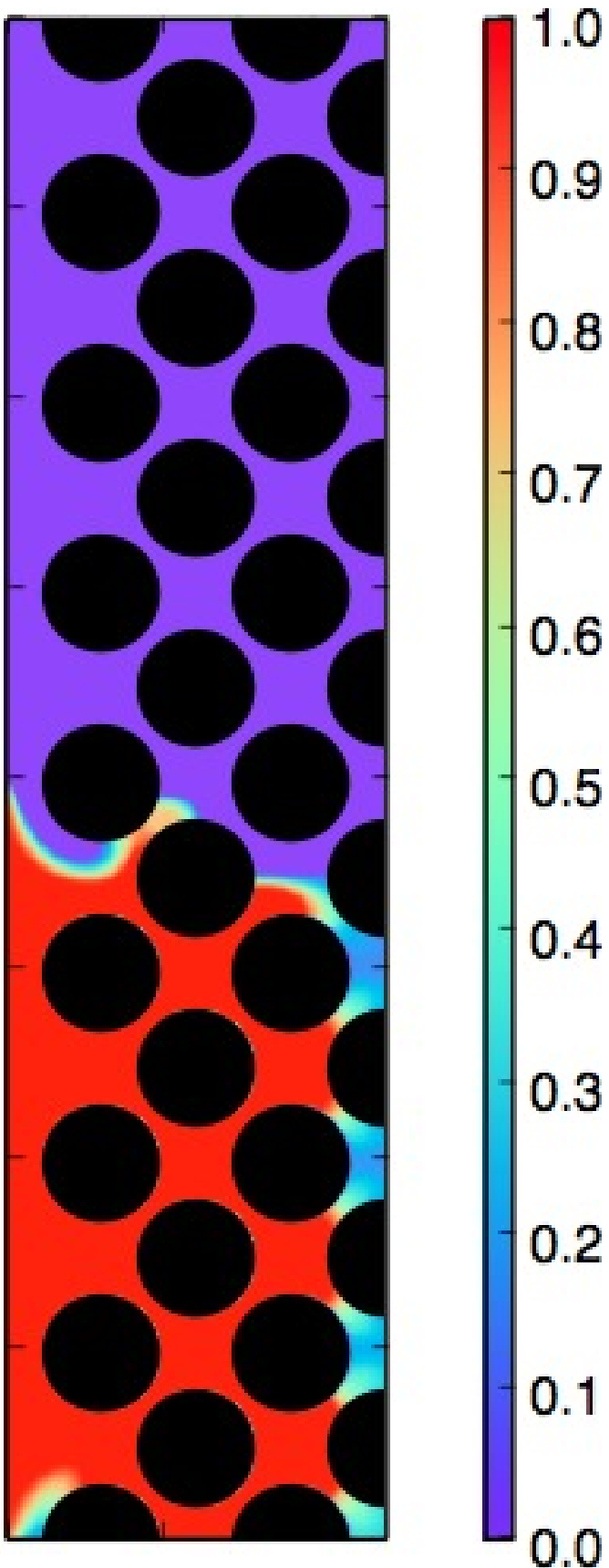}}
  \caption{\textsf{Fast bimolecular reaction in a porous 
      medium:}~This figure shows the concentration profiles 
    of the reactants $A$ and $B$ under the single-relation-time 
    method using the $D2Q9$ lattice model. The diffusivity 
    tensor is isotropic, and no negative values for the 
    concentration are observed for the reactants. 
    \label{Fig:2DPorous_A_B}}
\end{figure}

\begin{figure}
  \centering
  \subfigure[$u_C$ at $t = 0.125$]{
    \includegraphics[clip,scale=0.43]{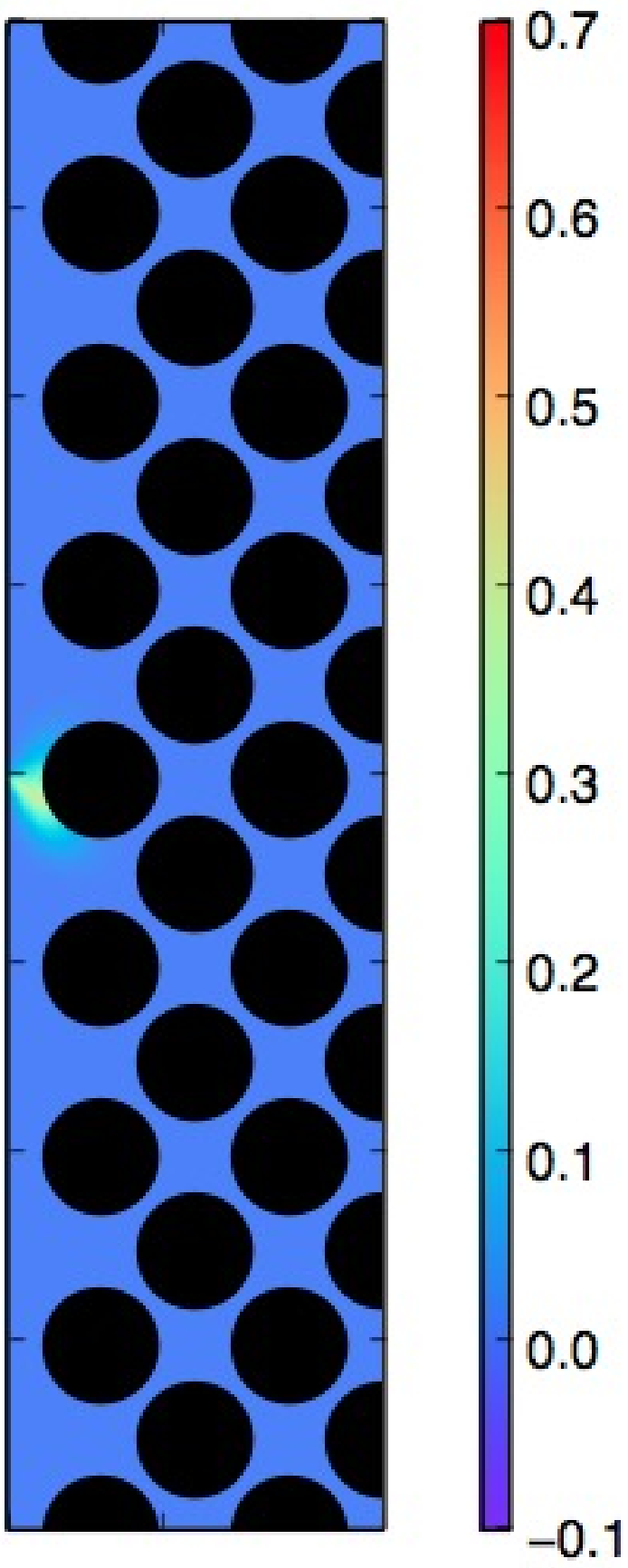}}
  \subfigure[$u_C$ at $t = 0.250$]{
    \includegraphics[clip,scale=0.43]{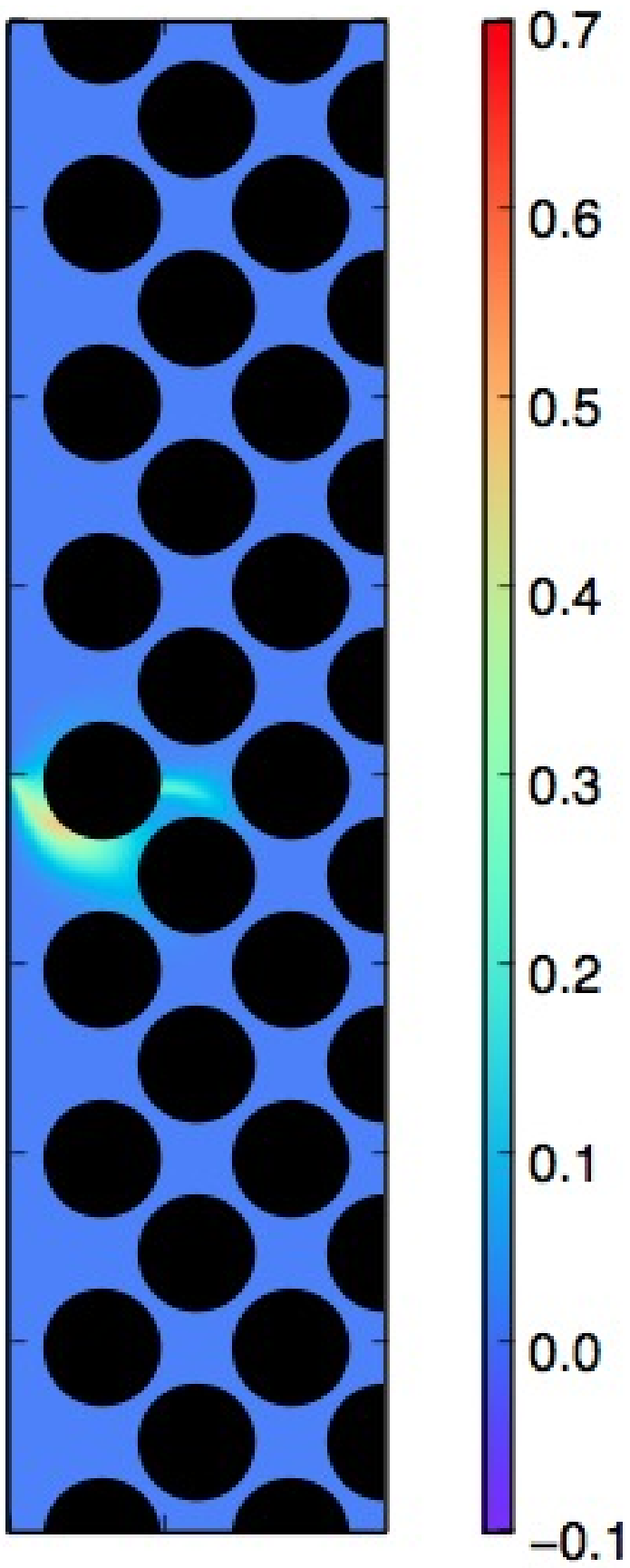}}
  \subfigure[$u_C$ at $t = 0.375$]{
    \includegraphics[clip,scale=0.43]{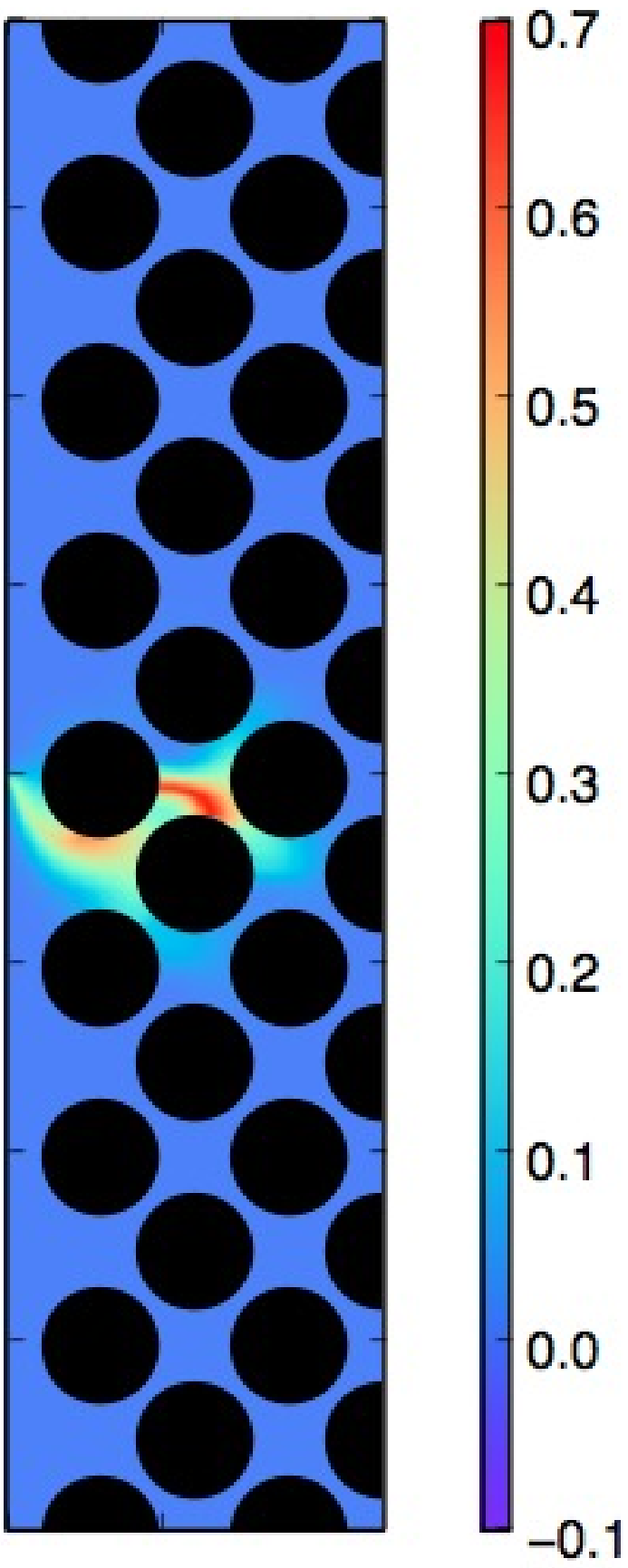}}
  \subfigure[$u_C$ at $t = 0.500$]{
    \includegraphics[clip,scale=0.43]{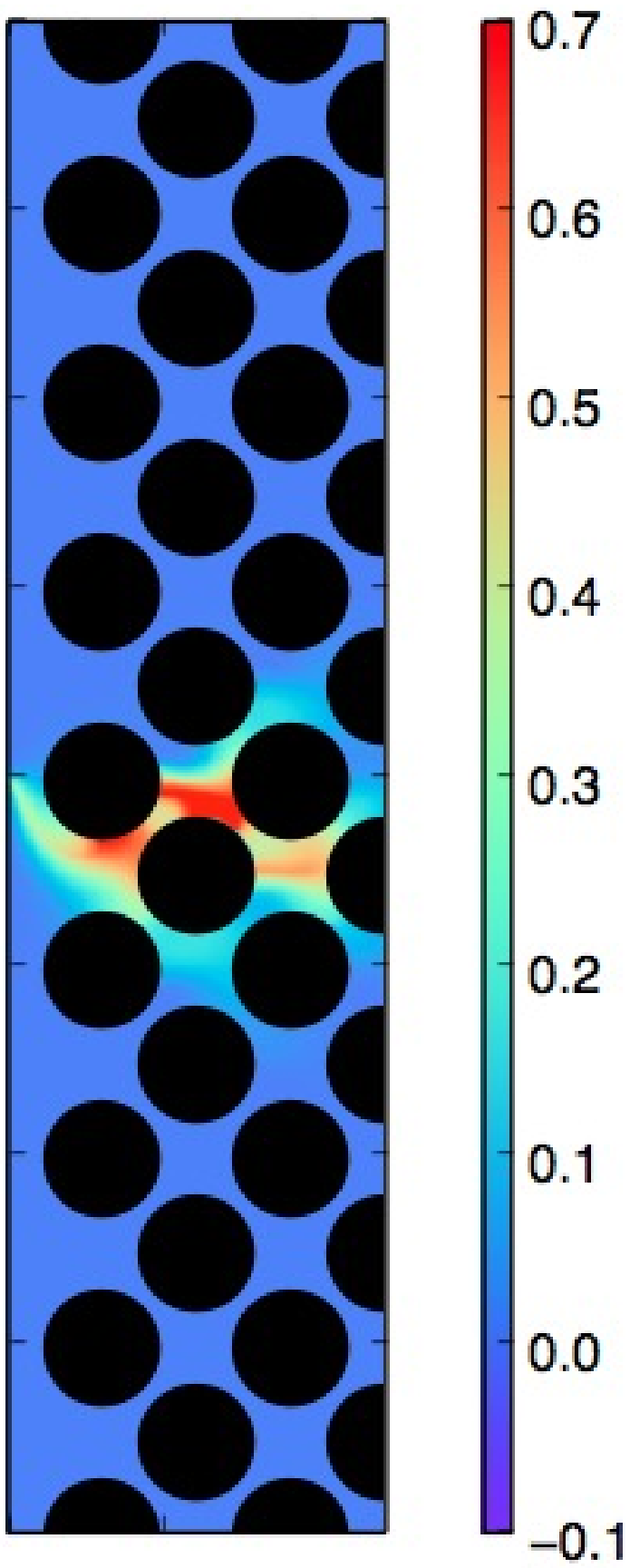}}
  \caption{\textsf{Fast bimolecular reaction in a porous 
      medium:}~This figure shows the concentration profiles 
    of the product $C$. \label{Fig:2DPorous_C}}
\end{figure}

\begin{figure}
\centering
\psfrag{minc}{$u_{\min}(t)$}
\psfrag{time}{time $t$}
\psfrag{cA}{Species $A$}
\psfrag{cB}{Species $B$}
\psfrag{cC1}{Species $C$, Case 1}
\psfrag{cC2}{Species $C$, Case 2}
\includegraphics[clip,scale=0.4]{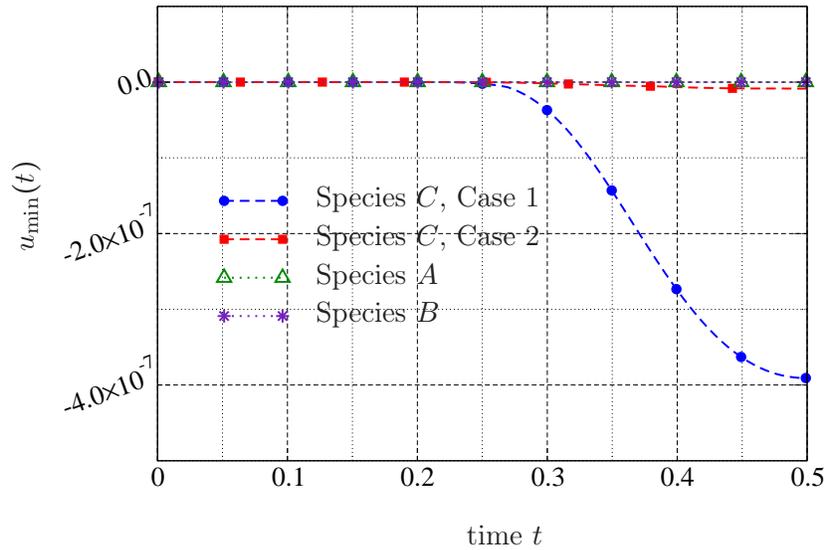}
\caption{\textsf{Fast bimolecular reaction in a porous 
    medium:}~This figure shows the minimum observed 
  values of the participating chemical species. 
  Concentrations of the chemical species $A$ and 
  $B$ did not violate the non-negative constraint. 
  However, small negative concentration of the 
  product $C$ are observed. These violations did 
  not diminish with time. But, the magnitude of 
  the violation can be decreased by refining the 
  discretization parameters. \label{Fig:2DPorous_min}}
\end{figure}

Table \ref{Tbl:2DPorous} provides the choices of 
$\Delta x$ and $\Delta t$ employed in the numerical 
simulation. We employed the single-relaxation-time 
lattice Boltzmann method using the $D2Q9$ lattice 
model for both flow and transport problems (i.e., 
see equation \eqref{Eqn:D2Q9Lattice}). 
The concentration of the chemical species $A$, $B$ and $C$ are
shown in Figures \ref{Fig:2DPorous_A_B}--\ref{Fig:2DPorous_C}. 
No negative values for the concentration are observed for the 
reactants $A$ and $B$, which is the result of equations (3.6a) 
and (3.6b) in 
\citep{Nakshatrala_Mudunuru_Valocchi_JCP_2013}. 
However, the concentration of the product $C$ exhibited 
negative values, which is shown in Figure \ref{Fig:2DPorous_min}. 
These negative values are very small in comparison with the 
maximum concentration in the domain. We also observed 
that refining the discretization parameters, $\Delta t$ and 
$\Delta x$, can reduce the violation of the non-negative 
constraint, unlike the previous problems with anisotropic 
diffusion tensor. 

\begin{table}
  \caption{\textsf{Fast bimolecular reaction in a porous 
      medium:}~Discretization parameters used in the 
    numerical experiment.}
  \label{Tbl:2DPorous}
  \begin{tabular}{| c | c | c |} \hline
    Case 	& 	$\Delta x$ 				& 	$\Delta t$ 			  \\ \hline
    1		&	 $6.25 \times 10^{-3}$ 	& 	$9.75 \times 10^{-4}$ \\ \hline
    2		&	 $5.00 \times 10^{-3}$	&	$6.25 \times 10^{-4}$ \\ \hline
  \end{tabular}
\end{table}

\subsubsection{Fast bimolecular reaction in an anisotropic 
  and heterogeneous medium} 
\begin{figure}
  \centering
  \psfrag{nf}{zero flux}
  \psfrag{lx}{$L_x$}
  \psfrag{ly}{$L_y/2$}
  \psfrag{cpA}{$\mathrm{u}_{A}^{\mathrm{p}}$}
  \psfrag{cpB}{$\mathrm{u}_{B}^{\mathrm{p}}$}
  \psfrag{c0A}{$\mathrm{u}_{A}(\mathbf{x},\mathrm{t}=0)=0$}
  \psfrag{c0B}{$\mathrm{u}_{B}(\mathbf{x},\mathrm{t}=0)=0$}
  \psfrag{c0C}{$\mathrm{u}_{C}(\mathbf{x},\mathrm{t}=0)=0$}
  \includegraphics[clip,scale=1]{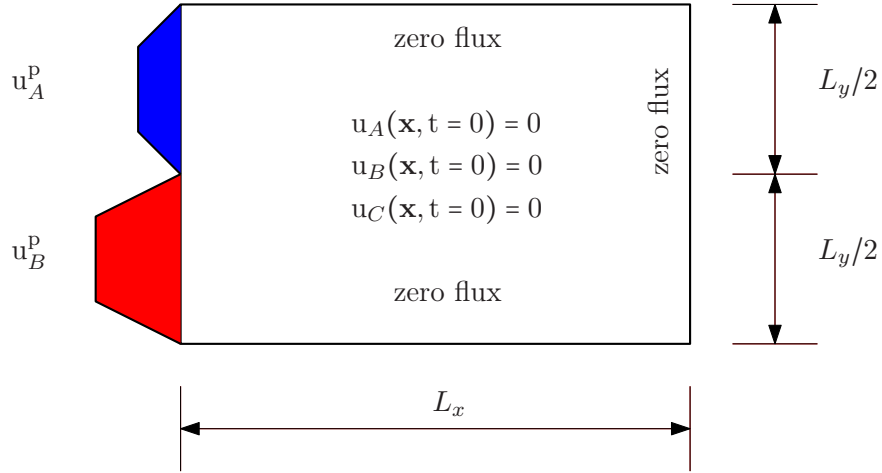}
  \caption{\textsf{Fast bimolecular reaction in 
      anisotropic and heterogeneous medium:}~This 
    figure provides a pictorial description of the 
    test problem. The reactants $A$ and $B$ undergo 
    transport (i.e., both advection and diffusion) 
    and reacts to give product $C$, which in turn 
    gets transported. We have taken $L_x = 2$ and 
    $L_y = 1$ in the numerical experiment. 
    \label{Fig:2DChemReact}}
\end{figure}

\begin{figure}
  \centering
  \includegraphics[scale=0.5,clip]{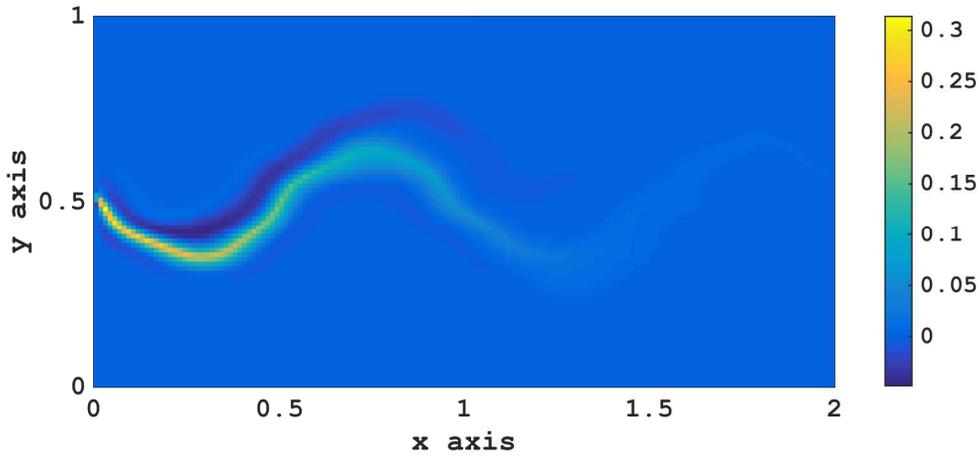}
  \caption{\textsf{Fast bimolecular reaction in 
      anisotropic and heterogeneous medium:}~This 
    figure shows the concentration profile of the 
    product $C$ at $t = \mathcal{T} = 0.25$ under the H-W 
    method. We have taken $\Delta x = 1.25 \times 
    10^{-2}$ and $\Delta t = 1.56 \times 10^{-6}$ 
    (see Case 3 in Table \ref{Tbl:2DChemReact}).
    \label{Fig:2DChemReact_uC}}
\end{figure}

\begin{figure}
  \centering
  \subfigure[Case 1]{
    \includegraphics[scale=0.45,clip]{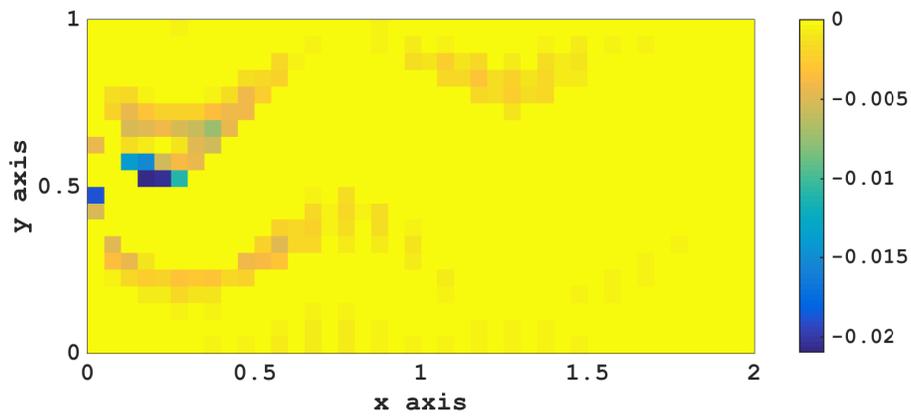}}
  \subfigure[Case 2]{
    \includegraphics[scale=0.45,clip]{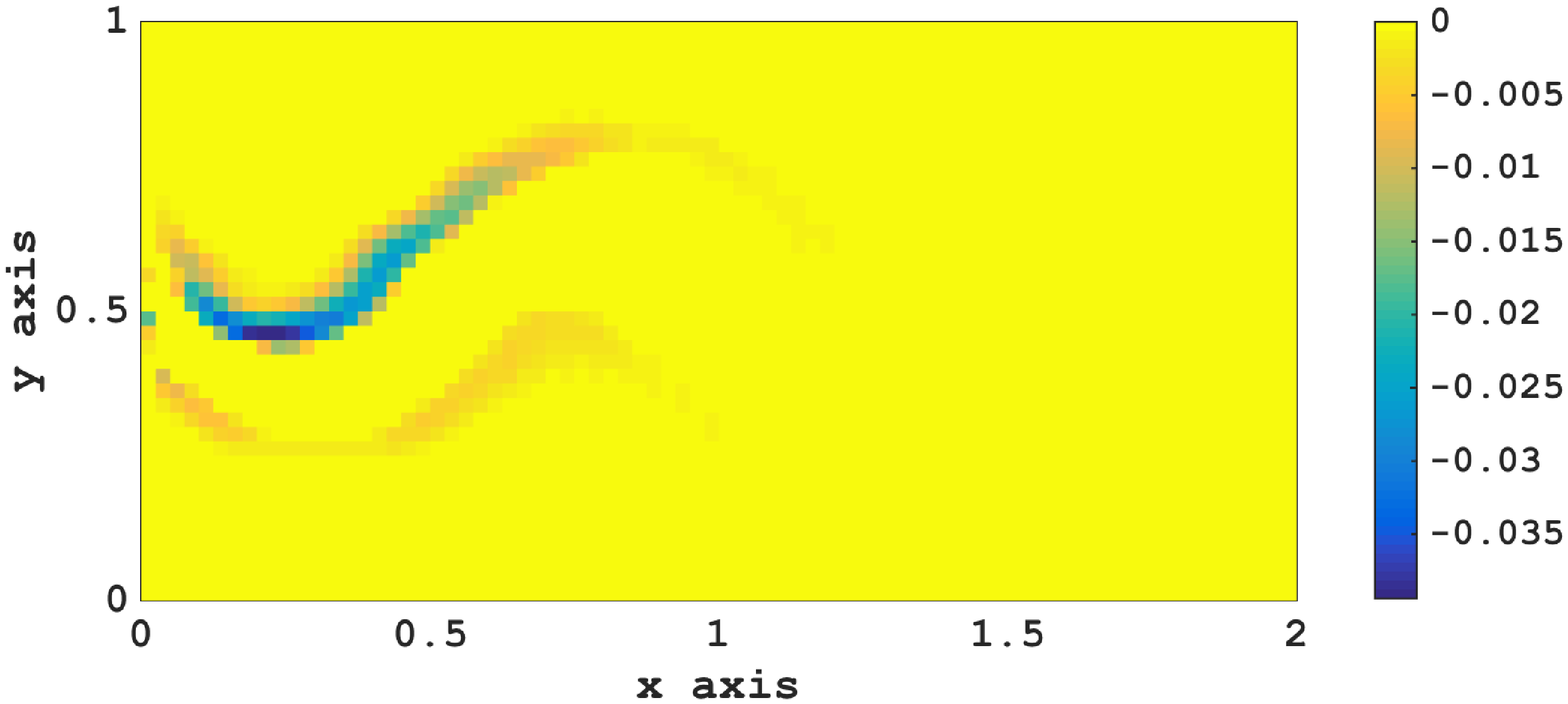}}
  \subfigure[Case 3]{
    \includegraphics[scale=0.45,clip]{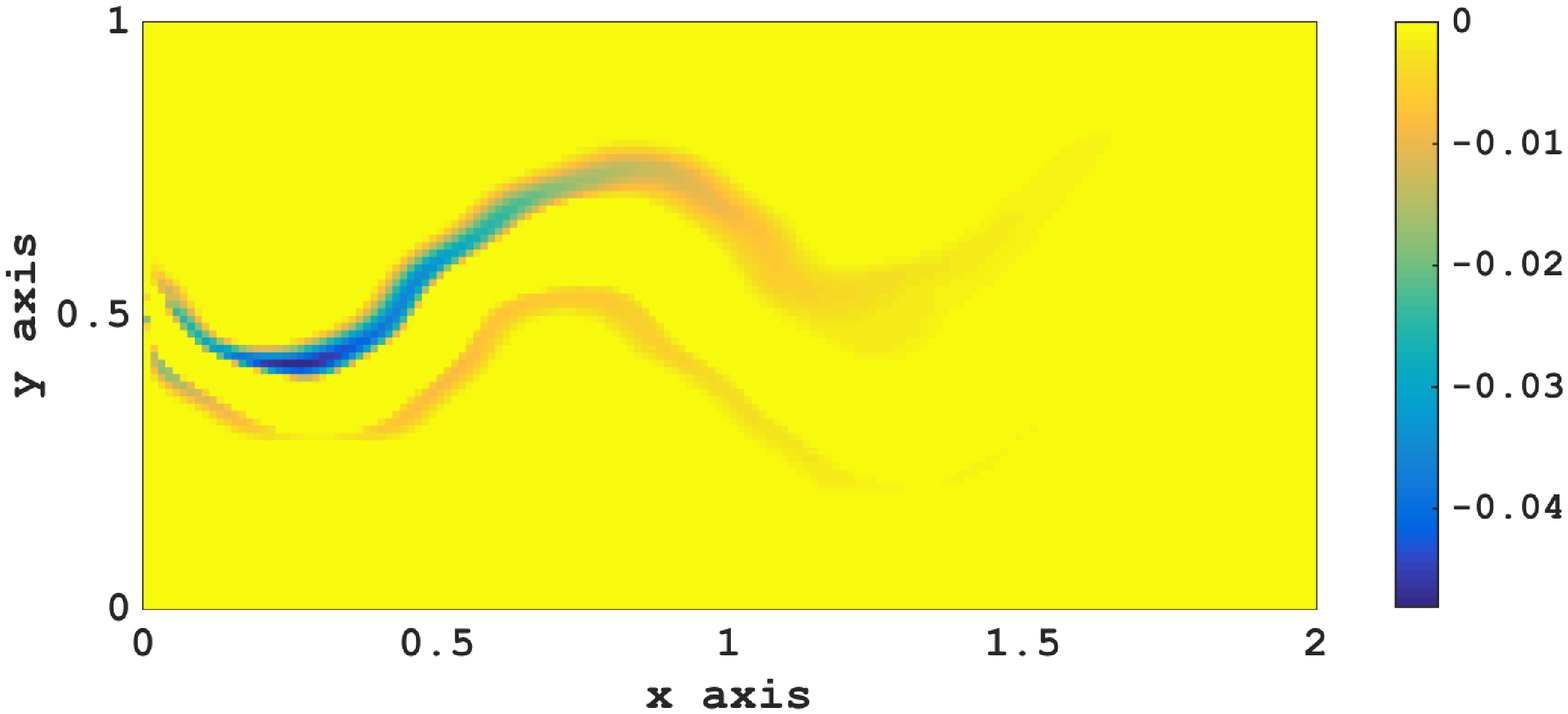}}
  \caption{\textsf{Fast bimolecular reaction in anisotropic 
      and heterogeneous medium:}~This figure shows the 
    regions where the concentration of the product 
    $C$ is negative at $t = 0.25$ under the H-W 
    method. The violations appear at large number of 
    lattice nodes, and are spread across the computational 
    domain.  
    \label{Fig:2DChemReact_negC}}
\end{figure}

\begin{figure}
  \centering
  \psfrag{time}{$t$}
  \psfrag{minc}{$u_{\min}\left( t \right)$}
  \psfrag{c1}{Case 1}
  \psfrag{c2}{Case 2}
  \psfrag{c3}{Case 3}
  \includegraphics[scale=0.4,clip]{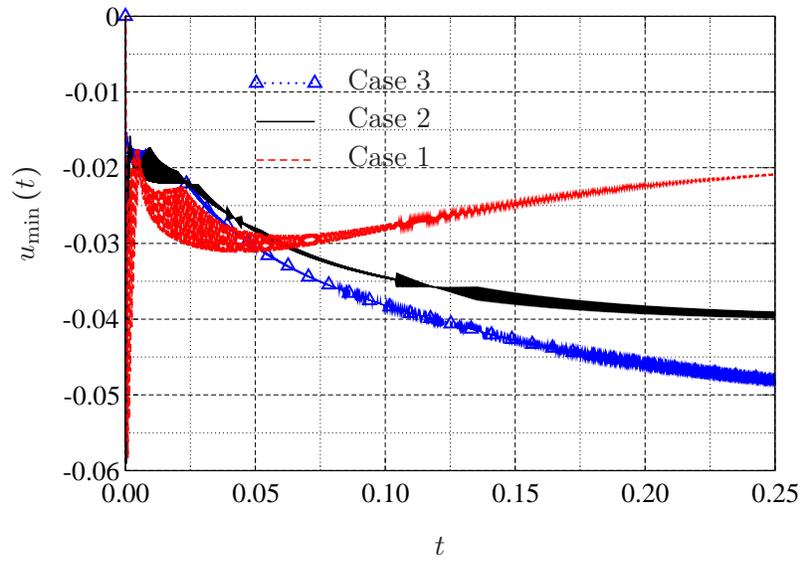}
  \caption{\textsf{Fast bimolecular reaction in anisotropic 
      and heterogeneous medium:}~ The minimum concentration 
    of the product $C$ is plotted against time for various 
    cases whose simulation parameters are provided in Table 
    \ref{Tbl:2DPorous}. \emph{This figure clearly illustrates 
      that refining the discretization parameters does not 
      eliminate the violation of the non-negative constraint 
      under the H-W method.} \label{Fig:2DChemReact_MinC}}
\end{figure}

We will consider the domain given in Figure 
\ref{Fig:2DChemReact}. Dimensions of the domain 
are $L_\mathrm{x} = 2$ and $L_\mathrm{y} = 1$. The 
advection velocity will be derived from the 
following stream function:
\begin{align}
	\psi \left( \mathrm{x}, \mathrm{y}\right) = -\mathrm{y}
	-\sum_{k = 1}^{3} \alpha_k \mathrm{cos}\left( 
	\frac{p_k \pi \mathrm{x}}{L_{\mathrm{x}}} - 
	\frac{\pi}{2}\right) \mathrm{sin}\left( 
	\frac{q_k \pi \mathrm{y}}{L_{\mathrm{y}}}\right)
\end{align}
where the parameters are given by 
\begin{align}
  (p_1,p_2,p_3) = (4.0,5.0,10.0), \; 
  (q_1,q_2,q_3) = (1.0,5.0,10.0), \; 
  (\alpha_1,\alpha_2,\alpha_3) = (0.08,0.02,0.01) 
\end{align}
The components of the advection velocity are calculated as: 
\begin{align}
	\mathrm{v}_{\mathrm{x}}\left( \mathrm{x}, \mathrm{y}\right) = -
	\frac{\partial \psi \left( \mathrm{x}, \mathrm{y}\right)}
	{\partial \mathrm{y}}, \quad
	\mathrm{v}_{\mathrm{y}}\left( \mathrm{x}, \mathrm{y}\right) = +
	\frac{\partial \psi \left( \mathrm{x}, \mathrm{y}\right)}
	{\partial \mathrm{x}} \quad
\end{align}
The dispersion tensor is taken as: 
\begin{align}
  \mathbf{D}\left( \mathrm{x}, \mathrm{y}\right) = 
  10^{-5} \mathbf{I} + \beta_{\mathrm{T}} \left\| \mathbf{v}\right\| 
  \mathbf{I} + \left( \beta_{\mathrm{L}} - \beta_{\mathrm{T}}\right) 
  \frac{\mathbf{v} \otimes \mathbf{v}}{\left\| \mathbf{v} \right\|}
\end{align}
where $\otimes$ denotes the tensor product, $\mathbf{I}$ is 
the $2 \times 2$ identity tensor, $\left\| \cdot \right\|$ 
is the 2-norm, $\beta_{\mathrm{T}} = 10^{-4}$, and $\beta_{\mathrm{L}} 
= 1$. The prescribed concentrations are $\mathrm{u}^{\mathrm{p}}_A = 
\mathrm{u}^{\mathrm{p}}_B = 1$, and the stoichiometry coefficients are
$n_A = 1$, $n_B = 2$ and $n_C = 1$. The time interval of interest 
is $\mathcal{T} = 0.25$. We employed the $D2Q9$ lattice model 
using the H-W method, and obtained the numerical solution for 
the fate of the product $C$. Discretization parameters for 
various cases are given in Table \ref{Tbl:2DChemReact}. 
\emph{
Figures \ref{Fig:2DChemReact_uC}--\ref{Fig:2DChemReact_MinC} 
clearly show that the H-W method violated the non-negative 
constraint, and the violations did not vanish either with time 
or with refinement of the discretization parameters.}

\begin{table}
\caption{\textsf{Fast bimolecular reaction in anisotropic and 
heterogeneous medium:}~Different discretization parameters and 
violation of the non-negative constraint.}
\label{Tbl:2DChemReact}
\begin{tabular}{| c | c | c | c | c | l |}
\hline
Case 	& 	$\Delta x$ 					& 	$\Delta t$ 	& $u_{\min}\left( \mathcal{T}\right)$ & $u_{\max} \left( \mathcal{T} \right)$ &
$N_{\mathrm{neg}}/N \times 100$ \\ \hline
1		& 		$5.00 \times 10^{-2}$	& 	$2.50 \times 10^{-4}$	& -0.0209 & 0.2704 & $275/861 \times 100 = 31.94\%$\\ \hline
2 		&		$2.50 \times 10^{-2}$	&	$6.25 \times 10^{-5}$	& -0.0394 & 0.3072 & $1283/3321 \times 100 = 38.63\%$\\ \hline
3 		& 		$1.25 \times 10^{-2}$   & 	$1.56 \times 10^{-6}$ 	& -0.0481 & 0.3136 & $5703/13041 \times 100 = 43.73\%$\\ \hline
\end{tabular}
\end{table}

\section{A THEORETICAL ANALYSIS}
\label{Sec:LBM_S5_Theoretical}
In this section, we will provide a simple criterion in 
terms of the discretization parameters to satisfy the 
non-negative constraint for one-dimensional problems. 
We will also limit our scope to pure diffusion equations 
(i.e., $\mathbf{v}(\mathbf{x},\mathrm{t}) = \mathbf{0})$, 
and $\partial \Omega = \Gamma^{\mathrm{D}}$. We will restrict 
the analysis to the $D1Q3$ lattice model. We initialize 
the discrete distributions $f_i$ at all lattice nodes 
as follows: 
\begin{align}
  f_i(x,t=0) = w_i \mathrm{u}_0(x)
\end{align}
Since $w_i > 0$ and $\mathrm{u}_0(x) \geq 0$, we have 
$f_i(x,0) \geq 0$. Furthermore, we assume that the 
Dirichlet boundary conditions will be discretized 
using the \emph{weighted splitting method} (see 
equation \eqref{Eqn:Dirichlet_BC_ad_hoc}). 
The weighted splitting method guarantees the 
non-negativity of distributions $f_i$ for a 
lattice node on the boundary provided that 
the prescribed concentration on $\Gamma^{\mathrm{D}}$ 
is non-negative. That is,
\begin{align}
  \text{if} \; \mathrm{u}^{\mathrm{p}}(x \in \Gamma^{\mathrm{D}}, t) 
  \geq 0 \;\text{then} \; f_i( x \in \Gamma^{\mathrm{D}},t) \geq 0
\end{align}
So far, we made sure that all the distributions at the 
previous time-level are non-negative, and the discretization
of the boundary conditions will not disrupt the non-negativity
of distributions. Since the distributions at time $t$ are
non-negative, equilibrium distributions $f_i^{\mathrm{eq}}$
will be non-negative in the calculation of the collision 
step at time $t + \Delta t$. That is,
\begin{align}
	f_i \left( x, t\right) \geq 0 \Rightarrow 
	u\left(\boldsymbol{x}, t\right) = 
	\sum_{i} f_i\left( x,t\right) \geq 0 \Rightarrow
	 f_i^{\mathrm{eq}}\left( \boldsymbol{x}, t\right) = 
		w_i u\left(\boldsymbol{x},t\right) \geq 0 
\end{align}
If all of these conditions are satisfied, restricting 
the value of the relaxation time $\tau$ in equation 
\eqref{Eqn:Collision} can lead to non-negative 
concentrations at all lattice nodes and for all 
time-levels. We require that
\begin{align}
	1 - 1/\tau \geq 0
\end{align}
Using equation \eqref{Eqn:relaxation_SRT} 
and the above inequality, one can obtain 
the following condition that ensures 
non-negativity:
\begin{align}
  \label{Eqn:LBM_NN_inequality}
  \Delta t \geq \frac{\Delta x^2}{6\mathrm{D}}
\end{align}
That is, if the discretization parameters, 
$\Delta t$ and $\Delta x$, satisfy inequality 
\eqref{Eqn:LBM_NN_inequality} then all the 
distributions $f_i$ will be non-negative. 
Non-negativity of all $f_i$'s implies the 
non-negativity of the concentration 
$u(\boldsymbol{x},t)$, which stems 
directly from equation 
\eqref{Eqn:Distribution_Concnetration}.

Note that this result is only valid for one-dimensional 
pure diffusion equation (i.e., the advection velocity 
is zero) and for $D1Q3$ lattice model. Furthermore, the 
above condition does not guarantee the preservation of 
the comparison principles.
Deriving similar conditions for more sophisticated 
lattice models in two and three dimensions and for 
multiple-relaxation-time methods will require a more 
rigorous analysis. 

Another noteworthy point is that we have put stronger 
conditions on the values of distributions $f_i$ in order 
to meet the non-negative constraint for the concentration. 
In other words, for nodal concentrations to be non-negative 
we made sure that all the distributions are non-negative 
(i.e., $f_i(\boldsymbol{x},t) \geq 0 \; \forall i$). However,
this condition can be relaxed by allowing some of the $f_i$ 
to be negative, but with an additional constraint that 
$\sum_{i = 0}^{m-1} f_i( \boldsymbol{x}, t) \geq 0$. We do 
not pursue such an approach here, but one can consider 
them in future developments of lattice Boltzmann methods.

\section{DISCUSSION AND CONCLUDING REMARKS}
\label{Sec:LBM_S6_CR}
The maximum and comparison principles are two important 
mathematical properties of diffusion-type equations. The 
non-negative constraint is an important physical constraint 
on the concentration in transport and reactive-transport 
equations. There are other properties that the solutions 
to diffusion-type equations satisfy under appropriate 
conditions on the input data; for example, the decay 
property. 
A main challenge in designing a predictive numerical 
formulation is to satisfy these mathematical principles 
and the physical constraints in the discrete setting.
In this paper, using representative numerical examples, 
we have systematically documented that the current 
LBM-based formulations do not satisfy the maximum 
principle, the comparison principle, the non-negative 
constraint, and the decay property.
We have also shown that the discretization of boundary 
conditions has an affect on the performance of the 
lattice Boltzmann method in meeting these properties. 
To this end, we proposed a new way of discretizing 
Dirichlet boundary conditions -- the weighted splitting 
method. We then derived a theoretical bound in terms 
of the time-step and lattice cell size that guarantees 
non-negative values for the concentration under the 
weighted splitting method for one-dimension problems.

For anisotropic diffusion problems, we considered 
two representative multiple-relaxation-time lattice Boltzmann 
methods proposed in \citep{Yoshida_Nagaoka_JCP_2010} and 
\citep{Huang_Wu_JCP_2014}. 
We have shown that these multiple-relaxation-time methods 
can give unphysical negative values for the concentration. 
It needs to be emphasized that stability conditions for 
the lattice Boltzmann method (i.e., Courant-Fredrichs-Lewy 
conditions) have satisfied in all our numerical experiments. 
\emph{This implies that meeting stability conditions 
 alone does not guarantee the preservation of the 
 mentioned mathematical principles in the discrete 
 setting.}
The main findings of the paper about LBM-based 
formulations can be summarized as follows:
\begin{enumerate}[(a)]
\item \textsf{One-dimensional problems:} 
  For a given time-step, one can eliminate 
  the violation of the non-negative constraint 
  and the maximum principle by refining the 
  lattice cell size. For a given lattice cell 
  size, the violation of the non-negative 
  constraint and the maximum principle 
  cannot be eliminated by decreasing the 
  time-step. Both these trends are similar 
  to the finite element method,
  which have been reported in 
  \citep{Nakshatrala_Nagarajan_Shabouei_Arxiv_2013}.
\item \textsf{Critical time-step:} Based on a 
  simple theoretical analysis, we obtained the 
  following bound for the time-step and 
  lattice cell size to meet the 
  non-negative constraint under LBM for 1D problems: 
  \begin{align}
    \Delta t \geq \frac{\Delta x^2}{6D}
  \end{align}
  One can obtain exactly the same bound under the 
  single-field Galerkin finite element method based 
  on the backward Euler time-stepping scheme for 1D 
  problems \citep{Nakshatrala_Nagarajan_Shabouei_Arxiv_2013}. 
  This is an interesting result given the fact
  that the underlying basis of the lattice Boltzmann method (which solves 
  Boltzmann equation to obtain distributions at lattice 
  nodes) is completely different from that of the finite 
  element method (which is based on a weak formulation). 
\item \textsf{Isotropic vs. anisotropic diffusion:} 
  The violations of the non-negative constraint and 
  the maximum principle are smaller in magnitude 
  and smaller in terms of spatial extent when the 
  diffusion is isotropic. Also, for a given time-step, 
  one can decrease 
  these violations by refining the lattice cell size 
  in the case of isotropic diffusion. On the other 
  hand, neither decreasing the time-step nor refining 
  the lattice cell size will eliminate the violation 
  of the non-negative constraint for anisotropic 
  diffusion. 
\item \textsf{Convex vs. non-convex domains:}
  The magnitudes of the violation of the non-negative 
  constraint are larger for non-convex domains. However, 
  it needs to be emphasized that one may 
  have violations even on convex domains. 
\item \textsf{Comparison principle:} The lattice 
  Boltzmann method violates the comparison 
  principle even in a simple setting like 
  1D problems or isotropic diffusion. The 
  violation of the comparison principle 
  will be more predominant in the case 
  of anisotropic diffusion. 
\item \textsf{Decay property:} The LBM-based 
  formulations, in general, violate the decay property. 
\item The lattice Boltzmann method does not posses 
  a variational structure, similar to the one proposed 
  by the finite element method. Due to this reason, the 
  non-negative formulations proposed under the finite 
  element method 
  (e.g., \citep{Nakshatrala_Nagarajan_Shabouei_Arxiv_2013}) 
  cannot be directly extended to the lattice Boltzmann 
  method.  
\item The only procedure that is available to meet 
  the non-negative constraint under the lattice 
  Boltzmann method is the clipping procedure, 
  which basically chops off the negative values. 
  But, this procedure fixes neither the violation 
  of the decay property nor the violation of the 
  comparison principle. Moreover, this method does 
  not have any physical or mathematical basis, and 
  it is rather \emph{ad hoc}.
\end{enumerate}

One should be wary of violations of the non-negative 
constraint, and the maximum and comparison principles 
in the numerical simulations using LBM. 
In the case of \emph{isotropic} diffusion, the authors 
suggest investigating the occurrence of the mentioned 
violations, if any of these violations occur, they can 
be significantly reduced by refining the lattice cell 
size and the time-step in accordance with the CFL 
condition. However, in the case of \emph{anisotropic} 
diffusion, no clear-cut guideline for reducing the 
violations exists. As demonstrated earlier, refining 
the discretization parameters may not improve the 
numerical solution. 
A future research direction could be development of
LBM-based formulations for transient diffusion-type 
equations (i.e., diffusion, advection-diffusion and 
advection-diffusion-reaction equations) that respect 
the maximum and comparison principles, and meet the 
non-negative constraint. This paper can serve as a 
source of benchmark problems for such a research 
endeavor.

\section*{ACKNOWLEDGMENTS}
The authors acknowledge the support from the Department of 
Energy through Subsurface Biogeochemical Research Program. 
Neither the United States Government nor any agency thereof, 
nor any of their employees, makes any warranty, express or 
implied, or assumes any legal liability or responsibility 
for the accuracy, completeness, or usefulness of any 
information. The opinions expressed in this paper are 
those of the authors and do not necessarily reflect 
that of the sponsor(s). 

\bibliographystyle{plainnat}
\bibliography{References/Master_References,References/Books}


\begin{figure}
  \centering
  \psfrag{time}{time}
  \psfrag{max c}{$u_{\max}\left( t \right)$}
  \psfrag{c1}{$\Delta t = 10^{-3}$}
  \psfrag{c2}{$\Delta t = 10^{-4}$}
  \psfrag{c3}{$\Delta t = 10^{-5}$}
  \includegraphics[clip,scale=0.4]{Figures/1D_1/max_c_dt.eps}
  \caption{\textsf{One-dimensional problem with 
      uniform initial condition:}~The maximum nodal 
    concentration is plotted against time for various 
    time-steps. In all the cases, the lattice cell size 
    is taken as $\Delta x = 0.1$. For the prescribed data, 
    the maximum principle asserts that the concentration 
    should be between $0$ and $1$ in the entire domain 
    and for all times. \emph{This figure shows that 
      the single-relaxation-time method violates the 
      maximum principle. Moreover, reducing the time-step 
      keeping the lattice cell size fixed may not eliminate 
      the violation of the maximum principle.} 
    \label{Fig:1D_1_Max_c_dt}}
\end{figure}

\end{document}